\documentclass[12pt,thmsb,titlepage]{article}
\usepackage{amsfonts, amsthm}
\usepackage{amssymb}
\usepackage{amsmath}
\usepackage{natbib}
\usepackage{graphicx}
\usepackage{sw20elba}
\usepackage{caption}
\usepackage{subcaption}

\renewcommand{\theequation}{\thesection.\arabic{equation}}
\input{tcilatex}

\begin{document}

\author{Humberto Moreira and Marcelo J. Moreira \thanks{%
This paper expands upon and supersedes the corresponding sections of our
working paper \textquotedblleft Contributions to the Theory of Optimal
Tests.\textquotedblright\ We thank Benjamin Mills and Gustavo de Castro for
outstanding research assistance, and are particularly indebted to Gustavo
for suggesting and implementing the designs for power comparisons reported
here. We thank Leandro Gorno, Patrik Guggenberger, Alexei Onatski, and Lucas
Vilela for helpful comments; Jose Diogo Barbosa, Felipe Flores, and Leonardo
Salim for suggestions on earlier drafts of this paper; and Isaiah Andrews
and Jose Olea for productive discussions and sharing numerical codes. We
gratefully acknowledge the research support of CNPq, FAPERJ, and NSF (via
grant SES-0819761).} \\
%EndAName
\emph{FGV/EPGE}}
\title{\textbf{Optimal Two-Sided Tests for Instrumental Variables Regression
with Heteroskedastic and Autocorrelated Errors}}
\date{This version: \today }
\maketitle

\begin{abstract}
\hspace{0.25in}This paper considers two-sided tests for the parameter of an
endogenous variable in an instrumental variable (IV) model with
heteroskedastic and autocorrelated errors. We develop the finite-sample
theory of weighted-average power (WAP) tests with normal errors and a known
long-run variance. We introduce two weights which are invariant to
orthogonal transformations of the instruments; e.g., changing the order in
which the instruments appear. While tests using the MM1 weight can be
severely biased, optimal tests based on the MM2 weight are naturally
two-sided when errors are homoskedastic.

We propose two boundary conditions that yield two-sided tests whether errors
are homoskedastic or not. The locally unbiased (LU)\ condition is related to
the power around the null hypothesis and is a weaker requirement than
unbiasedness. The strongly unbiased (SU) condition is more restrictive than
LU, but the associated WAP tests are easier to implement. Several tests are
SU in finite samples or asymptotically, including tests robust to weak IV
(such as the Anderson-Rubin, score, conditional quasi-likelihood ratio, and
I. \citeauthor{Andrews15}' (\citeyear{Andrews15}) PI-CLC tests) and
two-sided tests which are optimal when the sample size is large and
instruments are strong.

We refer to the WAP-SU tests based on our weights as MM1-SU and MM2-SU
tests. Dropping the restrictive assumptions of normality and known variance,
the theory is shown to remain valid at the cost of asymptotic
approximations. The MM2-SU test is optimal under the strong IV asymptotics,
and outperforms other existing tests under the weak IV asymptotics.
\end{abstract}

%\thanks{Humberto Moreira acknowledges financial support from CNPq and Marcelo J.
%Moreira gratefully acknowledges research support of the National Science
%Foundation via grant number SES-0819761. The authors thank Alexei Onatski
%and Jack Porter for comments. The first version of this paper entitled
%\textquotedblleft Contributions to the Theory Optimal
%Tests.\textquotedblright}

\newpage

\section{Introduction}

\setcounter{equation}{0}\hspace{0.25in}In an instrumental variable (IV)
model, researchers often rely on asymptotic approximations when making
inference on the structural coefficients. These approximations, however, can
be poor when instruments are weakly correlated with the endogenous
regressors as explained by \citet{NelsonStartz90}, \citet{BoundJaegerBaker95}%
, \citet{Dufour97}, and \citet{StaigerStock97}. The goal is to find reliable
econometric methods regardless of how strong the instruments are.

There has been some progress in the IV model with one endogenous variable
and $k$ instruments when errors are homoskedastic. \citet{AndersonRubin49}
propose a test statistic which has an asymptotic chi-square-$k$ distribution
regardless of how weak the instruments are. \citeauthor{Moreira01} (%
\citeyear{Moreira01}, \citeyear{Moreira09a}) shows that the Anderson-Rubin
statistic is optimal in the just-identified model, but points out potential
power gains when there exists more than one instrument. \citet{Kleibergen02}
and \citet{Moreira02} show that a score (LM) test statistic has a standard
chi-square-one distribution whether the instruments are weak or not. %
\citet{Moreira03} proposes to replace the critical value number by
conditional quantiles of test statistics. These conditional tests are
similar by construction, hence have correct size. He applies the conditional
method to the likelihood ratio (LR) statistic and the two-sided Wald
statistic. \citet{AndrewsMoreiraStock06} (hereinafter, AMS06) show that the
conditional likelihood ratio (CLR) test satisfies natural orthogonal
invariance conditions and is nearly optimal. \citet{AndrewsMoreiraStock07}
find that conditional Wald (CW) tests, however, have poor behavior and
object to their use in empirical work. \citet{MillsMoreiraVilela14} show
that the bad performance of CW tests is due to the asymmetric distribution
of one-sided Wald statistics when instruments are weak. By extending %
\citeauthor{Moreira03}'s (\citeyear{Moreira03}) conditional approach, they
find approximately unbiased Wald tests whose power is comparable to the CLR
test.

While use of the IV model with homoskedastic errors was important to advance
the literature on weak identification, the IV model with heteroskedastic and
autocorrelated (HAC) errors is considerably more relevant for applied
researchers. Some of the theoretical findings for homoskedastic errors are
easily extended for more complicated stochastic processes, whereas others
are not. Important work by \citet{StockWright00}, \citet{GuggenbergerSmith05}%
, \citet{Kleibergen06}, \citet{Otsu06}, and \citet{AndrewsMikusheva15},
among others, extends the tests conceived for the simple homoskedastic IV
model to the generalized method of moments (GMM) and generalized empirical
likelihood (GEL) frameworks. Their tests are of course applicable to the
HAC-IV model, but it is unknown whether these adaptations are optimal. The
purpose of this paper is exactly this: to develop a theory of optimal
two-sided tests for the HAC-IV model.

We are able to find a statistic that is pivotal and independent of a second
statistic, which is sufficient and complete for the instruments'
coefficients under the null. We show that the invariance argument of AMS06
for homoskedastic errors is only applicable if a (long-run) variance has a
Kronecker product structure. This limitation has profound consequences for
the behavior of weighted-average power (WAP) tests. We choose two priors for
the structural parameter and the instruments' coefficients and denote the
associated test statistics MM1 and MM2. The priors are chosen to illustrate
the effect of a poor weight choice on the power of WAP\ tests. Although
priors vanish asymptotically as in the Bernstein-von Mises theorem, the
associated tests can behave quite differently in finite samples (or under
the weak-instrument asymptotics). When a variance matrix has a Kronecker
product structure, both test statistics are orthogonally invariant, but only
MM2 satisfies an additional \emph{sign} invariance argument that preserves
the two-sided hypothesis testing problem. As a consequence, a WAP similar
test based on the MM1 statistic can behave as a one-sided test and have poor
power even with homoskedastic errors (this problem is analogous to the
conditional Wald tests documented by \citet{AndrewsMoreiraStock07}) while
the WAP similar test using the MM2 statistic has overall good power with a
Kronecker-product variance matrix. Other weight choices face the same
difficulties as the MM1 statistic for the HAC-IV model, including the
recently proposed WAP similar test by \citet{Olea15}, denoted ECS (HAC-IV).

When the (long-run) variance matrix does not have a Kronecker product
representation and the model is identified, the Anderson-Rubin test (among
other equivalent tests) is the uniformly most powerful unbiased test. In the
over-identified model, we show theoretically that it is possible to find a
weight so that the test is approximately unbiased and admissible. The lack
of invariance, however, makes it harder to construct such weights. In
practice, we endogeneize this search by imposing in the WAP maximization
problem a boundary condition based on the local power around the null
hypothesis. This locally unbiased (LU) condition is a weaker requirement
than unbiasedness, so it does not rule out admissibility. The WAP-LU tests
are found with non-linear algorithms, which makes it difficult to implement
them. We then propose a stronger requirement than LU, denoted the strongly
unbiased (SU) condition. The resulting class of tests includes several
two-sided tests robust to weak IV, including the Anderson-Rubin, score,
(pseudo) likelihood ratio tests by \citet{Kleibergen06} and %
\citet{AndrewsGuggenberger14b}, and I. \citeauthor{Andrews15}' (%
\citeyear{Andrews15}) PI-CLC tests. Two-sided optimal tests also satisfy the
SU condition asymptotically when the sample size is large and instruments
are strong. The WAP-SU tests have power close to the WAP-LU\ tests based on
the MM1 and MM2 weights, with the advantage being that the WAP-SU tests are
easy to implement with a standard linear programming software package. We
refer to the WAP-SU tests based on our weights as MM1-SU and MM2-SU tests.

We follow I. \citet{Andrews15} and implement numerical simulations based on %
\citet{Yogo04}. We choose, however, \citeauthor{Yogo04}'s (\citeyear{Yogo04}%
) design where the endogenous variable is the real stock return and the
instruments are genuinely weak. We find that, as our theory predicts, the
WAP similar tests can be quite erratic. In some designs, they behave as
usual two-sided tests and have good power. In other designs they behave as
one-sided tests and have power near zero. We do not recommend the MM1 and
MM2 similar tests for empirical researchers. The MM2-SU test, however,
outperforms other tests (including the MM1-SU test) and when it occasionally
has less power than competing tests, the power loss is small. We recommend
the use of the MM2-SU\ test in empirical work. Our asymptotic analysis is
quite general and encompasses all WAP similar and WAP-SU tests whose weight
does not depend strongly on the sample size.

The remainder of this paper is organized as follows. Section \ref{IV Model
and Statistics} introduces the HAC-IV model and presents the test
statistics, including the MM1 and MM2 statistics. Sections \ref{Maximization
Sec} and \ref{Conditions Sec} discuss the power maximization problem and the
WAP-LU\ and WAP-SU\ tests. Section \ref{Numerical Sec} presents power curves
and the role of LU\ and SU conditions in obtaining WAP tests with overall
good power. Section \ref{Asymptotic Sec} develops an asymptotic framework
that encompasses the weak IV and strong IV asymptotics. Section \ref%
{Application Sec} revisits the work of I. \citet{Andrews15} and %
\citet{Yogo04} on testing the intertemporal rate of substitution, with one
important modification. Section \ref{Conclusion Sec} contains concluding
remarks. All proofs are given in the appendices.

\section{The IV Model and Statistics \label{IV Model and Statistics}}

Consider the instrumental variable model 
\begin{eqnarray*}
y_{1} &=&y_{2}\beta +u \\
y_{2} &=&Z\pi +v_{2},
\end{eqnarray*}%
where $y_{1}$ and $y_{2}$ are $n$ $\times $ $1$ vectors of observations on
two endogenous variables, $Z$ is an $n\times k$ matrix of nonrandom
exogenous variables having full column rank, and $u$ and $v_{2}$ are $n$ $%
\times $ $1$ unobserved disturbance vectors having mean zero. The goal here
is to test the null hypothesis $H_{0}:\beta =\beta _{0}$ against the
alternative hypothesis $H_{1}:\beta \neq \beta _{0}$, treating $\pi $ as a
nuisance parameter. We do not not include covariates in this model, but we
note that can be easily handled by the usual projection arguments; see AMS06.

We look at the reduced-form model for $Y=\left[ y_{1},y_{2}\right] $:%
\begin{equation}
Y=Z\pi a^{\prime }+V,  \label{(reduced-form IV)}
\end{equation}%
where $a=\left( \beta ,1\right) ^{\prime }$ and $V=\left[ v_{1},v_{2}\right]
=\left[ u+v_{2}\beta ,v_{2}\right] $ is the $n\times 2$ matrix of
reduced-form errors. We allow the errors to be heteroskedastic and
autocorrelated. Let $P_{1}=Z\left( Z^{\prime }Z\right) ^{-1/2}$ and let $%
\left[ P_{1},P_{2}\right] \in \mathcal{O}_{n}$, the group of $n\times n$
orthogonal matrices. Pre-multiplying the reduced-form model (\ref%
{(reduced-form IV)}) by $\left[ P_{1},P_{2}\right] ^{\prime }$, we obtain
the pair of statistics $P_{1}^{\prime }Y$ and $P_{2}^{\prime }Y$. In this
section, we assume that $\left( Z^{\prime }Z\right) ^{-1/2}Z^{\prime }V$ is
normally distributed with known variance matrix $\Sigma $ (this assumption
will be relaxed later at the cost of asymptotic approximations). The
statistic $P_{2}^{\prime }Y$ is ancillary and we do not have previous
knowledge about the correlation structure on $V$. In consequence, we
consider tests based on $R=P_{1}^{\prime }Y$:%
\begin{equation*}
R=\mu a^{\prime }+\left( Z^{\prime }Z\right) ^{-1/2}Z^{\prime }V,
\end{equation*}%
where $\mu =\left( Z^{\prime }Z\right) ^{1/2}\pi $.

It is convenient to find the one-to-one transformation of $R$ given by the
pair%
\begin{eqnarray}
S &=&\left[ \left( b_{0}^{\prime }\otimes I_{k}\right) \Sigma \left(
b_{0}\otimes I_{k}\right) \right] ^{-1/2}\left( b_{0}^{\prime }\otimes
I_{k}\right) \overline{R}\text{ and}  \label{(Defns of S and T)} \\
T &=&\left[ \left( a_{0}^{\prime }\otimes I_{k}\right) \Sigma ^{-1}\left(
a_{0}\otimes I_{k}\right) \right] ^{-1/2}\left( a_{0}^{\prime }\otimes
I_{k}\right) \Sigma ^{-1}\overline{R},  \notag
\end{eqnarray}%
where $\overline{R}=vec\left[ \left( Z^{\prime }Z\right) ^{-1/2}Z^{\prime }Y%
\right] $, $a_{0}=\left( \beta _{0},1\right) ^{\prime }$ and $b_{0}=\left(
1,-\beta _{0}\right) ^{\prime }$. The pair $S$ and $T$ have three important
properties: \emph{(i)} they are independent; \emph{(ii)} $S$ is pivotal; and 
\emph{(iii)} $T$ is complete and sufficient for $\mu $ under the null. More
specifically, the statistics $S$ and $T$ have distribution%
\begin{eqnarray}
S &\sim &N\left( \left( \beta -\beta _{0}\right) C_{\beta _{0}}\mu
,I_{k}\right) \text{ and }T\sim N\left( D_{\beta }\mu ,I_{k}\right) \text{,
where}  \label{(Dist S and T)} \\
C_{\beta _{0}} &=&\left[ \left( b_{0}^{\prime }\otimes I_{k}\right) \Sigma
\left( b_{0}\otimes I_{k}\right) \right] ^{-1/2}\text{ and}  \notag \\
D_{\beta } &=&\left[ \left( a_{0}^{\prime }\otimes I_{k}\right) \Sigma
^{-1}\left( a_{0}\otimes I_{k}\right) \right] ^{-1/2}\left( a_{0}^{\prime
}\otimes I_{k}\right) \Sigma ^{-1}\left( a\otimes I_{k}\right) .  \notag
\end{eqnarray}%
The joint density $f_{\beta ,\mu }\left( s,t\right) $ is given by%
\begin{eqnarray*}
f_{\beta ,\mu }\left( s,t\right) &=&\left( 2pi\right) ^{-k/2}\exp \left( -%
\frac{\left\Vert s-\left( \beta -\beta _{0}\right) C_{\beta _{0}}\mu
\right\Vert ^{2}}{2}\right) \times \left( 2pi\right) ^{-k/2}\exp \left( -%
\frac{\left\Vert t-D_{\beta }\mu \right\Vert ^{2}}{2}\right) \\
&=&f_{\beta ,\mu }^{S}\left( s\right) \times f_{\beta ,\mu }^{T}\left(
t\right) ,
\end{eqnarray*}%
where $pi=3.1415...$ and $f_{\beta ,\mu }^{S}\left( s\right) $ and $f_{\beta
,\mu }^{T}\left( t\right) $ are the marginal densities for $S$ and $T$.

Examples of test statistics based on $S$ and $T$ are the Anderson-Rubin
(AR), the score or Lagrange multiplier (LM), and the quasi likelihood ratio
(LR) statistics. Anderson and Rubin (1949) propose to use a pivotal
statistic. In our model the Anderson-Rubin statistic is given by 
\begin{equation}
AR=S^{\prime }S.  \label{(AR stat)}
\end{equation}%
In Appendix A, we derive the $LM$ and $LR$ statistics under that the
assumption the errors are normal. For any full column rank matrix $X$, let $%
N_{X}=X\left( X^{\prime }X\right) ^{-1}X^{\prime }$ and $M_{X}=I-N_{X}$.
Then the $LM$ statistic simplifies to%
\begin{equation}
LM=S^{\prime }N_{C_{\beta _{0}}D_{\beta _{0}}^{-1}T}S\text{.}
\label{(LM stat)}
\end{equation}%
The likelihood ratio statistic is given by%
\begin{equation}
LR=\max_{a}\overline{R}^{\prime }\Sigma ^{-1/2}N_{\Sigma ^{-1/2}(a\otimes
I_{k})}\Sigma ^{-1/2}\overline{R}-T^{\prime }T.  \label{(LR stat)}
\end{equation}%
The $LR$ statistic is apparently not a simple function of $S$ and $T$ (which
makes it difficult to implement the test coupled with conditional critical
values). \citet{Kleibergen06} instead adapts the formula for the likelihood
ratio statistic derived by \citet{Moreira03} in the homoskedastic IV model
to the GMM framework. For the HAC-IV model, this quasi likelihood ratio
statistic becomes%
\begin{equation}
QLR=\frac{AR-r\left( T\right) +\sqrt{\left( AR-r\left( T\right) \right)
^{2}+4LM\cdot r\left( T\right) }}{2},  \label{(QLR stat)}
\end{equation}%
where $AR$ and $LM$ are defined in (\ref{(AR stat)}) and (\ref{(LM stat)}),
and $r\left( T\right) =T^{\prime }T$. \citet{AndrewsGuggenberger14b} use a
Kronecker product $\Omega \otimes \Phi $ (where $\Omega $ and $\Phi $ are
positive-definite matrices respectively with dimensions $2\times 2$ and $%
k\times k$) approximation to the variance $\Sigma $; see %
\citet{VanLoanPtsianis93} for more details on Kronecker product
approximations.

We now present two novel WAP\ statistics based on the weighted-average
density%
\begin{equation}
h_{\Lambda }\left( s,t\right) =\int f_{\beta ,\mu }\left( s,t\right) \text{ }%
d\Lambda \left( \beta ,\mu \right) .  \label{(WAP density)}
\end{equation}%
These weight functions use the Kronecker product $\Omega \otimes \Phi $
approximation to $\Sigma $ with the Frobenius norm (i.e.,\ the norm of a
matrix $X$ is given by $\left\Vert X\right\Vert =\sqrt{tr\left( X^{\prime
}X\right) }$). For the MM1 statistic $h_{1}\left( s,t\right) $, we choose $%
\Lambda \left( \beta ,\mu \right) $ to be $N\left( \beta _{0},1\right)
\times N\left( 0,\sigma ^{2}\Phi \right) $. For the MM2 statistic $%
h_{2}\left( s,t\right) $, we first define the identity $\tan \left( \theta
\right) \equiv d_{\beta \left( \theta \right) }/c_{\beta \left( \theta
\right) }$, where%
\begin{equation}
c_{\beta }=(\beta -\beta _{0})\cdot (b_{0}^{\prime }\Omega b_{0})^{-1/2}%
\text{ and }d_{\beta }=a^{\prime }\Omega ^{-1}a_{0}\cdot (a_{0}^{\prime
}\Omega ^{-1}a_{0})^{-1/2}.  \label{(c_b and d_b)}
\end{equation}%
We choose $\Lambda \left( \beta ,\mu \right) $ so that the prior for $\theta 
$ and $\mu $ are \textit{Unif}$\left[ -pi,pi\right] \times N\left(
0,\left\Vert l_{\beta \left( \theta \right) }\right\Vert ^{-2}\zeta \cdot
\Phi \right) $, where $l_{\beta }=\left( c_{\beta },d_{\beta }\right)
^{\prime }$.

In Appendix A, we show that the MM1 and MM2 statistics are 
\begin{eqnarray}
h_{1}\left( s,t\right) &\hspace{-0.08in}=\hspace{-0.08in}&\left( 2pi\right)
^{-k-1/2}\int \left\vert \Psi _{\beta ,\sigma ^{2}}\right\vert ^{-1/2}\exp
\left( -\frac{\left( s^{\prime },t^{\prime }\right) \Psi _{\beta ,\sigma
^{2}}^{-1}\left( s^{\prime },t^{\prime }\right) ^{\prime }+\left( \beta
-\beta _{0}\right) ^{2}}{2}\right) d\beta  \label{(h densities)} \\
h_{2}\left( s,t\right) &\hspace{-0.08in}=\hspace{-0.08in}&\left( 2pi\right)
^{-\left( k+1\right) }\int_{-pi}^{pi}\left\vert \Psi _{\beta \left( \theta
\right) ,\left\Vert l_{\beta \left( \theta \right) }\right\Vert ^{-2}\zeta
}\right\vert ^{-1/2}\exp \left( -\frac{\left( s^{\prime },t^{\prime }\right)
\Psi _{\beta \left( \theta \right) ,\left\Vert l_{\beta \left( \theta
\right) }\right\Vert ^{-2}\zeta }^{-1}\left( s^{\prime },t^{\prime }\right)
^{\prime }}{2}\right) d\theta ,  \notag
\end{eqnarray}%
where the matrix $\Psi _{\beta ,\sigma ^{2}}$ is given by%
\begin{equation}
\Psi _{\beta ,\sigma ^{2}}=I_{2}\otimes I_{k}+\sigma ^{2}\left[ 
\begin{array}{cc}
\left( \beta -\beta _{0}\right) ^{2}C_{\beta _{0}}\Phi C_{\beta _{0}} & 
\left( \beta -\beta _{0}\right) C_{\beta _{0}}\Phi D_{\beta }^{\prime } \\ 
\left( \beta -\beta _{0}\right) D_{\beta }\Phi C_{\beta _{0}} & D_{\beta
}\Phi D_{\beta }^{\prime }%
\end{array}%
\right] .  \label{(Psi_b,sigma2)}
\end{equation}

\subsection{Kronecker Variance Matrix}

We consider here the special case where $\Sigma =\Omega \otimes \Phi $
exactly. This framework is particularly interesting for two reasons. First,
it encompasses the homoskedastic case by taking $\Phi $ to be the identity
matrix. We will show that the $S$ and $T$ statistics for general error
structure simplify to the original statistics of \citeauthor{Moreira01} (%
\citeyear{Moreira01}, \citeyear{Moreira09a}) for the homoskedastic model.
Second, the model where $\Sigma $ has a Kronecker product structure enjoys
natural invariance properties. Some statistics are invariant but others are
not. This has profound consequences for testing procedures based on these
statistics. Indeed, typical tests based on noninvariant statistics (such as
those using a constant or \citeauthor{Moreira03}'s (\citeyear{Moreira03})
conditional critical value function) behave as one-sided tests for parts of
the parameter space. We will illustrate this problem numerically in Section %
\ref{Numerical Sec}.

When $\Sigma =\Omega \otimes \Phi $, the statistics $S$ and $T$ defined in (%
\ref{(Defns of S and T)}) simplify to 
\begin{eqnarray}
S &=&\Phi ^{-1/2}(Z^{\prime }Z)^{-1/2}Z^{\prime }Yb_{0}\cdot (b_{0}^{\prime
}\Omega b_{0})^{-1/2}\text{ and }  \label{(Defns of S and T kron)} \\
T &=&\Phi ^{-1/2}(Z^{\prime }Z)^{-1/2}Z^{\prime }Y\Omega ^{-1}a_{0}\cdot
(a_{0}^{\prime }\Omega ^{-1}a_{0})^{-1/2}.  \notag
\end{eqnarray}%
Their distribution is given by%
\begin{equation}
S\sim N\left( c_{\beta }\Phi ^{-1/2}\mu ,I_{k}\right) \text{ and }T\sim
N\left( d_{\beta }\Phi ^{-1/2}\mu ,I_{k}\right) .
\label{(Dist S and T kron)}
\end{equation}%
AMS06 use invariance arguments for the special case $\Phi =I_{k}$. However,
the parameter $\mu _{\Phi }=\Phi ^{-1/2}\mu $ is unknown because $\mu $ is
unknown. Hence, AMS06's invariance argument applies to the new parameter $%
\mu _{\Phi }=\Phi ^{-1/2}\mu $. Specifically, let $g\in \mathcal{O}_{n}$ and
consider the transformation in the sample space%
\begin{equation*}
g\circ \left( S,T\right) =\left( gS,gT\right) .
\end{equation*}%
The induced transformation in the parameter space is%
\begin{equation*}
g\circ \left( \beta ,\mu _{\Phi }\right) =\left( \beta ,g\mu _{\Phi }\right)
.
\end{equation*}

Invariant tests depend on the data only through 
\begin{equation}
Q=\left[ 
\begin{array}{cc}
Q_{S} & Q_{ST} \\ 
Q_{ST} & Q_{T}%
\end{array}%
\right] =\left[ 
\begin{array}{cc}
S^{\prime }S & S^{\prime }T \\ 
S^{\prime }T & T^{\prime }T%
\end{array}%
\right] .  \label{(Q def)}
\end{equation}%
The density of $Q$ at $q$ for the parameters $\beta $ and $\lambda =\pi
^{\prime }\left( Z^{\prime }Z\right) ^{1/2}\Phi ^{-1}\left( Z^{\prime
}Z\right) ^{1/2}\pi $ is given by 
\begin{eqnarray*}
&&f_{\beta ,\lambda }(q_{S},q_{ST},q_{T})\overset{}{=}K_{0}\exp (-\lambda
(c_{\beta }^{2}+d_{\beta }^{2})/2)\left\vert q\right\vert ^{(k-3)/2} \\
&&\hspace{1.03in}\times \exp (-(q_{S}+q_{T})/2)(\lambda \xi _{\beta
}(q))^{-(k-2)/4}I_{(k-2)/2}(\sqrt{\lambda \xi _{\beta }(q)}),
\end{eqnarray*}%
where $K_{0}^{-1}=2^{(k+2)/2}pi^{1/2}\Gamma _{(k-1)/2}$, $\Gamma _{(\cdot )}$
is the gamma function, $I_{(k-2)/2}(\cdot )$ denotes the modified Bessel
function of the first kind, and 
\begin{equation}
\xi _{\beta }(q)=c_{\beta }^{2}q_{S}+2c_{\beta }d_{\beta }q_{ST}+d_{\beta
}^{2}q_{T}.  \label{(Defn of Xi_Beta)}
\end{equation}

The following proposition shows that the WAP densities $h_{1}\left(
s,t\right) $ and $h_{2}\left( s,t\right) $ are invariant when the covariance
matrix is a Kronecker product. Indeed, the Kronecker product approximation $%
\Omega \otimes \Phi $ to $\Sigma $ in the definition of the weights was
chosen exactly to guarantee the test statistics are \emph{orthogonal}
invariant.

AMS06 show there also exists a \emph{sign} transformation that preserves the
two-sided hypothesis testing problem. Consider the group $\mathcal{O}_{1}$,
which contains only two elements: $\overline{g}\in $ $\left\{ -1,1\right\} $%
. The group transformation in the sample is 
\begin{equation*}
\overline{g}\circ \left( Q_{S},Q_{ST},Q_{T}\right) =\left( Q_{S},\overline{g}%
\cdot Q_{ST},Q_{T}\right) ,
\end{equation*}%
whose maximal invariant is $Q_{S}$, $\left\vert Q_{ST}\right\vert $, and $%
Q_{T}$. This group yields a\ transformation in the parameter space. For $%
\overline{g}=-1$, AMS06 show that this transformation is%
\begin{eqnarray}
\overline{g}\circ \left( \beta ,\lambda \right) &=&\left( \beta _{0}-\frac{%
d_{\beta _{0}}(\beta -\beta _{0})}{d_{\beta _{0}}+2j_{\beta _{0}}(\beta
-\beta _{0})},\lambda \frac{(d_{\beta _{0}}+2j_{\beta _{0}}(\beta -\beta
_{0}))^{2}}{d_{\beta _{0}}^{2}}\right) ,\text{ where}  \notag \\
j_{\beta _{0}}\hspace{-0.08in} &=&\hspace{-0.08in}\frac{e_{1}^{\prime
}\Omega ^{-1}a_{0}}{(a_{0}^{\prime }\Omega ^{-1}a_{0})^{-1/2}}\text{ and }%
e_{1}=(1,0)^{\prime }.  \label{(Defn of Beta2* and Lambda2*)}
\end{eqnarray}%
(by the definition of a group, the parameter remains unaltered at $\overline{%
g}=1$). The transformation in (\ref{(Defn of Beta2* and Lambda2*)}) flips
the sign of $\beta -\beta _{0}$ for $\beta \neq \beta _{AR}$ defined as 
\begin{equation}
\beta _{AR}=\frac{\omega _{11}-\omega _{12}\beta _{0}}{\omega _{12}-\omega
_{22}\beta _{0}}\text{ where }\Omega =\left[ \omega _{i,l}\right] \text{.}
\label{(Defn of beta_AR)}
\end{equation}%
So the \emph{sign} transformation preserves the two-sided hypothesis testing
problem $H_{0}:\beta =\beta _{0}$ against $H_{1}:\beta \neq \beta _{0}$, but
not the one-sided, e.g., testing $H_{0}:\beta \leq \beta _{0}$ against $%
H_{1}:\beta >\beta _{0}$.

\bigskip

\begin{proposition}
\label{Invariant WAP HAC-IV Prop} The following holds when $\Sigma =\Omega
\otimes \Phi $:\newline
\emph{(i)} The weighted-average densities $h_{1}\left( s,t\right) $ and $%
h_{2}\left( s,t\right) $ are invariant to orthogonal transformations. That
is, they depend on the data only through $Q$; and\newline
\emph{(ii)} The weighted-average density $h_{2}\left( s,t\right) $ is
invariant to sign transformations. It depends on the data only through $%
Q_{S} $, $\left\vert Q_{ST}\right\vert $, and $Q_{T}$.\newline
\end{proposition}

\bigskip

The MM1 statistic is not \emph{sign} invariant. We can create a
weighted-average statistic that is sign invariant by replacing the weight in 
$h_{1}=\int f_{\beta _{0},\lambda }\left( q_{S},q_{ST},q_{T}\right) $ $%
d\Lambda _{1}\left( \beta ,\lambda \right) $ by 
\begin{equation}
\Lambda \left( \beta ,\lambda \right) =\frac{\Lambda _{1}\left( \beta
,\lambda \right) +\Lambda _{1}\left( \overline{g}\circ \left( \beta ,\lambda
\right) \right) }{2},  \label{(correct weight)}
\end{equation}%
for $\overline{g}=-1$. We note that 
\begin{equation*}
\int f_{\beta ,\lambda }(q_{S},q_{ST},q_{T})\text{ }d\Lambda \left( \beta
,\lambda \right) =\int \int f_{\beta ,\lambda }(q_{S},q_{ST},q_{T})\text{ }%
d\Lambda _{1}\left( \overline{g}\circ \left( \beta ,\lambda \right) \right) 
\text{ }\nu \left( d\overline{g}\right) ,
\end{equation*}%
where $\nu $ is the Haar probability measure on the group $\mathcal{O}_{1}$: 
$\nu \left( \left\{ 1\right\} \right) =\nu \left( \left\{ -1\right\} \right)
=1/2$. Because 
\begin{eqnarray*}
\int f_{\beta ,\lambda }(q_{S},-q_{ST},q_{T})\ d\Lambda \left( \beta
,\lambda \right) &=&\int f_{\left( -1\right) \circ \left( \beta ,\lambda
\right) }(q_{S},q_{ST},q_{T})\text{ }d\Lambda \left( \beta ,\lambda \right)
\\
&=&\int f_{\beta ,\lambda }(q_{S},q_{ST},q_{T})\text{ }d\Lambda \left( \beta
,\lambda \right) ,
\end{eqnarray*}%
the weighted-average statistic based on (\ref{(correct weight)}) only
depends on $q_{S},\left\vert q_{ST}\right\vert ,q_{T}$. But the MM2
statistic is already sign invariant for having chosen a clever prior for $%
\beta $ and $\mu $. In fact, the MM2 prior was chosen so that the final
statistic is sign invariant. Tests based on $h_{2}\left( s,t\right) $ are
naturally two-sided tests for the null $H_{0}:\beta =\beta _{0}$ against the
alternative $H_{1}:\beta \neq \beta _{0}$ when $\Sigma =\Omega \otimes \Phi $%
. This important property does not hold for standard tests based on $%
h_{1}\left( s,t\right) $. The WAP test (denoted ECS-HACIV) proposed recently
by \citet{Olea15} is not \emph{sign} invariant either. Sections \ref%
{Numerical Sec} and \ref{Application Sec} present numerical simulations
showing that all these WAP similar tests can behave like one-sided tests for
some parameter values. In the next section, we will discuss ways to
circumvent this problem whether $\Sigma $ has a Kronecker product structure
or not.

\section{Weighted-Average Power Tests \label{Maximization Sec}}

So far, we have only described test statistics. Coupled with critical
values, we obtain the test procedures commonly used in the literature. The
Anderson-Rubin test rejects the null when $AR>c\left( k\right) $, where $%
c\left( d\right) $ is the $1-\alpha $ quantile of a chi-square distribution
with $d$ degrees of freedom. The LM test rejects the null when $LM>c\left(
1\right) $. The conditional tests reject the null when each test statistic $%
\psi \left( S,T\right) >\kappa \left( T\right) $. Each critical value
function $\kappa \left( T\right) $ is the null conditional quantile of $\psi 
$ given $T=t$; see \citet{Moreira03} for details (we omit the dependence of
the critical value function on the statistic $\psi $ when there is no
ambiguity). For example, the CQLR test rejects the null when the QLR
statistic defined in (\ref{(QLR stat)}) is larger than the conditional
critical value.

Our goal in this section is to find optimal tests. Specifically, a test is
defined to be a measurable function $\phi \left( s,t\right) $ that is
bounded by $0$ and $1$. For a given outcome, the test rejects the null with
probability $\phi \left( s,t\right) $ and accepts the null with probability $%
1-\phi \left( s,t\right) $, e.g., the Anderson-Rubin test is simply $I\left(
AR>c\left( k\right) \right) $ where $I\left( \cdot \right) $ is the
indicator function. The test is said to be nonrandomized if $\phi $ only
takes values $0$ and $1$; otherwise, it is called a randomized test. We note
that%
\begin{equation*}
E_{\beta ,\mu }\phi \left( S,T\right) \equiv \int \phi \left( s,t\right)
f_{\beta ,\mu }\left( s,t\right) \text{ }d\left( s,t\right)
\end{equation*}%
is the probability of rejecting the null when the parameters are $\beta $
and $\mu $. The object $E_{\beta ,\mu }\phi \left( S,T\right) $ taken as a
function of $\beta $ and $\mu $ gives the power curve for the test $\phi $.
In particular, $E_{\beta _{0},\mu }\phi \left( S,T\right) $ gives the null
rejection probability. By Tonelli's theorem, we can write 
\begin{equation}
E_{\Lambda }\phi \left( S,T\right) =\int E_{\beta ,\mu }\phi \left(
s,t\right) d\Lambda \left( \beta ,\mu \right) =\int \phi \left( s,t\right)
h_{\Lambda }\left( s,t\right) \text{ }d\left( s,t\right) ,
\end{equation}%
where $h_{\Lambda }\left( s,t\right) $ is defined in (\ref{(WAP density)}).
Hence, $E_{\Lambda }\phi \left( S,T\right) $ is the weighted-average power
for the measure $\Lambda \left( \beta ,\mu \right) $.

A natural first step is to find tests that maximize WAP and have size no
larger than $\alpha $. That is,%
\begin{equation}
\max_{0\leq \phi \leq 1}E_{\Lambda }\phi \left( S,T\right) \text{, where }%
E_{\beta _{0},\mu }\phi \left( S,T\right) \leq \alpha ,\forall \mu .
\end{equation}%
Since the parameter $\mu $ is unknown, finding a WAP test with correct size
is nontrivial. The task entails finding a least favorable distribution $%
\Lambda _{0}$ to construct the WAP test as described in Section 3.8 of %
\citet{LehmannRomano05}. This test rejects the null when the likelihood
ratio is large:%
\begin{equation}
\frac{h_{\Lambda }\left( s,t\right) }{\int f_{\beta _{0},\mu }^{T}\left(
t\right) \text{ }d\overline{\Lambda }\left( \mu \right) }>\kappa ,
\end{equation}%
where $\kappa \cdot \overline{\Lambda }$ is really a Lagrange multiplier in
an infinite-dimensional space; see Lemma 3 of \citet{MoreiraMoreira10} for
details\footnote{%
Also available as Lemma 2 in the most recent version, %
\citet{MoreiraMoreira13}. Both versions are available on Marcelo Moreira's
website: http://www.fgv.br/professor/mjmoreira/}. For a parameter $\mu $ of
small dimension, we can apply numerical algorithms to approximate the WAP
test (such as the one by \citet{ElliottMuellerWatson15} or the linear
programming algorithm of \citet{MoreiraMoreira13}).

The task of finding tests with correct size is simplified if we can find
optimal similar tests:%
\begin{equation}
\max_{0\leq \phi \leq 1}E_{\Lambda }\phi \left( S,T\right) \text{, where }%
E_{\beta _{0},\mu }\phi \left( S,T\right) =\alpha ,\forall \mu .
\label{(WAP similar)}
\end{equation}%
Because the statistic $T$ is sufficient and complete under the null, any
similar test is conditionally similar (for almost all levels $T=t$). Hence,
we can solve%
\begin{equation*}
\max_{0\leq \phi \leq 1}E_{\Lambda }\phi \left( S,t\right) \text{, where }%
E_{\beta _{0}}\phi \left( S,t\right) =\alpha .
\end{equation*}%
The WAP similar test rejects the null when%
\begin{equation}
\frac{h_{\Lambda }\left( s,t\right) }{f_{\beta _{0}}^{S}\left( s\right)
\cdot h_{\Lambda }^{T}\left( t\right) }>\kappa \left( t\right) ,
\label{(WAP similar hT)}
\end{equation}%
where $\kappa \left( t\right) $ is a conditional critical value function and 
$h_{\Lambda }^{T}\left( t\right) =\int h_{\Lambda }\left( s,t\right) $ $ds$.
By Tonelli's theorem, 
\begin{eqnarray*}
h_{\Lambda }^{T}\left( t\right) &=&\int \int f_{\beta ,\mu }\left(
s,t\right) \text{ }d\Lambda \left( \beta ,\mu \right) \text{ }ds \\
&=&\int \int f_{\beta ,\mu }\left( s,t\right) ds\text{ }d\Lambda \left(
\beta ,\mu \right) \\
&=&\int f_{\beta ,\mu }^{T}\left( t\right) \text{ }d\Lambda \left( \beta
,\mu \right) .
\end{eqnarray*}

For arbitrary weights $\Lambda $, neither the WAP test with correct size nor
the WAP similar test is guaranteed to have overall good power in finite
samples\footnote{%
As the geneticist and statistician Anthony W. F. \citet[p.
60]{Edwards92} remarks, \textquotedblleft It is sometimes said, in defence
of the Bayesian concept, that the choice of prior distribution is
unimportant in practice, because it hardly influences the posterior
distribution at all when there are moderate amounts of data. The less said
about this `defence' the better.\textquotedblright}. Take for a moment the
case where $\Sigma =\Omega \otimes \Phi $. The WAP tests based on $%
h_{1}\left( s,t\right) $ can have very low power for some parameter values.
Because the WAP test with correct size and the WAP similar test based on the
MM1 weight are not sign invariant, they can actually behave like one-sided
tests for parts of the parameter space.

This issue is analogous to the problem with conditional Wald tests found by %
\citet{AndrewsMoreiraStock07} which leads them to give a very specific
recommendation: \textquotedblleft \textit{The evident conclusion for applied
work is that researchers choosing among these tests (including conditional
Wald) should use the CLR test. The strong asymptotic bias and often low
power of the conditional Wald tests indicate that they can yield misleading
inferences and are not useful, even as robustness checks.}%
\textquotedblright\ For our purposes we can of course circumvent this
problem by replacing $h_{1}\left( s,t\right) $ by a sign invariant weight
given by (\ref{(correct weight)}) or by the density $h_{2}\left( s,t\right) $%
. However, this solution relies on model symmetries (i.e., sign invariance)
and only works for Kronecker covariance matrices.

On the other hand, \citet{MillsMoreiraVilela14} find approximately unbiased
Wald tests which have overall good power. Their procedure only works for the
model with homoskedastic errors, but it does hint that imposing additional
constraints can actually help to obtain optimal tests with overall good
power for general $\Sigma $.

\section{Two-Sided Boundary Conditions \label{Conditions Sec}}

The WAP similar test based on $h_{2}\left( s,t\right) $ is a two-sided test
in the homoskedastic case precisely because the sign-group of
transformations preserves the two-sided testing problem when $\Sigma =\Omega
\otimes \Phi $. More specifically, because this test depends only on $Q_{S}$%
, $\left\vert Q_{ST}\right\vert $, and $Q_{T}$ it is locally unbiased; see
Corollary 1 of \citet{AndrewsMoreiraStock06b}. When errors are
autocorrelated and heteroskedastic, however, the covariance $\Sigma $
typically does not have a Kronecker product structure. In this case, the WAP
similar test (or a WAP test with correct size) based on $h_{2}\left(
s,t\right) $ may not have good power for parts of the parameter space. Worse
yet, when the covariance matrix lacks Kronecker product structure, there is
actually no sign invariance argument to accommodate two-sided testing.

\bigskip

\begin{proposition}
\label{No sign invariance Prop} Assume that we cannot write $\Sigma $ as $%
\Omega \otimes \Phi $ for a $2\times 2$ matrix $\Omega $ and a $k\times k$
matrix $\Phi $, both symmetric and positive definite. Then for the data
group of transformations $\left[ S,T\right] \rightarrow \left[ \pm S,T\right]
$, there exists no group of transformations in the parameter space which
preserves the testing problem.
\end{proposition}

\bigskip

Proposition \ref{No sign invariance Prop} asserts that we cannot simplify
the two-sided hypothesis testing problem using sign invariance arguments. It
is then much more difficult to find a weight so that the test is, loosely
speaking, two-sided. An unbiasedness condition instead adjusts the weights
automatically (whether $\Sigma $ has a Kronecker product or not). Hence, we
can seek approximately optimal unbiased tests.

An important property of WAP tests is admissibility. Theorem \ref%
{Admissibility Thm} below shows that the WAP unbiased tests are admissible.
The proof follows exactly the same steps as the proof for admissibility of
WAP similar tests of \citet{MoreiraMoreira13} (see Comment 1 after their
Theorem 4)\footnote{\citet{Olea15} provides an alternative proof that
similar tests are admissible by contradiction.}. For completeness, we
provide a proof in the appendix for the following theorem.

\bigskip

\begin{theorem}
\label{Admissibility Thm} Let $\left( \beta ,\mu \right) \in \mathbb{B}%
\times \mathbb{P}$, where both sets compact. Assume that the weight $\Lambda 
$ appearing in (\ref{(WAP density)}) has full support on $\mathbb{B}\times 
\mathbb{P}$. Then there exists a sequence of Bayes' tests $\phi _{m}\left(
s,t\right) $ which weakly converges (in the weak* topology to the $\mathcal{L%
}_{\infty }(\mathbb{R}^{2k})$ space) to the WAP unbiased test. In
particular, the WAP\ unbiased test is admissible.
\end{theorem}

\textbf{Comments: 1. }The weak convergence guarantees, for example, that the
limiting power function of $\phi _{m}\left( s,t\right) $ is the power
function of the WAP\ unbiased test. See \citet{MoreiraMoreira13} for details
on weak convergence of tests.

\textbf{2. }The theorem assumes the parameter space is compact. It may be
possible to drop this assumption with some additional technical conditions;
see \citet{Lehmann52}. The compactness assumption, however, may not be
overly restrictive in practice. First, one could argue that we can pin down
a region large enough in which the parameter lies. Second, the usual
mathematical and statistical software packages have limited numerical
accuracy, so for all practical purposes the weight $\Lambda $ in the average
density $h_{\Lambda }\left( s,t\right) $ has support in a compact set.

\bigskip

Proposition \ref{No sign invariance Prop} shows that there is no sign group
structure which preserves the null and alternative. This makes the task of
finding a weight function $h_{\Lambda }\left( s,t\right) $ which yields a
WAP\ unbiased test difficult with HAC errors. Instead of seeking a weight
function $\Lambda $ so that the WAP test is approximately unbiased, we can
select an arbitrary weight and find the optimal test among unbiased tests;
see \citet{MoreiraMoreira13}. In practice, it would be computationally
intensive to handle so many constraints of the form $E_{\beta ,\mu }\phi
\left( S,T\right) \geq E_{\beta _{0},\mu _{0}}\phi \left( S,T\right) $ for
any scalar $\beta $ and $k$-dimensional vectors $\mu $ and $\mu _{0}$,
especially when $k$ is large. Instead we choose two different restrictions.
The first condition is based on the local power around the null hypothesis.
It is a weaker condition than unbiasedness, so it does not rule out
admissibility. The second condition is a stronger requirement but is easier
to implement. Better yet, numerical simulations will show it yields little
power reduction compared to the first condition. Both conditions and their
associated WAP tests are presented next.

\subsection{Locally Unbiased (LU) Condition \label{LU subsec}}

If the test is unbiased, the derivative of the power function must be equal
to zero under the null. The next proposition uses this fact and completeness
of $T$ to provide a necessary condition for a test to be unbiased. This
locally unbiased (LU) condition states that the test must be similar and
uncorrelated with linear combinations (which depend on the instruments'
coefficient $\mu $) of the pivotal statistic $S$.

\bigskip

\begin{proposition}
\label{LU Prop} A test is said to be locally unbiased (LU) if%
\begin{equation}
E_{\beta _{0},\mu }\phi \left( S,T\right) =\alpha \text{ and }E_{\beta
_{0},\mu }\phi \left( S,T\right) S^{\prime }C_{\beta _{0}}\mu =0\text{, }%
\forall \mu .  \tag{LU}  \label{(LU eq)}
\end{equation}%
If a test is unbiased, then it is LU.
\end{proposition}

\bigskip

In the case $k=1$ where the model is exactly identified, we have an
optimality result for any choice of $\Lambda $. The Anderson-Rubin test is
the uniformly most powerful unbiased (UMPU) test and has power function
depending on the noncentrality parameter $\left( \beta -\beta _{0}\right)
^{2}C_{\beta _{0}}^{2}\mu ^{2}$. We can prove this result directly from
Theorem 2-(a) of \citeauthor{Moreira01} (\citeyear{Moreira01}, %
\citeyear{Moreira09a}) for homoskedastic errors (with the scalar $\mu $ and
matrix $\Omega $ being replaced by $\mu _{\Phi }$ and $\Sigma $). As this
setup resembles the just-identified model with homoskedastic errors,
optimality of the Anderson-Rubin test for HAC errors and $k=1$ follows
straightforwardly.

\bigskip

\begin{proposition}
\label{Just Ident Prop} If $k=1$, the Anderson-Rubin test is the uniformly
most powerful unbiased test and has a power function given by%
\begin{equation*}
P_{\beta ,\mu }\left( AR>c\left( 1\right) \right) =1-G\left( c\left(
1\right) ;\frac{\left( \beta -\beta _{0}\right) ^{2}\mu ^{2}}{b_{0}^{\prime
}\Sigma b_{0}}\right) ,
\end{equation*}%
where $G\left( \cdot ;\delta ^{2}\right) $ is the noncentral $\chi
^{2}\left( 1\right) $ distribution function with noncentrality parameter $%
\delta ^{2}$. Furthermore, the LM\ and CQLR tests are equivalent to the
Anderson-Rubin test, and are also optimal.
\end{proposition}

\bigskip

Following Proposition \ref{LU Prop}, the WAP-LU test solves 
\begin{equation}
\underset{0\leq \phi \leq 1}{\max }E_{\Lambda }\phi \left( S,T\right) \text{%
, where }E_{\beta _{0},\mu }\phi \left( S,T\right) =\alpha \text{ and }%
E_{\beta _{0},\mu }\phi \left( S,T\right) S^{\prime }C_{\beta _{0}}\mu
=0,\forall \mu .  \label{(WAP-LU)}
\end{equation}%
The optimal tests based on $h_{1}\left( s,t\right) $ and $h_{2}\left(
s,t\right) $ are denoted respectively MM1-LU and MM2-LU tests. In the
just-identified model, the MM1-LU test is shown to be the uniformly most
powerful unbiased test. The MM2-LU\ test is equivalent to the MM2 similar
test and is also optimal.

\bigskip

\begin{proposition}
\label{Just Ident LU prop} The following hold when $k=1$:\newline
\emph{(a)} The MM2-LU\ and MM2 similar tests are equivalent and uniformly
most powerful unbiased tests.\newline
\emph{(b)} Both MM1-LU and MM2-LU tests are uniformly most powerful unbiased
tests.
\end{proposition}

\textbf{Comments: 1. }The MM2 similar test automatically satisfies the LU
condition when $k=1$. Hence, the MM2-LU and MM2 similar tests are equivalent
when the model is exactly identified.

\textbf{2.} The MM1 similar test is not locally unbiased even when $k=1$.
Close inspection of the weighted density $h_{1}\left( s,t\right) $ shows
that $d_{\beta }/c_{\beta }$ is the relative contribution of the one-sided $%
S\cdot T$ statistic to the $AR=S^{2}$ statistic. If $\Sigma $ is close to
being singular (that is, $\left\vert \Sigma \right\vert $ is near zero), the
ratio $d_{\beta }/c_{\beta }$ can diverge to infinity. The MM1 test can then
behave as a one-sided test. We will illustrate this problem numerically in
Section \ref{Numerical Sec}.

\bigskip

In the case $k>1$ where the model is overidentified, we no longer have a
uniformly most powerful unbiased test. However, we can still find WAP tests
which are locally unbiased. Relaxing both constraints in (\ref{(WAP-LU)})
assures us the existence of Lagrange multipliers; see %
\citet{MoreiraMoreira13}. Therefore, we solve the approximated maximization
problem:%
\begin{eqnarray}
\underset{0\leq \phi \leq 1}{\max }E_{\Lambda }\phi \left( S,T\right) \text{%
, where }\alpha -\epsilon &\leq &E_{\beta _{0},\mu }\phi \left( S,T\right)
\leq \alpha +\epsilon ,\forall \mu
\label{(relaxed optimization Boundary 1 eq)} \\
\text{and }E_{\beta _{0},\mu _{l}}\phi \left( S,T\right) S^{\prime }C_{\beta
_{0}}\mu _{l} &=&0,\text{ for }l=1,...,m,  \notag
\end{eqnarray}%
when $\epsilon $ is small and the number of discretizations $m$ is large.
The optimal test rejects the null hypothesis when%
\begin{equation}
h_{\Lambda }\left( s,t\right) -s^{\prime }C_{\beta
_{0}}\sum_{l=1}^{m}c_{l}^{\epsilon }\mu _{l}f_{\beta _{0},\mu _{l}}\left(
s,t\right) >\int f_{\beta _{0},\mu }\left( s,t\right) \text{ }d\Lambda
_{\epsilon }\left( \mu \right) ,  \label{(WAP-LU test)}
\end{equation}%
where the measure $\Lambda _{\epsilon }$ and the scalars $c_{l}^{\epsilon }$%
, $l=1,...,m$, are multipliers associated to boundary constraints in the
maximization problem (\ref{(relaxed optimization Boundary 1 eq)}).

We can use $f_{\beta _{0},\mu }\left( s,t\right) =$ $f_{\beta
_{0}}^{S}\left( s\right) \times f_{\beta _{0},\mu }^{T}\left( t\right) $ to
write (\ref{(WAP-LU test)}) as%
\begin{equation}
\frac{h_{\Lambda }\left( s,t\right) }{f_{\beta _{0}}^{S}\left( s\right) }%
-s^{\prime }C_{\beta _{0}}\sum_{l=1}^{m}c_{l}^{\epsilon }\mu _{l}f_{\beta
_{0},\mu _{l}}^{T}\left( t\right) >\int f_{\beta _{0},\mu }^{T}\left(
t\right) \text{ }d\Lambda _{\epsilon }\left( \mu \right) .
\end{equation}%
Letting $\epsilon \downarrow 0$, the optimal test rejects the null
hypothesis when 
\begin{equation}
\frac{h_{\Lambda }\left( s,t\right) }{f_{\beta _{0}}^{S}\left( s\right) }%
-s^{\prime }C_{\beta _{0}}\sum_{l=1}^{m}c_{l}\mu _{l}f_{\beta _{0},\mu
_{l}}^{T}\left( t\right) >\kappa \left( t\right) ,
\end{equation}%
where $\kappa \left( t\right) $ is the conditional $1-\alpha $ quantile of%
\begin{equation}
\frac{h_{\Lambda }\left( S,t\right) }{f_{\beta _{0}}^{S}\left( S\right) }%
-S^{\prime }C_{\beta _{0}}\sum_{l=1}^{m}c_{l}\mu _{l}f_{\beta _{0},\mu
_{l}}^{T}\left( t\right) .
\end{equation}%
This representation is very convenient as we can find 
\begin{equation}
\kappa \left( t\right) =\lim_{\epsilon \downarrow 0}\int f_{\beta _{0},\mu
}^{T}\left( t\right) \text{ }d\Lambda _{\epsilon }\left( \mu \right)
\end{equation}%
by numerical approximations of the conditional distribution instead of
searching for an infinite-dimensional multiplier $\Lambda _{\epsilon }$. We
then search for the values $c_{l}$ so that 
\begin{equation}
E_{\beta _{0},\mu _{l}}\phi \left( S,T\right) S^{\prime }C_{\beta _{0}}\mu
_{l}=\int \phi \left( s,t\right) s^{\prime }C_{\beta _{0}}\mu _{l}f_{\beta
_{0}}^{S}\left( s\right) f_{\beta _{0},\mu _{l}}^{T}\left( t\right) =0,
\end{equation}%
by taking into consideration that $\kappa \left( t\right) $ depends on $%
c_{l} $, $l=1,...,m$. We can find $c_{l}$, $l=1,...,m$ with a nonlinear
numerical algorithm\footnote{%
The two-step procedure just described is the usual \emph{substitution method}
for a system of equations, but here we have an uncountable number of
equations and unknowns.}.

As an alternative procedure, we consider a condition stronger than the LU
condition which is simpler to implement numerically. This strategy turns out
to be useful because it provides a simple way to implement tests with
overall good power. We explain this alternate condition next.

\subsection{Strongly Unbiased (SU) Condition \label{SU subsec}}

The LU\ condition asserts that the test $\phi $ is uncorrelated with a
linear combination indexed by the instruments' coefficients $\mu $ and the
pivotal statistic $S$. We note that the LU\ condition trivially holds if 
\begin{equation}
E_{\beta _{0},\mu }\phi \left( S,T\right) =\alpha \text{ and }E_{\beta
_{0},\mu }\phi \left( S,T\right) S=0,\forall \mu .  \tag{SU}  \label{(SU eq)}
\end{equation}%
That is, the test $\phi $ is uncorrelated with the $k$-dimensional statistic 
$S$ itself under the null. This strongly unbiased (SU)\ condition states
that the test $\phi \left( S,T\right) $ is uncorrelated with $S$ for all
instruments' coefficients $\mu $. The WAP-SU test based on the weight $%
\Lambda $ solves%
\begin{equation}
\underset{0\leq \phi \leq 1}{\max }E_{\Lambda }\phi \left( S,T\right) \text{%
, where }E_{\beta _{0},\mu }\phi \left( S,T\right) =\alpha \text{ and }%
E_{\beta _{0},\mu }\phi \left( S,T\right) S=0,\forall \mu .  \label{(WAP-SU)}
\end{equation}%
The optimal tests based on $h_{1}\left( s,t\right) $ and $h_{2}\left(
s,t\right) $ are denoted respectively MM1-SU and MM2-SU tests.

When $k=1$, the LU\ and SU conditions are equivalent (hence, the MM1-SU and
MM2-SU tests are uniformly most powerful unbiased). When $k>1$, the
following lemma proves the LU\ condition is strictly weaker than the SU
condition. Hence, finding WAP similar tests that satisfy the SU instead of
the LU\ condition in theory may entail unnecessary power losses. In
practice, numerical simulations in Section \ref{Numerical Sec} indicate that
there is little power gain --if any-- by using the LU\ instead of the SU\
condition (with the MM1-SU and MM2-SU tests having the advantage of being
easier to implement).

\bigskip

\begin{lemma}
\textbf{\label{LU not SU Lemma} }Define the integral%
\begin{equation*}
F_{\phi }(\mu _{1},\mu _{2})=E_{\beta _{0},D_{\beta _{0}}^{-1}\mu _{2}}\phi
\left( s,t\right) s^{\prime }C_{\beta _{0}}\mu _{1}=\int \phi \left(
s,t\right) s^{\prime }C_{\beta _{0}}\mu _{1}\cdot f_{\beta _{0}}^{S}\left(
s\right) f_{\beta _{0},D_{\beta _{0}}^{-1}\mu _{2}}^{T}\left( t\right) \text{
}d\left( s,t\right) .
\end{equation*}%
For $k>1$, there exists a test function $\phi :\left[ S,T\right] \rightarrow %
\left[ 0,1\right] $ such that $F_{\phi }(\mu _{1},\mu _{1})=0$ for all $\mu
_{1}$, and $F_{\phi }(\mu _{1},\mu _{2})\neq 0$, for some $\mu _{1}$ and $%
\mu _{2}$.
\end{lemma}

\bigskip

Because the statistic $T$ is complete, we can carry on power maximization in
(\ref{(WAP-SU)}) for each level of $T=t$:%
\begin{equation}
\underset{0\leq \phi \leq 1}{\max }E_{\Lambda }\phi \left( S,t\right) \text{%
, where }E_{\beta _{0}}\phi \left( S,t\right) =\alpha \text{ and }E_{\beta
_{0}}\phi \left( S,t\right) S=0,  \label{(boundary SU eq)}
\end{equation}%
where the expectation is taken with respect to $S$ only. The WAP-SU test
rejects the null when 
\begin{equation*}
\frac{h_{\Lambda }\left( s,t\right) }{f_{\beta _{0}}^{S}\left( s\right)
\cdot h_{\Lambda }^{T}\left( t\right) }>\kappa \left( s,t\right) ,
\end{equation*}%
where the function $\kappa \left( s,t\right) =\overline{\kappa }_{0}\left(
t\right) +s^{\prime }\overline{\kappa }_{1}\left( t\right) $ is such that
the optimal test satisfies the SU condition. The term $h_{\Lambda
}^{T}\left( t\right) $ can be absorbed in the critical value function. For
numerical stability, however, we recommend keeping it so that the numerator
and denominator are of the same order of magnitude.

In practice, we can find $\overline{\kappa }_{0}\left( t\right) $ and $%
\overline{\kappa }_{1}\left( t\right) $ using linear programming based on
simulations for the statistic $S$. Consider the approximated problem 
\begin{eqnarray*}
\max_{0\leq x^{\left( j\right) }\leq 1} &&\text{ }J^{-1}\sum_{j=1}^{J}x^{%
\left( j\right) }\frac{h_{\Lambda }\left( s^{\left( j\right) },t\right) }{%
h_{\Lambda }^{T}\left( t\right) }\exp \left( s^{\left( j\right) \prime
}s^{\left( j\right) }/2\right) \left( 2pi\right) ^{k/2} \\
\text{s.t.} &&\text{ }J^{-1}\sum_{j=1}^{J}x^{\left( j\right) }=\alpha \text{
and} \\
&&\text{ }J^{-1}\sum_{j=1}^{J}x^{\left( j\right) }s_{l}^{^{\left( j\right)
}}=0,\text{ for }l=1,...,k.
\end{eqnarray*}%
Each $j$-th draw of $S$ is iid standard-normal: 
\begin{equation*}
S^{\left( j\right) }=\left[ 
\begin{array}{c}
S_{1}^{\left( j\right) } \\ 
\vdots \\ 
S_{k}^{^{\left( j\right) }}%
\end{array}%
\right] \sim N\left( 0,I_{k}\right) .
\end{equation*}%
We note that for the linear programming, the only term which depends on $T=t$
is $h_{\Lambda }\left( s^{\left( j\right) },t\right) /h_{\Lambda }^{T}\left(
t\right) $. The multipliers for this linear programming problem are the
critical value functions $\overline{\kappa }_{0}\left( t\right) $ and $%
\overline{\kappa }_{1}\left( t\right) $. To speed up the numerical
algorithm, we can use the same sample $S^{\left( j\right) }$, $j=1,...,J,$
for every level $T=t$.

Finally, we use the WAP test found in (\ref{(boundary SU eq)}) to find a
useful \textit{two-sided power envelope}. The next proposition finds the
optimal test for any given alternative which satisfies the SU\ condition.

\bigskip

\begin{proposition}
\label{POSU Prop} The optimal SU test for a point alternative $\left( \beta
,\mu \right) $ rejects the null hypothesis when 
\begin{equation}
\frac{\left( s^{\prime }C_{\beta _{0}}\mu \right) ^{2}}{\mu C_{\beta
_{0}}^{2}\mu }>c(1).  \label{(POSU test eq)}
\end{equation}%
This test is denoted the Point Optimal Strongly Unbiased (POSU) test and has
power given by%
\begin{equation*}
P_{\beta ,\mu }\left( \frac{\left( s^{\prime }C_{\beta _{0}}\mu \right) ^{2}%
}{\mu C_{\beta _{0}}^{2}\mu }>c\left( 1\right) \right) =1-G\left( c\left(
1\right) ;\left( \beta -\beta _{0}\right) ^{2}\mu ^{\prime }C_{\beta
_{0}}^{2}\mu \right) ,
\end{equation*}%
where $G\left( \cdot ;\delta ^{2}\right) $ is the noncentral $\chi
^{2}\left( 1\right) $ distribution function with noncentrality parameter $%
\delta ^{2}$.
\end{proposition}

\textbf{Comments: 1. }The POSU test does not depend on $\beta $ but does
depend on the direction of the vector $C_{\beta _{0}}\mu $.

\textbf{2.} When $k=1$, the Anderson-Rubin and POSU\ tests are the same.

\bigskip

The power plot of $1-G\left( c\left( 1\right) ;\left( \beta -\beta
_{0}\right) ^{2}\mu ^{\prime }C_{\beta _{0}}^{2}\mu \right) $ as $\beta $
and $\mu $ change yields the two-sided power envelope. This power envelope
is the two-sided analogue of the one-sided power envelope among similar
tests. This power upper bound, based on the Point Optimal Similar (POS) test
for the alternative $\left( \beta ,\mu \right) $, is given by the plot of $%
1-\Phi \left( \sqrt{c\left( 1\right) }-\left\vert \beta -\beta
_{0}\right\vert \sqrt{\mu C_{\beta _{0}}^{2}\mu }\right) $, where $\Phi
\left( \cdot \right) $ is the standard normal distribution.

\section{Numerical Evaluation of WAP Tests \label{Numerical Sec}}

In this section, we provide numerical simulations for WAP\ tests based on
the MM statistics. The MM tests are WAP similar tests based on $h_{1}\left(
s,t\right) $ and $h_{2}\left( s,t\right) $. The MM-LU\ and MM-SU tests also
satisfy respectively the locally unbiased and strongly unbiased conditions.
The goal in this section is to numerically illustrate the importance of
using two-sided conditions to obtain tests with overall good power.

We can write 
\begin{equation*}
\Omega =\left[ 
\begin{array}{cc}
\omega _{11}^{1/2} & 0 \\ 
0 & \omega _{22}^{1/2}%
\end{array}%
\right] P_{\Omega }\left[ 
\begin{array}{cc}
1+\rho & 0 \\ 
0 & 1-\rho%
\end{array}%
\right] P_{\Omega }^{\prime }\left[ 
\begin{array}{cc}
\omega _{11}^{1/2} & 0 \\ 
0 & \omega _{22}^{1/2}%
\end{array}%
\right] ,
\end{equation*}%
where $P_{\Omega }$ is an orthogonal matrix and $\rho =\omega _{12}/\omega
_{11}^{1/2}\omega _{22}^{1/2}$. For the numerical simulations, we specify $%
\omega _{11}=\omega _{22}=1$.

We use the decomposition of $\Omega $ to perform numerical simulations for a
class of covariance matrices:%
\begin{equation*}
\Sigma =P_{\Omega }\left[ 
\begin{array}{cc}
1+\rho & 0 \\ 
0 & 0%
\end{array}%
\right] P_{\Omega }^{\prime }\otimes diag\left( \varsigma _{1}\right)
+P_{\Omega }\left[ 
\begin{array}{cc}
0 & 0 \\ 
0 & 1-\rho%
\end{array}%
\right] P_{\Omega }^{\prime }\otimes diag\left( \varsigma _{2}\right) ,
\end{equation*}%
where $\varsigma _{1}$ and $\varsigma _{2}$ are $k$-dimensional vectors.

We consider two possible choices for $\varsigma _{1}$ and $\varsigma _{2}$.
For the first design, we set $\varsigma _{1}=\varsigma _{2}=\left(
1/\varepsilon -1,1,...,1\right) ^{\prime }$. The covariance matrix then
simplifies to a Kronecker product: $\Sigma =\Omega \otimes diag\left(
\varsigma _{1}\right) $. For the non-Kronecker design, we set $\varsigma
_{1}=\left( 1/\varepsilon -1,1,...,1\right) ^{\prime }$ and $\varsigma _{2}=$
$\left( 1,...,1,1/\varepsilon -1\right) ^{\prime }$. This setup captures the
data asymmetry in extracting information about the parameter $\beta $ from
each instrument. For small $\varepsilon $, the angle between $\varsigma _{1}$
and $\varsigma _{2}$ is nearly $90^{\circ }$. We report numerical
simulations for $\varepsilon =\left( k+1\right) ^{-1}$. As $k$ increases,
the vector $\varsigma _{1}$ becomes orthogonal to $\varsigma _{2}$ in the
non-Kronecker design.

We set the parameter $\mu =\left( \lambda ^{1/2}/\sqrt{k}\right) 1_{k}$ for $%
k=2,5,10,20$ and $\rho =-0.5,0.2,0.5,0.9$. We choose $\lambda
/k=0.5,1,2,4,8,16$, which span the range from weak to strong instruments. We
focus on tests with significance level 5\% for testing $\beta _{0}=0$. To
conserve space, we report here only power plots for $k=5$, $\rho =0.9$, and $%
\lambda /k=2,8$. The full set of simulations is available on Marcelo
Moreira's website.

We present plots for the power envelope and power functions against various
alternative values of $\beta $ and $\lambda $. All results reported here are
based on 1,000 Monte Carlo simulations. We plot power as a function of the
rescaled alternative $\left( \beta -\beta _{0}\right) \lambda ^{1/2}$, which
reflects the difficulty in making inference on $\beta $ for different
instruments' strength.

\begin{figure}[tbh]
\caption{Power Comparison (Kronecker Variance)}
\label{fig:Kronecker}\centering \bigskip \minipage{0.5\textwidth} \centering %
\includegraphics[width=5.5cm]{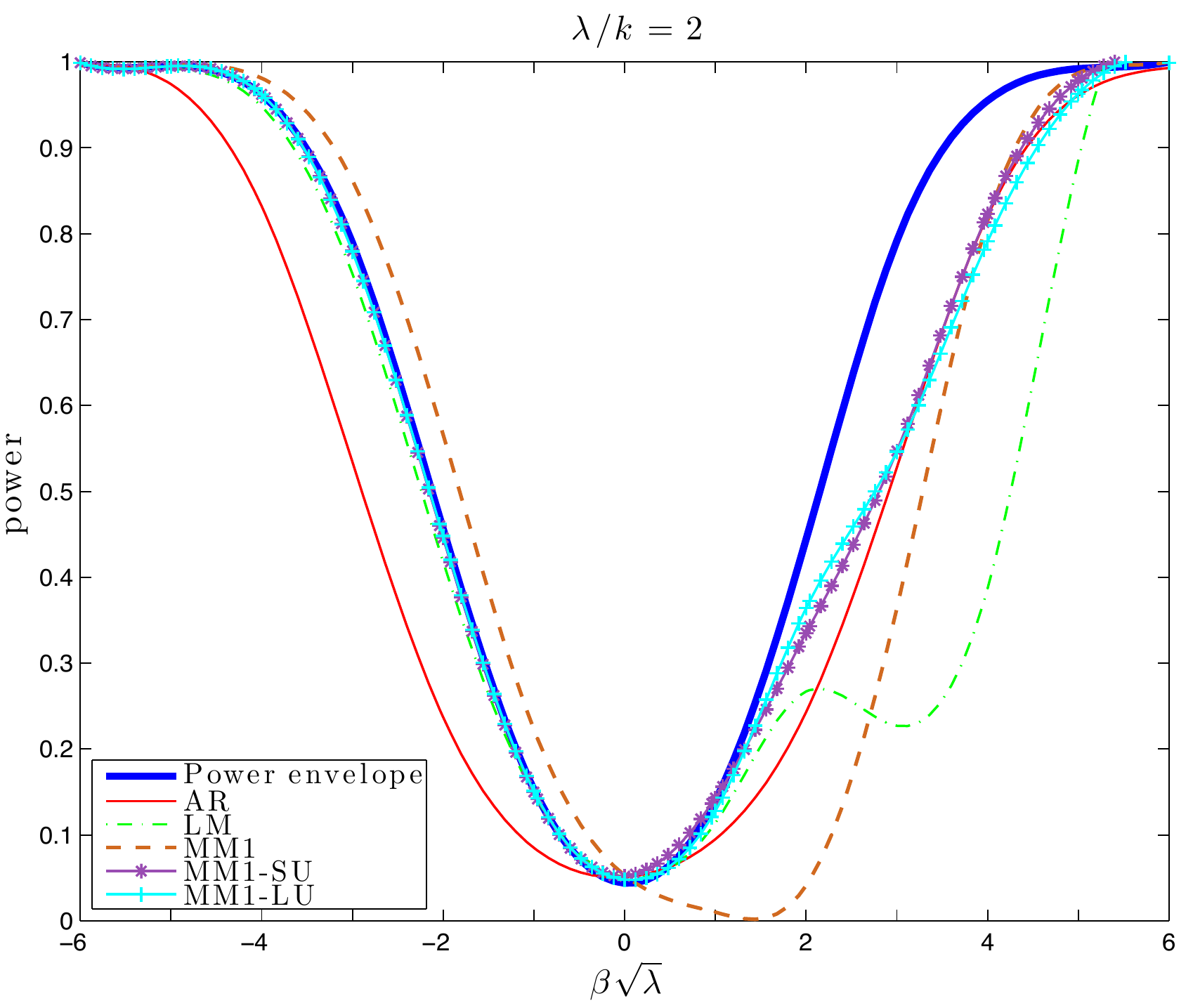} \endminipage%
\hfill \minipage{0.5\textwidth} \centering %
\includegraphics[width=5.5cm]{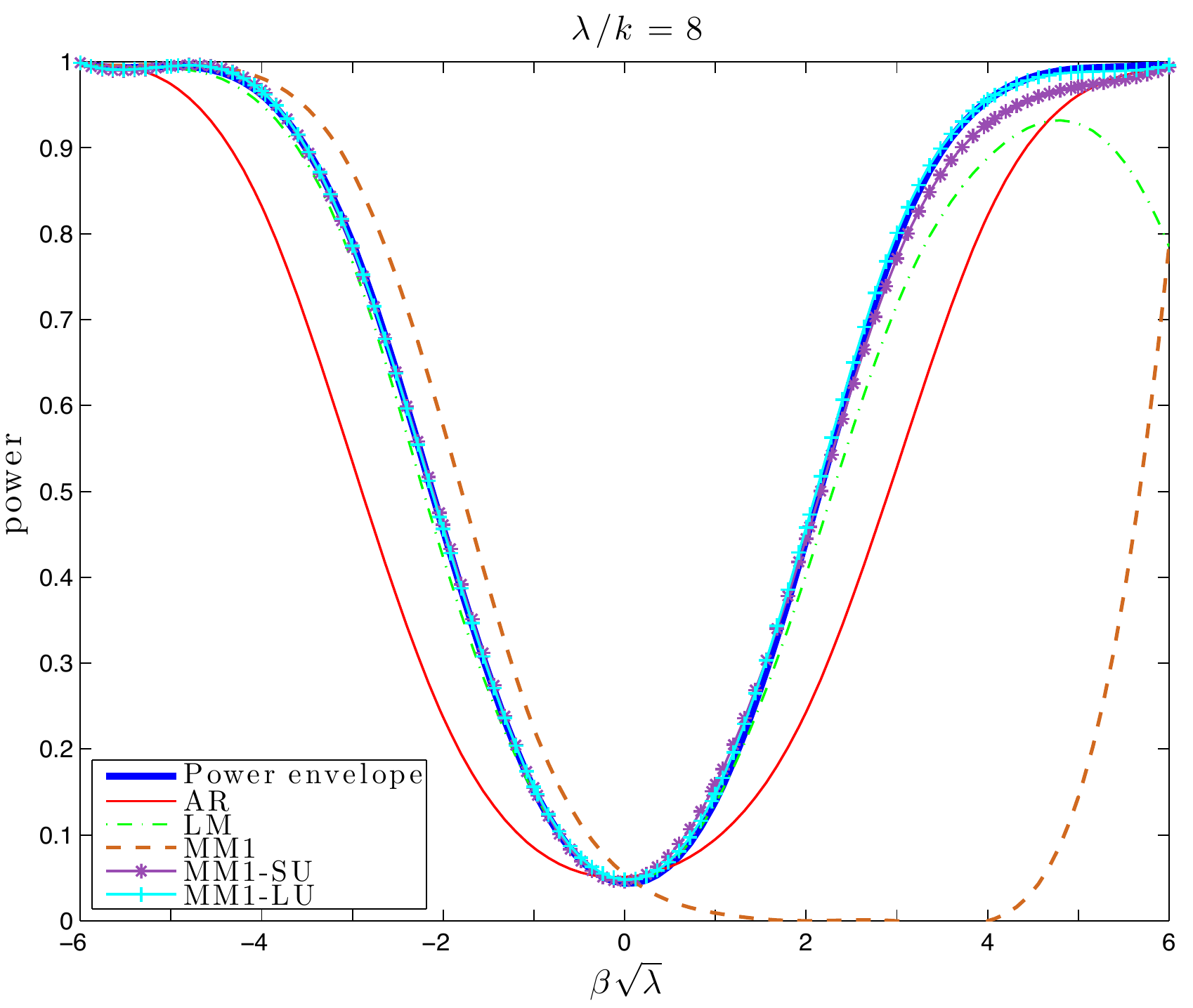} \endminipage %
\hfill \minipage{0.5\textwidth} \bigskip \centering %
\includegraphics[width=5.5cm]{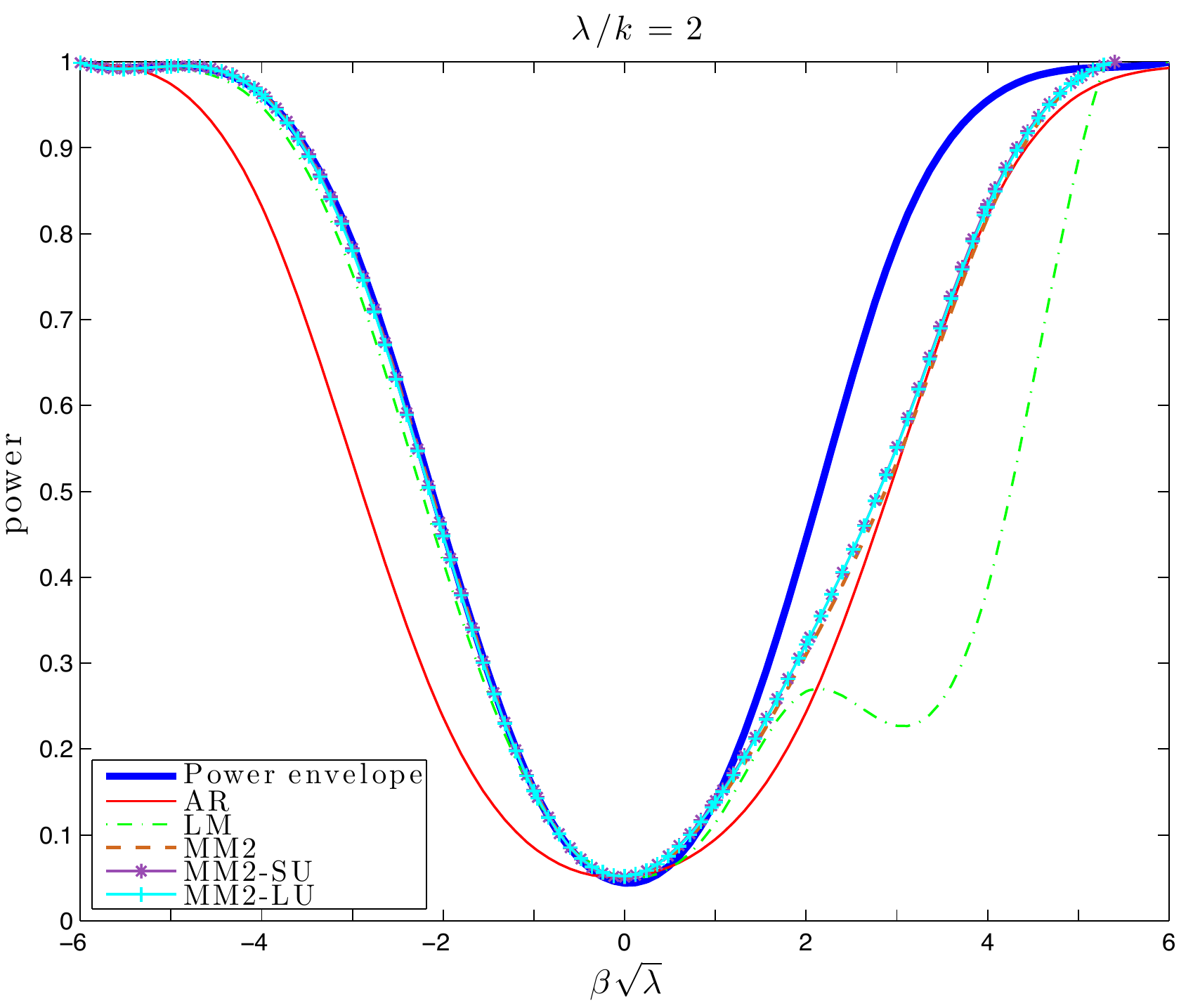} \endminipage%
\hfill \minipage{0.5\textwidth} \bigskip \centering %
\includegraphics[width=5.5cm]{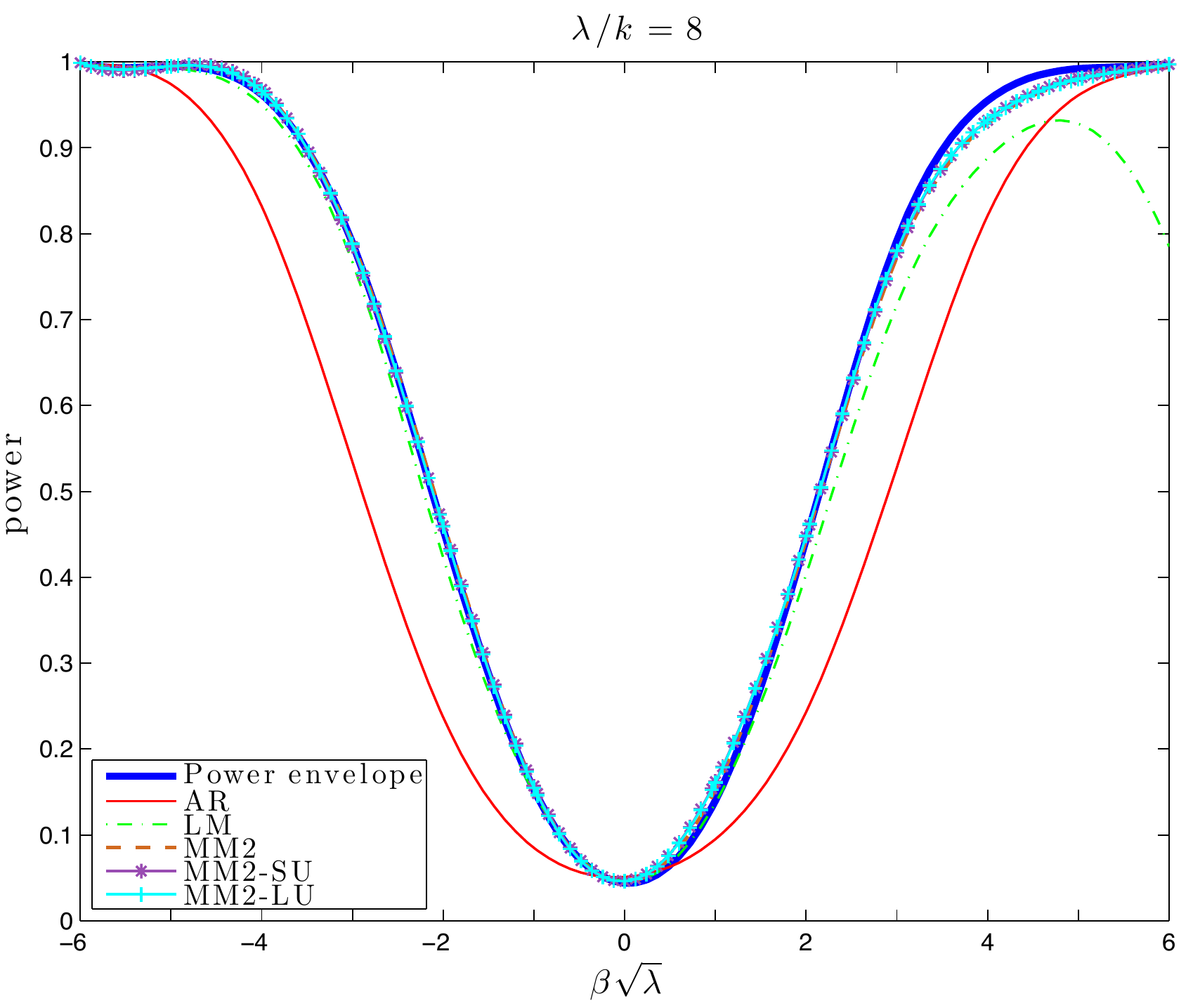} \endminipage %
\hfill
\end{figure}

Figure \ref{fig:Kronecker} reports numerical results for the Kronecker
product design. All four pictures present the power envelope and power
curves for two existing tests, the Anderson-Rubin ($AR$) and score ($LM$)
tests.

The first two graphs plot the power curves for the three WAP tests based on
the MM1 statistic with $\sigma ^{2}=10$. All three tests reject the null
when the $h_{1}\left( s,t\right) $ statistic is larger than an adjusted
critical value function. In practice, we approximate these critical value
functions with 10,000 replications. The MM1 test sets the critical value
function to be the 95\% empirical quantile of $h_{1}\left( S,t\right) $. The
MM1-SU\ test uses a conditional linear programming algorithm to find its
critical value function. The MM1-LU test uses a nonlinear optimization
package.

The AR test has power considerably lower than the power envelope when
instruments are both weak ($\lambda /k=2$) and strong ($\lambda /k=8$). The
LM test does not perform well when instruments are weak, and its power
function is not monotonic even when instruments are strong. These two facts
about the AR\ and LM tests are well documented in the literature; see %
\citet{Moreira03} and AMS06. The figure also reveals some salient findings
for the tests based on the MM1 statistic. First, all MM1-based tests have
correct size. Second, the MM1 similar test can have large bias to the point
that it has zero power for parts of the parameter space. Hence, a naive
choice for the density can yield a WAP test which can have overall poor
power. We can eliminate this problem by imposing an unbiased condition when
selecting an optimal test. The MM1-SU\ test is easy to implement and has
power closer to the power upper bound. When instruments are weak, its power
lies moderately below the reported power envelope. This is expected as the
number of parameters is too large\footnote{%
The MM1-SU power is nevertheless close to the two-sided power envelope for
orthogonally invariant tests as in AMS06 (which is applicable to this
design, but not reported here).}. When instruments are strong, its power is
virtually the same as the power envelope.

To support the use of the MM1-SU\ test we also consider the MM1-LU test,
which imposes a weaker unbiased condition. Close inspection of the graphs
show that the derivative of the power function of the MM1 test is different
from zero at $\beta =\beta _{0}$. This observation suggests that the power
curve of the WAP test would change considerably if we were to force the
power derivative to be zero at $\beta =\beta _{0}$. Indeed, we implement the
MM1-LU test where the locally unbiased condition is true at only one point,
the true parameter $\mu $. This parameter is of course unknown to the
researcher and this test is not feasible. However, by considering the
locally unbiased condition for other values of the instruments'
coefficients, the WAP test would be smaller ---not larger. The power curves
of MM1-LU\ and MM1-SU\ tests are very close, which shows that there is not
much to be gained by relaxing the strongly unbiased condition.

The last two graphs plot the power curves for the three WAP tests based on
the MM2 statistic with $\zeta =10$. By using the density $h_{2}\left(
s,t\right) $, we avoid the pitfalls for the MM1 test. Recall that $%
h_{2}\left( s,t\right) $ is invariant to those data transformations which
preserve the two-sided hypothesis testing problem. Hence, the MM2 similar
test is unbiased and has overall good power without imposing any additional
unbiased conditions. The graphs illustrate this theoretical finding, as the
MM2, MM2-SU, and MM2-LU tests have numerically the same power curves. This
conclusion changes dramatically when the covariance matrix is no longer a
Kronecker product.

\begin{figure}[tbh]
\caption{Power Comparison (Non-Kronecker Variance)}
\label{fig:Non-Kronecker}\centering \bigskip \minipage{0.5\textwidth} %
\centering \includegraphics[width=5.5cm]{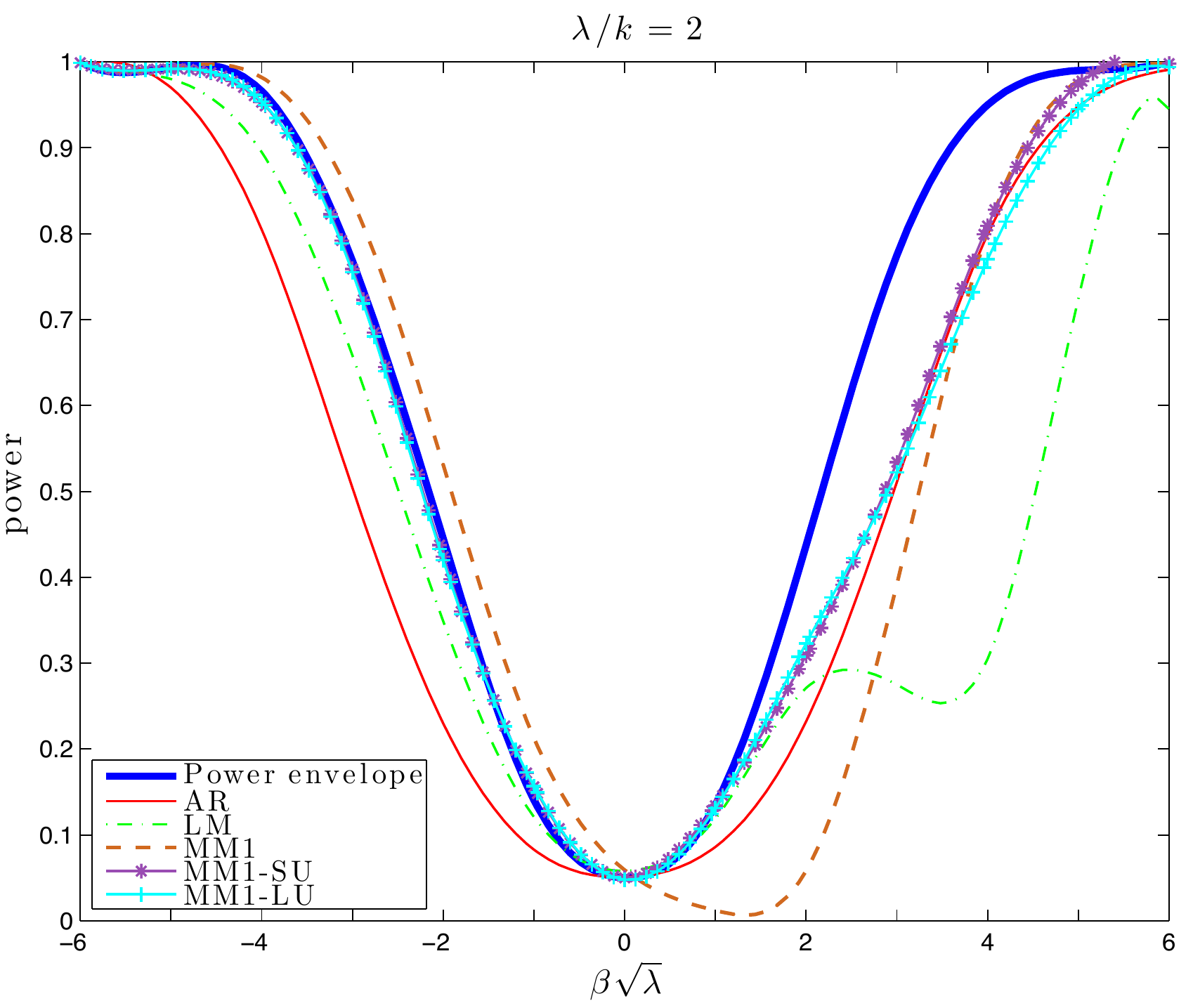} %
\endminipage\hfill \minipage{0.5\textwidth} \centering %
\includegraphics[width=5.5cm]{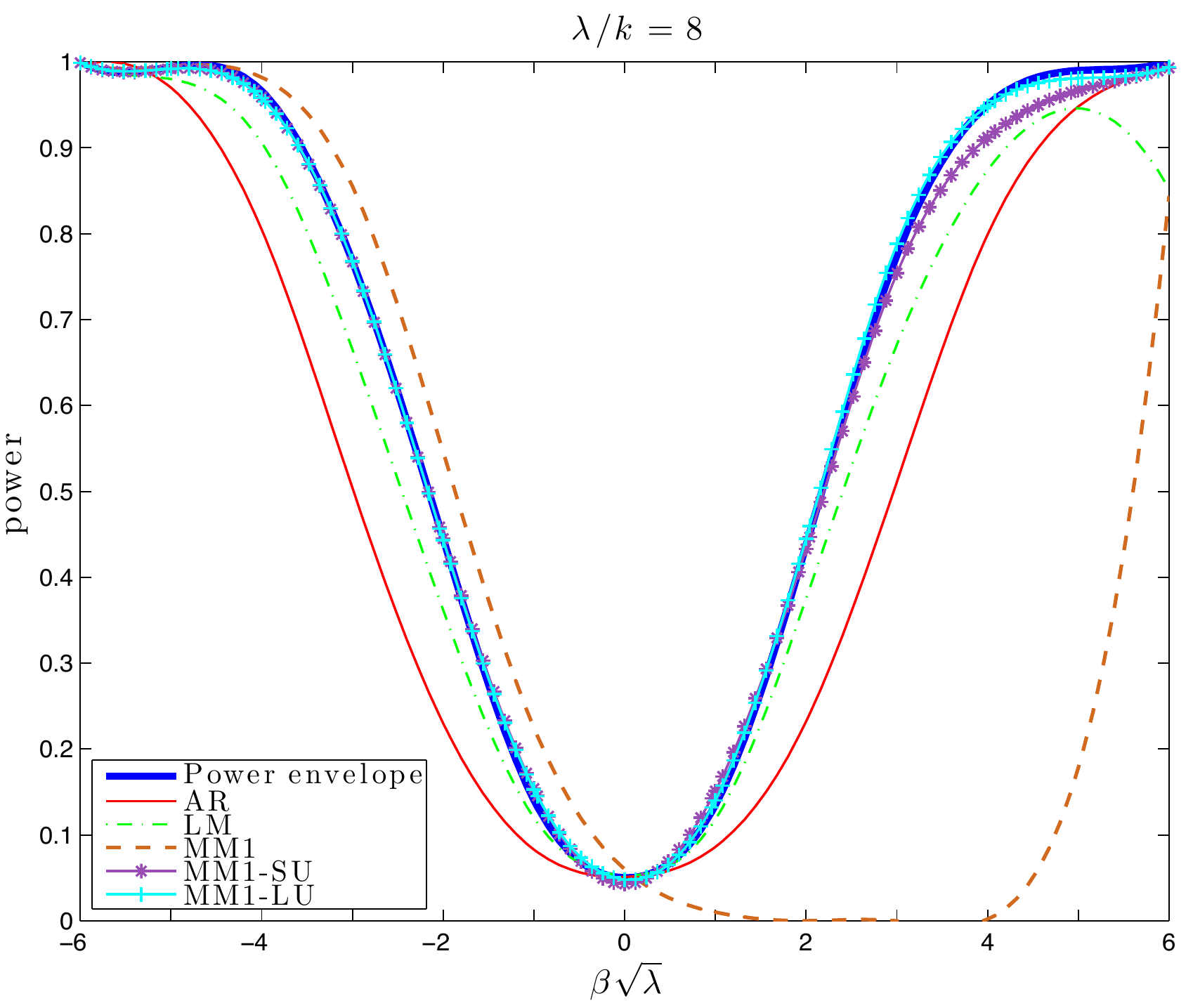} \endminipage %
\hfill \minipage{0.5\textwidth} \bigskip \centering %
\includegraphics[width=5.5cm]{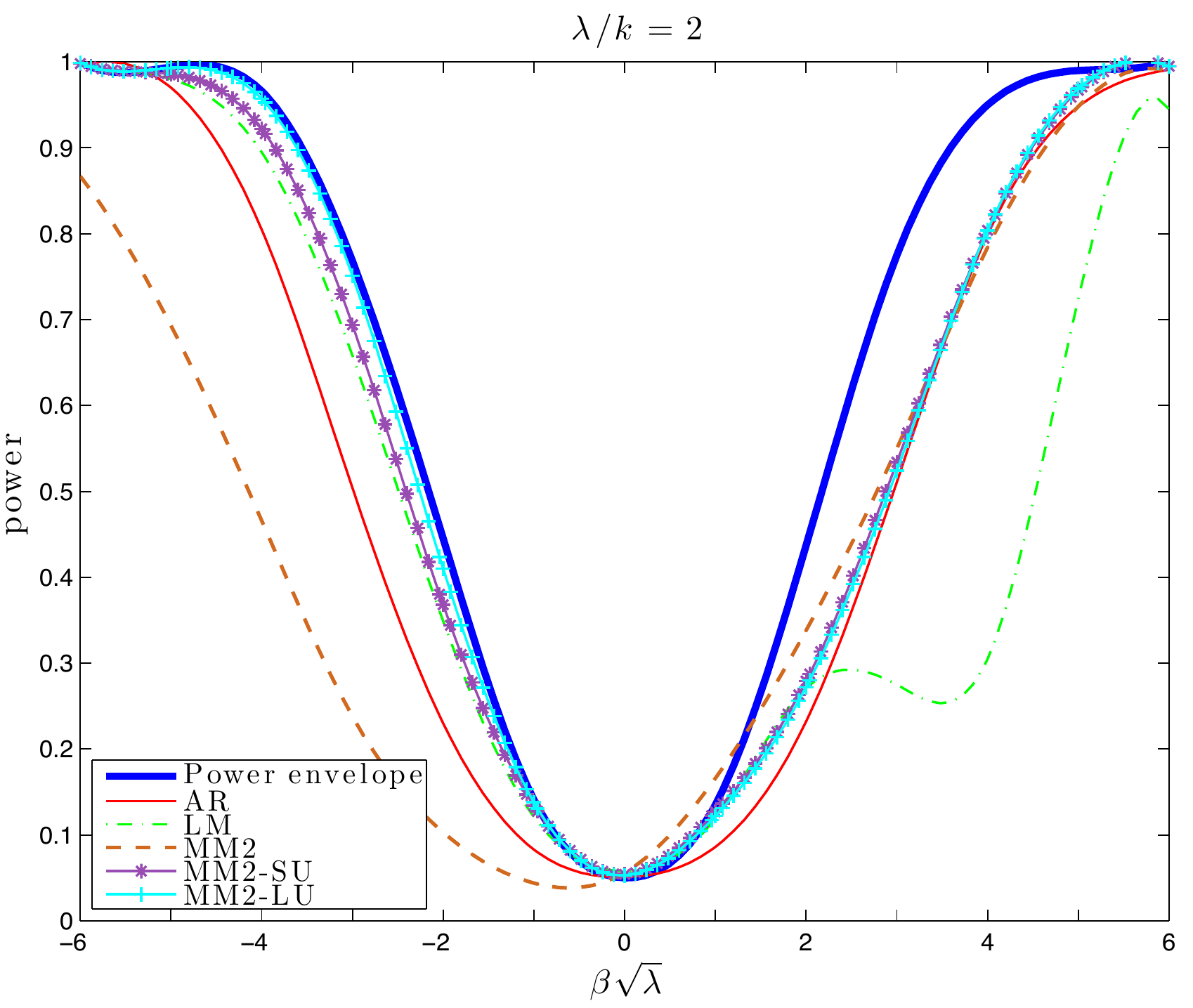} \endminipage%
\hfill \minipage{0.5\textwidth} \bigskip \centering %
\includegraphics[width=5.5cm]{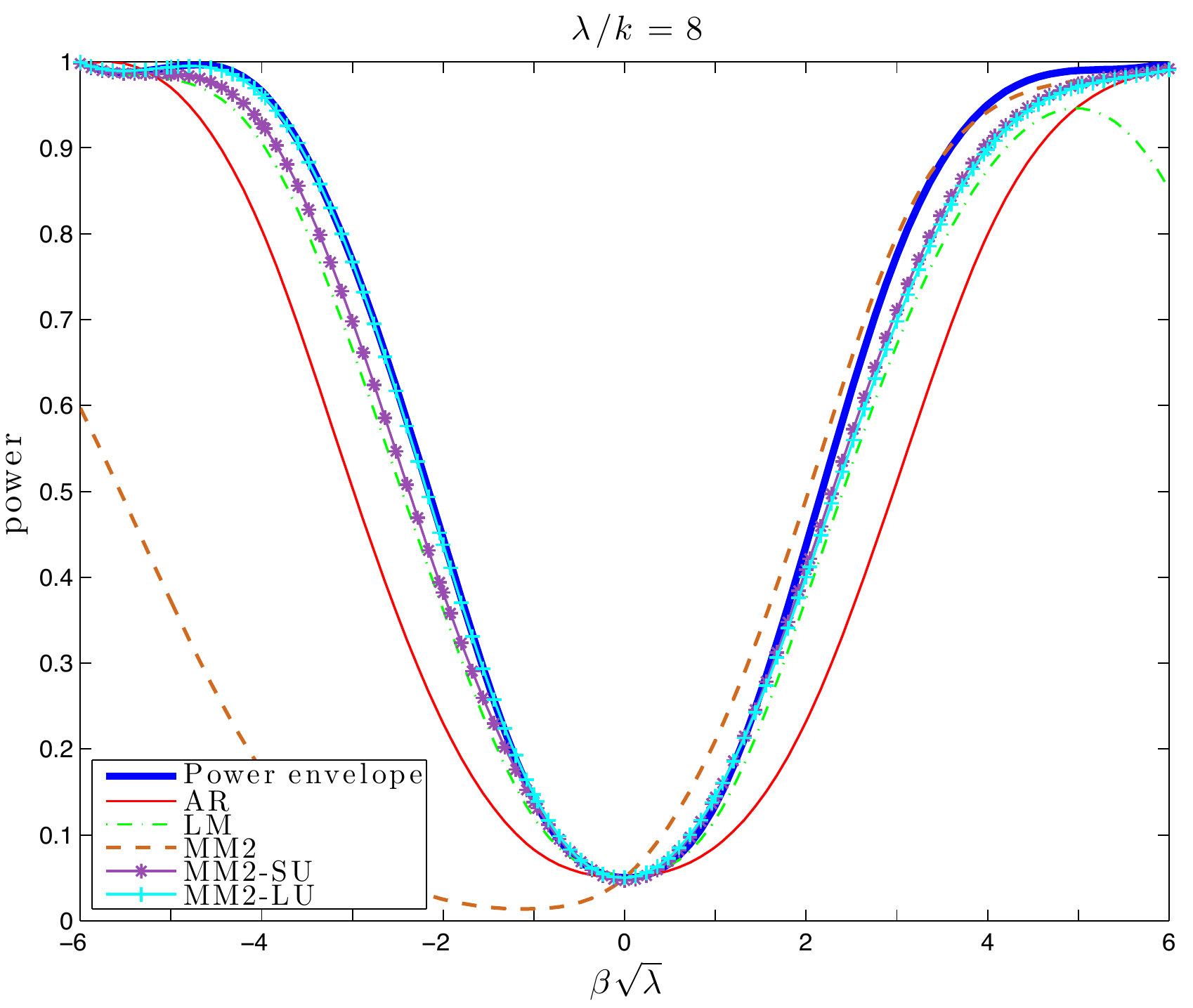} \endminipage %
\hfill
\end{figure}

Figure \ref{fig:Non-Kronecker} presents the power curves for all reported
tests for the non-Kronecker design. Both MM1 and MM2 tests are severely
biased and have overall bad power. For each design, we can make the tests
approximately unbiased by choosing the $\sigma ^{2}$ and $\zeta $ parameters
large enough. However, this unbiasedness control is pointwise in the
parameter space. We can always find a design such that each test behaves as
a one-sided test and has very low power in parts of the parameter space.
Hence, the strong asymptotic bias and often-low power of the conditional
Wald tests found by \citet{AndrewsMoreiraStock07} also hold for the MM1
(even for the homoskedastic IV model) and MM2 similar tests (only for the
HAC-IV model). These WAP similar tests are highly biased with power equal to
zero in some parts of the parameter space. Therefore, just as %
\citet{AndrewsMoreiraStock07} object to the use of conditional Wald tests,
we do not recommend the MM1 and MM2 similar tests for empirical researchers.

Proposition \ref{No sign invariance Prop} shows that we cannot find a group
of data transformations which preserve the two-sided testing problem with
heteroskedastic-autocorrelated errors. Hence, a choice for the density for
the WAP test based on symmetry considerations is not obvious. The correct
density choice can be particularly difficult due to the large
parameter-dimension (the coefficients $\mu $ and covariance $\Sigma $).
Instead, we can endogenize the weight choice so that the WAP test will be
automatically unbiased. This is done by the MM1-LU and MM2-LU tests. These
two tests perform as well as the MM1-SU and MM2-SU\ tests. Because the
latter two tests are easy to implement, we recommend their use in empirical
practice.

\section{Asymptotic Theory \label{Asymptotic Sec}}

All theoretical and numerical results so far do not rely on the sample size $%
n$ at all as we have assumed the statistics $S$ and $T$ to be exactly
normally distributed with known variance $\Sigma $. In this section we relax
this assumption at the cost of asymptotic approximations.

Let $z_{i}$ and $v_{i}$ denote the $i$-th row of $Z$ and $V$, respectively,
written as column vectors of dimensions $k$ and $2$. We make the following
two assumptions as the sample size $n$ grows.

\bigskip

\noindent \textbf{Assumption 1. }$n^{-1}Z^{\prime
}Z=n^{-1}\sum_{i=1}^{n}z_{i}z_{i}^{\prime }\rightarrow _{p}D_{Z}$ for some
positive definite $k\times k$ matrix $D_{Z}$.\smallskip

\bigskip

\noindent \textbf{Assumption 2. }$n^{-1/2}\sum_{i=1}^{n}\left( v_{i}\otimes
z_{i}\right) \rightarrow _{d}N(0,\Sigma _{\infty })$ for some positive
definite $2k\times 2k$ matrix $\Sigma _{\infty }$.

\bigskip

Assumption 1 holds under Birkhoff's Ergodic Theorem. Assumption 2 holds
under suitable conditions by a central limit theorem (CLT). It also assumes
that the long-run covariance matrix of $\Sigma _{\infty }$ is positive
definite, as is usual in the literature. We no longer omit the dependence of 
$\Sigma $ on the sample size $n$ and, hereinafter, write $\Sigma _{n}$.
Assumption 2 asserts that $\Sigma _{\infty }$ is the limit of $\Sigma _{n}$
as $n$ grows. Let $\widehat{\Sigma }_{n}$ be a consistent estimator of $%
\Sigma _{\infty }$ based on $\{\left( \widehat{v}_{i}\otimes z_{i}\right)
:i\leq n\}$, where $\widehat{v}_{i}$ are reduced-form residuals. There are
many HAC estimators in the literature that can be used for this purpose;
see, e.g., \citet{NeweyWest87} and \citet{Andrews91}. For brevity, we do not
provide an explicit set of conditions under which one or more of these HAC
estimators is consistent; see \citet{Jansson02} for details. We note,
however, that the presence of weak instruments does not complicate standard
proofs of the consistency of HAC estimators. Indeed, the convergence for
most estimators holds uniformly over all true parameters $\beta $ and $\pi $.

We now introduce feasible versions of $S_{n}$ and $T_{n}$ with the variance $%
\Sigma _{n}$ replaced by the estimator $\widehat{\Sigma }_{n}$:%
\begin{eqnarray}
\widehat{S}_{n} &=&\left[ \left( b_{0}^{\prime }\otimes I_{k}\right) 
\widehat{\Sigma }_{n}\left( b_{0}\otimes I_{k}\right) \right] ^{-1/2}\left(
b_{0}^{\prime }\otimes I_{k}\right) \overline{R}_{n}\text{ and}
\label{(S^ and T^ defn)} \\
\widehat{T}_{n} &=&\left[ \left( a_{0}^{\prime }\otimes I_{k}\right) 
\widehat{\Sigma }_{n}^{-1}\left( a_{0}\otimes I_{k}\right) \right]
^{-1/2}\left( a_{0}^{\prime }\otimes I_{k}\right) \widehat{\Sigma }_{n}^{-1}%
\overline{R}_{n},  \notag
\end{eqnarray}%
where $\overline{R}_{n}=vec\left[ \left( Z^{\prime }Z\right)
^{-1/2}Z^{\prime }Y\right] $. Likewise, we define the feasible statistic $%
\widehat{\psi }_{n}$ as $\psi \left( S,T,\Sigma ,D_{Z}\right) $ with the
arguments being replaced by their sample analogues:%
\begin{equation}
\widehat{\psi }_{n}=\psi (\widehat{S}_{n},\widehat{T}_{n},\widehat{\Sigma }%
_{n},\widehat{D}_{Z})\text{, where }\widehat{D}_{Z}=n^{-1}Z^{\prime }Z.
\label{(Psi^ and Dz^ defn)}
\end{equation}

\bigskip

\noindent \textbf{Assumption 3. }The prior distribution for $\left( \beta
,\pi \right) $ is absolutely continuous to the Lebesgue measure in $\mathbb{R%
}^{k+1}$. Its density 
\begin{equation*}
w(\beta ,\pi ,\widehat{D}_{Z})=w_{1}(\left. \pi \right\vert \beta ,\widehat{D%
}_{Z})\cdot w_{2}(\beta ,\widehat{D}_{Z})
\end{equation*}%
has full support and is a continuous function of $\pi $ and $\beta $.

\bigskip

Assumption 3 allows the density $w(\beta ,\pi ,\widehat{D}_{Z})$ to depend
on the data through $\widehat{D}_{Z}$. This generalization allows us to
cover all tests considered here and asymptotically behaves as $w(\beta ,\pi
,D_{Z})$ (and so we will omit the dependence of the weights on $\widehat{D}%
_{Z}$ out of convenience). Although the conditional density $w_{1}(\left.
\pi \right\vert \beta )$ does not depend on $\beta $ for the MM1 tests, it
does depend on $\beta $ for the MM2 tests. Assumption 3 also guarantees that
the priors for $\beta $ and $\pi $ are not dogmatic and will vanish
asymptotically as in the Bernstein-von Mises theorem. If we set the prior on 
$\mu $, then the associated prior on $\pi $ $=\left( Z^{\prime }Z\right)
^{1/2}\mu $ depends on the sample size. For example, the MM statistics
introduced in (\ref{(h densities)}) use the prior $\mu \sim N\left( 0,\sigma
^{2}\Phi \right) $. For the associated prior on $\pi \sim N\left( 0,\left(
\sigma ^{2}/n\right) \widehat{D}_{Z}^{-1/2}\Phi \widehat{D}%
_{Z}^{-1/2}\right) $ not to be sensitive to the sample size, the parameters $%
\sigma ^{2}$ and $\zeta $ present in the MM1 and MM2 statistics must
eventually grow at the rate $n$. We make the dependence of $\Lambda \left(
\beta ,\mu \right) $ on the sample size $n$ explicit and, hereinafter, use
the notation $\Lambda _{n}$.

We now analyze the asymptotic behavior of the WAP similar and WAP-SU tests.
Recall that both of these types of tests depend on the test statistic%
\begin{equation}
\frac{h_{\Lambda _{n}}\left( s,t\right) }{f_{\beta _{0}}^{S}\left( s\right)
\cdot h_{\Lambda _{n}}^{T}\left( t\right) }.
\label{(WAP original statistic)}
\end{equation}%
When instruments are weak, the numerator and denominator have the same order
of magnitude. When instruments are strong, the integrands in the weighted
densities $h_{\Lambda _{n}}\left( s,t\right) $ and $h_{\Lambda
_{n}}^{T}\left( t\right) $ grow exponentially fast and we can apply the
Laplace approximation. Because both densities involve $k+1$ integrals, the
test statistic in (\ref{(WAP original statistic)}) is again well-behaved.
The caveat is that a simple, closed-form approximation for $h_{\Lambda
_{n}}^{T}\left( t\right) $ does not seem available under strong instruments.
The WAP similar and WAP-SU tests, however, remain the same if we standardize
(\ref{(WAP original statistic)}) by any function of $t$. We replace $%
h_{\Lambda _{n}}^{T}\left( t\right) $ by $\left( 1+\left\Vert t\right\Vert
\right) ^{-1}h_{\Lambda _{\beta _{0},n}}^{T}\left( t\right) $, where 
\begin{equation}
h_{\Lambda _{\beta _{0},n}}^{T}\left( t\right) =\int f_{\beta _{0},\left(
Z^{\prime }Z\right) ^{1/2}\pi }^{T}\left( t\right) w\left( \beta _{0},\pi
\right) d\pi .  \label{(weighted null density)}
\end{equation}

The WAP similar and WAP-SU tests reject the null when%
\begin{equation}
WAP=\frac{h_{\Lambda _{n}}\left( S,T\right) }{f_{\beta _{0}}^{S}\left(
S\right) \cdot \left( 1+\left\Vert T\right\Vert \right) ^{-1}h_{\Lambda
_{\beta _{0},n}}^{T}\left( T\right) }  \label{(WAP statistic)}
\end{equation}%
is larger than $\kappa _{n}\left( t\right) $ and $\kappa _{n}\left(
s,t\right) $, respectively\footnote{%
The use of a Laplace approximation of the ratio of weighted average under
the alternative and the null is standard under the usual asymptotics. What
is perhaps not standard is the additional term to absorb different rates and
unify nonstandard asymptotics. Indeed, if we were to replace $h_{\Lambda
}^{T}\left( t\right) $ only by $h_{\Lambda _{0,n}}^{T}\left( t\right) $, the
numerator and denominator in (\ref{(WAP original statistic)}) would have
different orders of magnitude under strong instruments.}.

Whether the instruments are weak or strong, we are able to obtain an
approximation to (\ref{(WAP statistic)}). Define%
\begin{eqnarray*}
n\cdot Q_{n}(\beta ,\pi ) &=&\frac{1}{2}\left\Vert \Sigma ^{-1/2}\left( 
\overline{R}-(a\otimes \left( Z^{\prime }Z\right) ^{1/2}\pi )\right)
\right\Vert ^{2} \\
&=&\frac{1}{2}\left\Vert [S:T]-\left[ (\beta -\beta _{0})C_{\beta _{0}}:%
\text{ }D_{\beta }\right] (I_{2}\otimes \left( Z^{\prime }Z\right) ^{1/2}\pi
)\right\Vert ^{2}.
\end{eqnarray*}%
In Appendix B shows that the WAP statistic is asymptotically equivalent to%
\begin{equation}
\frac{\int \exp \left( -n\cdot Q_{n}\left( \beta ,\pi \left( \beta \right)
\right) \right) w\left( \beta ,\pi \left( \beta \right) \right) \left\vert
\left( a^{\prime }\otimes \widehat{D}_{Z}^{1/2}\right) \Sigma
_{n}^{-1}\left( a\otimes \widehat{D}_{Z}^{1/2}\right) \right\vert
^{-1/2}d\beta }{\exp \left( -\frac{S^{\prime }S}{2}\right) \left[
1+\left\Vert T\right\Vert \right] ^{-1}w\left( \beta _{0},\pi \left( \beta
_{0}\right) \right) \left\vert \left( a_{0}^{\prime }\otimes \widehat{D}%
_{Z}^{1/2}\right) \Sigma _{n}^{-1}\left( a_{0}\otimes \widehat{D}%
_{Z}^{1/2}\right) \right\vert ^{-1/2}},  \label{(WAP 1st Laplace)}
\end{equation}%
where the constrained maximum likelihood estimator (MLE) for $\pi $ is%
\begin{eqnarray}
\pi \left( \beta \right)  &=&\left( Z^{\prime }Z\right) ^{-1/2}\left[
(a^{\prime }\otimes I_{k})\Sigma _{n}^{-1}(a\otimes I_{k})\right]
^{-1}(a^{\prime }\otimes I_{k})\Sigma _{n}^{-1}\overline{R}\text{ and} \\
\overline{R} &=&\Sigma _{n}^{1/2}\left[ 
\begin{array}{c}
\left[ \left( b_{0}^{\prime }\otimes I_{k}\right) \Sigma _{n}\left(
b_{0}\otimes I_{k}\right) \right] ^{-1/2}\left( b_{0}^{\prime }\otimes
I_{k}\right) \Sigma _{n}^{1/2} \\ 
\left[ \left( a_{0}^{\prime }\otimes I_{k}\right) \Sigma _{n}^{-1}\left(
a_{0}\otimes I_{k}\right) \right] ^{-1/2}\left( a_{0}^{\prime }\otimes
I_{k}\right) \Sigma _{n}^{-1/2}%
\end{array}%
\right] ^{\prime }\left[ 
\begin{array}{c}
S \\ 
T%
\end{array}%
\right] .  \notag
\end{eqnarray}

The same approximation (\ref{(WAP 1st Laplace)}) holds for the $\widehat{WAP}
$ statistic where we replace $S$, $T$, and $\Sigma $ by their feasible
versions given in (\ref{(S^ and T^ defn)}). The resulting approximation to
the $\widehat{WAP}$ statistic is a function of $\widehat{S}_{n}$, $\widehat{T%
}_{n}$, $\Sigma _{n}$, and $\widehat{D}_{Z}$. The critical values for the
WAP conditional tests and WAP-SU tests, respectively $\kappa _{n}\left(
t\right) $ and $\kappa _{n}\left( s,t\right) $, are taken under the
assumption that the $k$-dimensional vector $\widehat{S}_{n}$ has a standard
normal distribution (in practice, these critical values are also functions
of the consistent estimators $\widehat{\Sigma }_{n}$ and $\widehat{D}_{Z}$
as well, but we omit this dependence out of convenience). For example, for a
given weight density $w\left( \beta ,\pi \right) $, the critical function $%
\kappa _{n}\left( t\right) $ is simply the $1-\alpha $ quantile of (\ref%
{(WAP 1st Laplace)}) given $T=t$.

We now find the asymptotic distribution for the WAP\ tests under the WIV
asymptotics. We make the following assumption.

\bigskip

\noindent \textbf{Assumption WIV-FA}. (a) $\pi =C/n^{1/2}$ for some
non-stochastic vector $C$.

(b) $\beta $ is a fixed constant for all $n\geq 1.$

(c) $k$ is a fixed positive integer that does not depend on $n.$\smallskip

\bigskip

Under WIV, $\pi \left( \beta \right) $ is $o_{p}\left( 1\right) $ and the
WAP statistics behave the same as if the weights were simply $w\left( \beta
,0\right) $. As $n\rightarrow \infty $, the finite-sample critical value
functions $\kappa _{n}\left( t\right) $ and $\kappa _{n}\left( s,t\right) $
respectively converge to their asymptotic counterparts $\kappa _{\infty
}\left( t\right) $ and $\kappa _{\infty }\left( s,t\right) $, which are
based on (\ref{(WAP 1st Laplace)}) with $w\left( \beta ,\pi \left( \beta
\right) \right) $ replaced by $w\left( \beta ,0\right) $. We then obtain the
following convergence by the continuous mapping theorem and the joint
distribution%
\begin{eqnarray}
\left[ 
\begin{array}{c}
S_{\infty } \\ 
T_{\infty }%
\end{array}%
\right] &\sim &N\left( \left[ 
\begin{array}{c}
\left( \beta -\beta _{0}\right) C_{\beta _{0},\infty } \\ 
D_{\beta _{0},\infty }%
\end{array}%
\right] \left( D_{Z}\right) ^{1/2}C,I_{2k}\right) \text{, where}
\label{(S and T WIV-FA)} \\
C_{\beta _{0},\infty } &=&\left[ \left( b_{0}^{\prime }\otimes I_{k}\right)
\Sigma _{\infty }\left( b_{0}\otimes I_{k}\right) \right] ^{-1/2}\text{ and}
\notag \\
D_{\beta _{0},\infty } &=&\left[ \left( a_{0}^{\prime }\otimes I_{k}\right)
\Sigma _{\infty }^{-1}\left( a_{0}\otimes I_{k}\right) \right] ^{-1/2}\left(
a_{0}^{\prime }\otimes I_{k}\right) \Sigma _{\infty }^{-1}\left( a\otimes
I_{k}\right) .  \notag
\end{eqnarray}

\bigskip

\begin{theorem}
\label{Weak IV Thm} Under Assumptions \emph{W\emph{IV-FA} }and \emph{1-3}:%
\newline
\emph{(i)} $\left( \widehat{S}_{n},\widehat{T}_{n}\right) \rightarrow
_{d}\left( S_{\infty },T_{\infty }\right) ;$\newline
\emph{(ii)} $P\left( WAP\left( \widehat{S}_{n},\widehat{T}_{n}\right)
>\kappa _{n}\left( \widehat{T}_{n}\right) \right) \rightarrow P\left(
WAP\left( S_{\infty },T_{\infty }\right) >\kappa _{\infty }\left( T_{\infty
}\right) \right) ;$ and\newline
\emph{(iii)} $P\left( WAP\left( \widehat{S}_{n},\widehat{T}_{n}\right)
>\kappa _{n}\left( \widehat{S}_{n},\widehat{T}_{n}\right) \right)
\rightarrow P\left( WAP\left( S_{\infty },T_{\infty }\right) >\kappa
_{\infty }\left( S_{\infty },T_{\infty }\right) \right) .$
\end{theorem}

\bigskip

Both WAP conditional and WAP-SU tests have asymptotic null rejection
probabilities being equal to $\alpha $. The asymptotic power of the WAP
tests has a complicated form under WIV asymptotics. We can, of course, rely
on numerical simulations to compare their performance with other available
tests. In Section \ref{Application Sec}, we present power plots for testing
the intertemporal elasticity of substitution based on the designs of %
\citet{Yogo04}.

For strong instruments with local alternatives (SIV-LA), we consider the
Pitman drift where $\beta $ is local to the null value $\beta _{0}$ as $%
n\rightarrow \infty $.

\bigskip

\noindent \textbf{Assumption SIV-LA. }(a)\textbf{\ }$\beta =\beta
_{0}+B/n^{1/2}$ for some constant $B\in \mathbb{R}.$

(b) $\pi $ is a fixed non-zero $k$-vector for all $n\geq 1.$

(c) $k$ is a fixed positive integer that does not depend on $n$.$\medskip $

\bigskip

Under the SIV-LA asymptotics, the WAP statistics are shown to be increasing
transformations of the $LR$ statistic. This result is general and holds for
any prior which satisfies Assumption 3.

\bigskip

\begin{theorem}
\label{Asy Eff Strong IV Thm} Suppose Assumptions \emph{\emph{SIV-LA} }and 
\emph{1-3 }hold. The long-run variance $\Sigma _{\infty }$ is known, or
unknown but consistently estimable by $\widehat{\Sigma }_{n}$. Then the WAP
similar and WAP-SU tests are asymptotically equivalent to the LR test given
in (\ref{(LR stat)}).
\end{theorem}

\textbf{Comment.} \textbf{1. }In the proof, we apply the Laplace
approximation twice, first with respect to the integral for $\pi $ and then
for $\beta $. For the MM1 and MM2 statistics, we can alternatively find a
simple expression after integrating out the prior for the instruments'
coefficients with $\sigma ^{2}$ or $\zeta $ growing at rate $n$ and then
applying the Laplace approximation for $\beta $. Both approaches coincide.

\textbf{2. }The SIV-LA behavior of the ECS (HAC-IV) test appears to be just
a special case of our theory using Laplace approximations.

\textbf{3. }For higher-order expansions, we can use Watson's lemma; for
references, we recommend \citet{Olver97} for deterministic functions and %
\citeauthor{OnatskiMoreiraHallin14a} (\citeyear{OnatskiMoreiraHallin14a}, %
\citeyear{OnatskiMoreiraHallin14b}) for random functions.

\textbf{4.} Because $T_{n}/n^{1/2}\rightarrow _{p}D_{\beta
_{0}}D_{Z}^{1/2}\pi $ under SIV-LA, $\left\Vert T_{n}\right\Vert $ diverges
to infinity w.p.1 (with probability approaching one). The critical value
functions for both the WAP conditional and WAP-SU tests collapse then to the 
$1-\alpha $ asymptotic (unconditional) quantile. As a result, the WAP
conditional and WAP-SU tests are asymptotically similar and efficient under
the SIV asymptotics.

\bigskip

The null rejection probability of WAP tests is $\alpha $ under WIV and SIV
asymptotics. Pointwise convergence of the null rejection probability, of
course, does not necessarily imply the size is asymptotically $\alpha $ (in
a uniform sense). \citet[p. 1037]{Moreira03} suggests to use \citet{Parzen54}
and \citet{Andrews86} to assure size is uniformly controlled. A series of
papers, including \citet{AndrewsChengGuggenberger11} and %
\citet{AndrewsGuggenberger14a}, develop several powerful methods to check
uniform size control and have been applied to many econometric models; see %
\citet{AndrewsGuggenberger10}, \citet{AndrewsGuggenberger14a}, and %
\citet{MillsMoreiraVilela14b}, among others. Conceivably, we can apply those
methods to the WAP statistics coupled with the critical value functions $%
\kappa _{n}\left( t\right) $ and $\kappa _{n}\left( s,t\right) $. This line
of research will be considered in a separate paper.

We can also analyze the WAP tests under strong instruments with fixed
alternatives (SIV-FA). We follow \citet{MillsMoreiraVilela14} and make the
following assumption.

\bigskip

\noindent \textbf{Assumption SIV-FA. }(a)\textbf{\ }$\beta =\beta _{0}+B$
for some nonzero $B\in \mathbb{R}.$

(b) $\pi $ is a fixed non-zero $k$-vector for all $n\geq 1.$

(c) $k$ is a fixed positive integer that does not depend on $n$.$\medskip $

\bigskip

It is natural to expect that the power converges to one if the parameter $%
\beta $ is fixed. However, not all tests have this property even in the IV
model with homoskedastic errors; see \citet{AndrewsMoreiraStock04} and %
\citet{MillsMoreiraVilela14} for examples. Hence, it is important to
establish consistency for the WAP tests.

If the parameter $\beta $ is fixed, the WAP statistics are proportional to
the exponential of $LR$. Because $LR/n$ converges to a non-zero constant,
the WAP tests are consistent. The next theorem formalizes this result.

\bigskip

\begin{theorem}
\label{Consistency Thm} Suppose Assumptions \emph{\emph{SIV-FA} }and \emph{%
1-3 }hold. The long-run variance $\Sigma _{\infty }$ is known, or unknown
but consistently estimable by $\widehat{\Sigma }_{n}$. Then the following
hold:\newline
\emph{(i)} $2.\left( \log \widehat{WAP}\right) /n=\widehat{LR}/n+o_{p}\left(
1\right) ;$ and\newline
\emph{(ii)} $\widehat{LR}/n=LR/n+o_{p}\left( 1\right) \rightarrow \gamma >0.$
\end{theorem}

\textbf{Comment:} If $D_{\beta }\neq 0$, the functions $\kappa_{n} \left(
t\right) $ and $\kappa_{n} \left( s,t\right) $ converge to a constant
obtained under SIV-FA. If $D_{\beta }=0$, the critical functions do not
converge. However, they are bounded, and so WAP tests are consistent.

\section{Power Comparison \label{Application Sec}}

In this section, we follow I. \citet{Andrews15} who calibrates designs for
power comparison based on the work of \citet{Yogo04} on the elasticity of
intertemporal substitution in eleven developed countries.

\citet{Yogo04} tests the effect of interest rates on the level of aggregate
demand in an IV model. He considers a linear regression in which asset
return affects consumption growth, and the reverse form of this regression.
In both equations, the endogenous variable (consumption or asset return) can
be correlated with the error (innovation). To remedy this problem, he
chooses four instruments: lagged values of nominal interest rate, inflation,
consumption growth, and log dividend-price ratio.

I. \citet{Andrews15} selects the real interest rate (\emph{rf} in %
\citeauthor{Yogo04}'s (\citeyear{Yogo04}) notation) as the endogenous
variable. Several tests perform well in his design, including MM2-SU,
PI-CLC, and (WAP similar) ECS tests. In fact, only in a few countries do
these tests have slightly different performance; see Section 7.2.1 of I. %
\citet{Andrews15}. The difficulty in assessing the relative performance of
each test arises because the instruments are not particularly weak in this
design. Indeed, the first-stage F-statistic reported by \citet{Yogo04} (see
his Table I) is below 10 in only four countries (Japan, Switzerland, United
Kingdom, and the United States). We instead join \citet{deCastro15} in
choosing the real stock return (\emph{re} in \citeauthor{Yogo04}'s (%
\citeyear{Yogo04}) notation) as the endogenous variable. The instruments are
considerably weaker in this design: the F-statistic is smaller than 4.18 in
all countries, and always less than the F-statistic for interest rate. Our
decision to use stock returns aims to highlight the differences between the
tests proposed for the HAC-IV model. Apart from using stock returns instead
of interest rates, our design is akin to that of I. \citet{Andrews15}. We
use the Newey-West estimator with three lags, and the resulting power curves
are based on 5,000 Monte Carlo simulations. In parallel to our asymptotic
theory, we choose the ratio of the tuning parameters $\sigma ^{2}$ and $%
\zeta $ to the sample size to be one-tenth for the MM1 and MM2 statistics,
respectively.

Figure 3 plots power curves for the two-sided power envelope, Anderson-Rubin
(AR), score (LM), WAP similar MM1, WAP similar MM2, and ECS (HAC-IV) tests.
Although the AR and LM tests are unbiased, the MM1, MM2, and ECS tests
perform unreliably. To illustrate the problem, we mention three countries.
For Australia, the MM1 and ECS tests have low power for parts of the
parameter space, while the MM2 test behaves more like a two-sided test. For
France, the ECS test performs well, while both MM1 and MM2 tests can have
low power. For the USA, the ECS test has power near zero and behaves more as
a one-sided test while the MM1 and MM2 tests are nearly unbiased. In some
countries, these three tests have power even lower than the Anderson-Rubin
test (e.g., the ECS test for Germany and Italy).

\begin{figure}[tbh]
\caption{Power Comparison (WAP similar tests)}%
\begin{subfigure}[b]{.5\textwidth}
  \centering
  \includegraphics[trim = 45mm 80mm 45mm 80mm, clip, width=5.5cm]{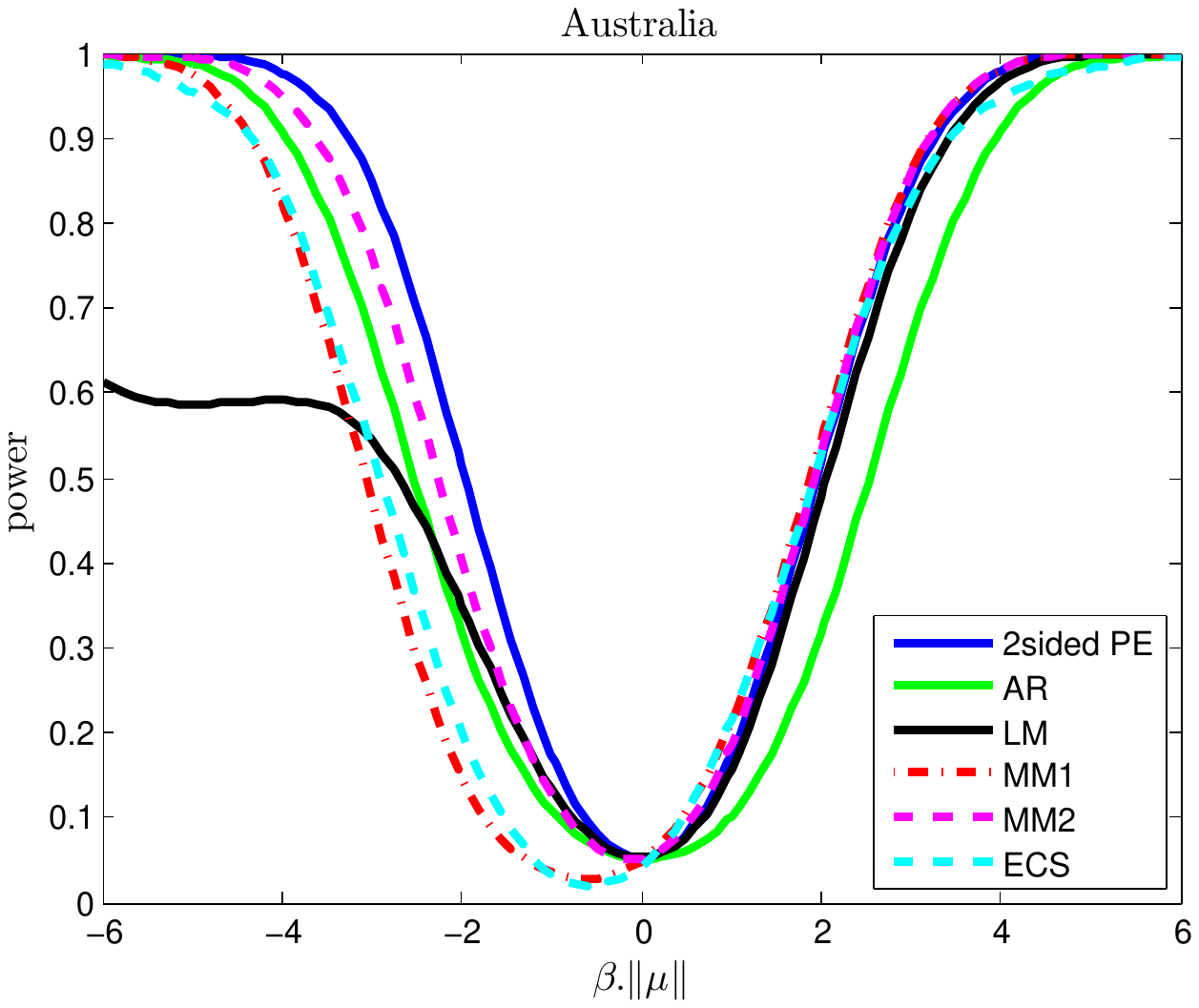}
\end{subfigure}%
\begin{subfigure}[b]{.5\textwidth}
  \centering
  \includegraphics[trim = 45mm 80mm 45mm 80mm, clip, width=5.5cm]{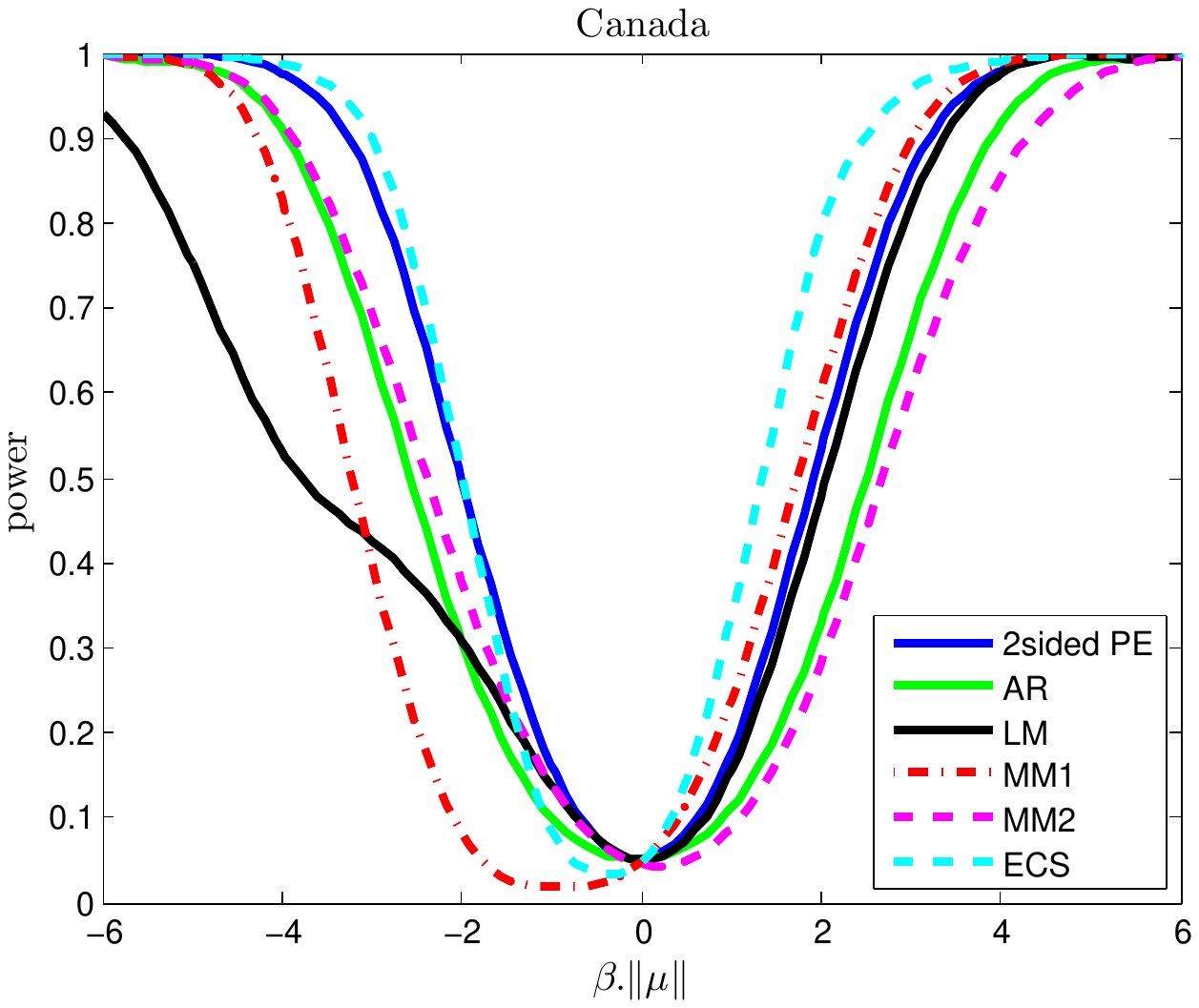}
\end{subfigure}
\begin{subfigure}[b]{.5\textwidth}
  \centering
  \includegraphics[trim = 45mm 80mm 45mm 80mm, clip, width=5.5cm]{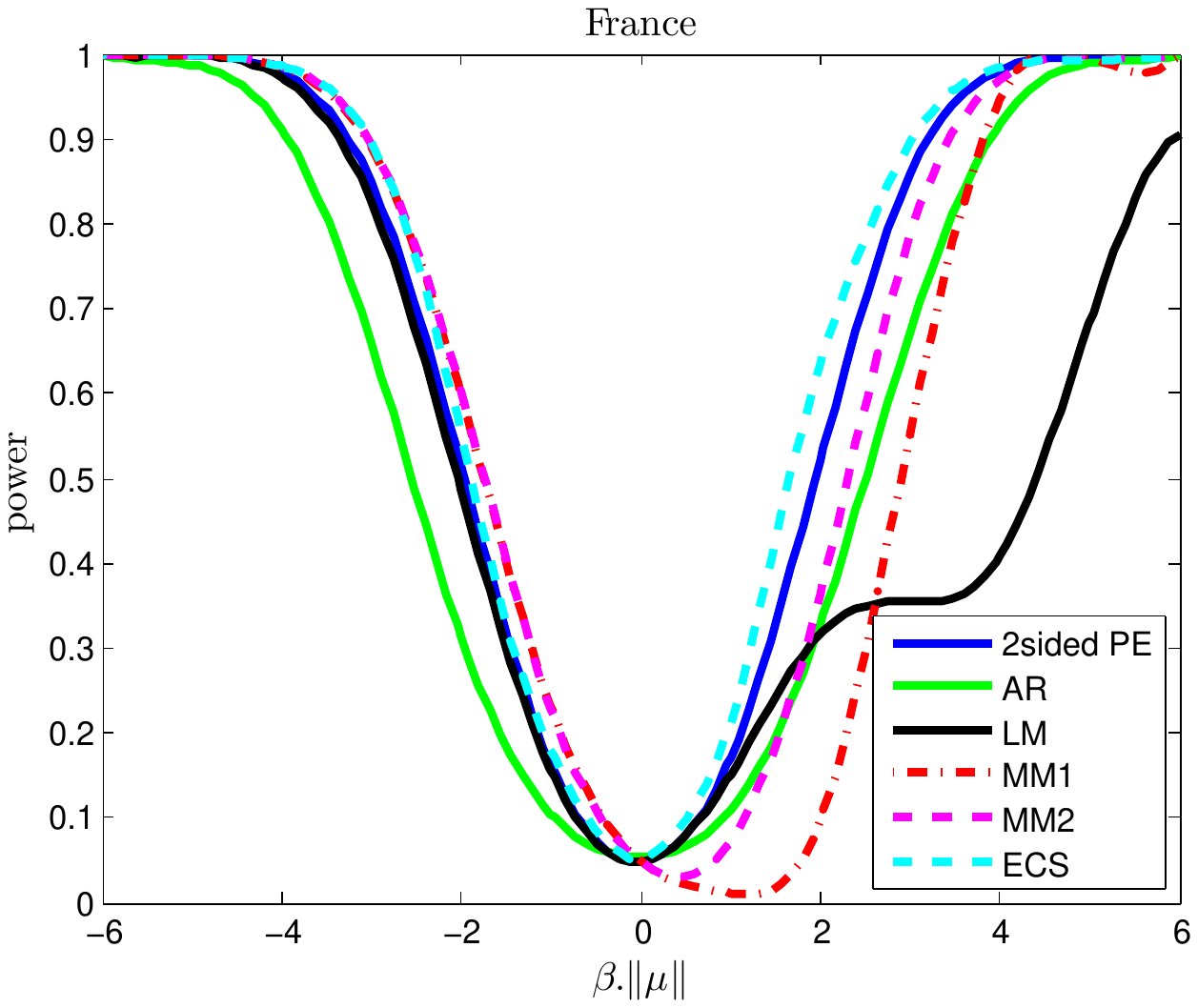}
\end{subfigure}%
\begin{subfigure}[b]{.5\textwidth}
  \centering
  \includegraphics[trim = 45mm 80mm 45mm 80mm, clip, width=5.5cm]{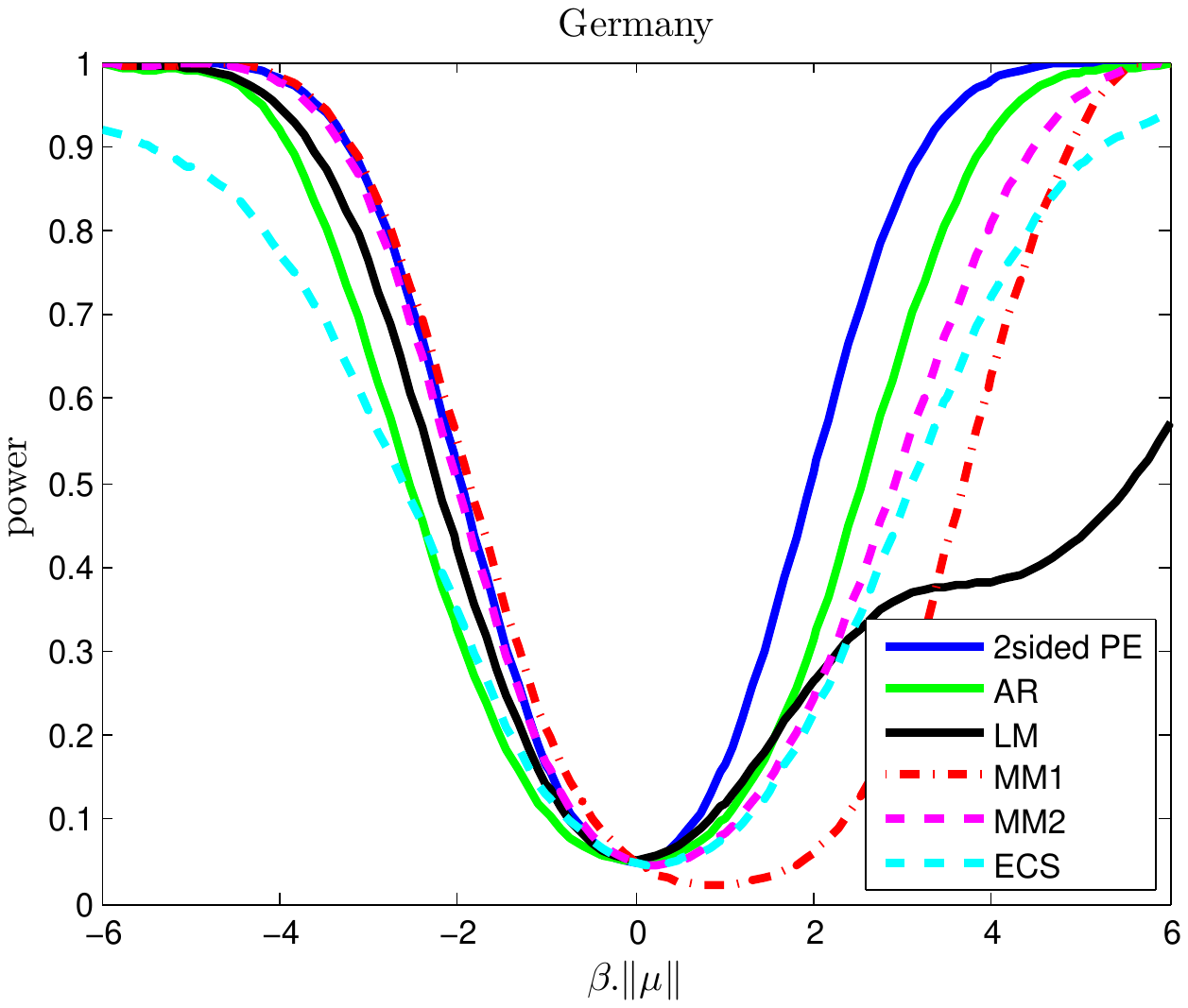}
\end{subfigure}
\begin{subfigure}[b]{.5\textwidth}
  \centering
  \includegraphics[trim = 45mm 80mm 45mm 80mm, clip, width=5.5cm]{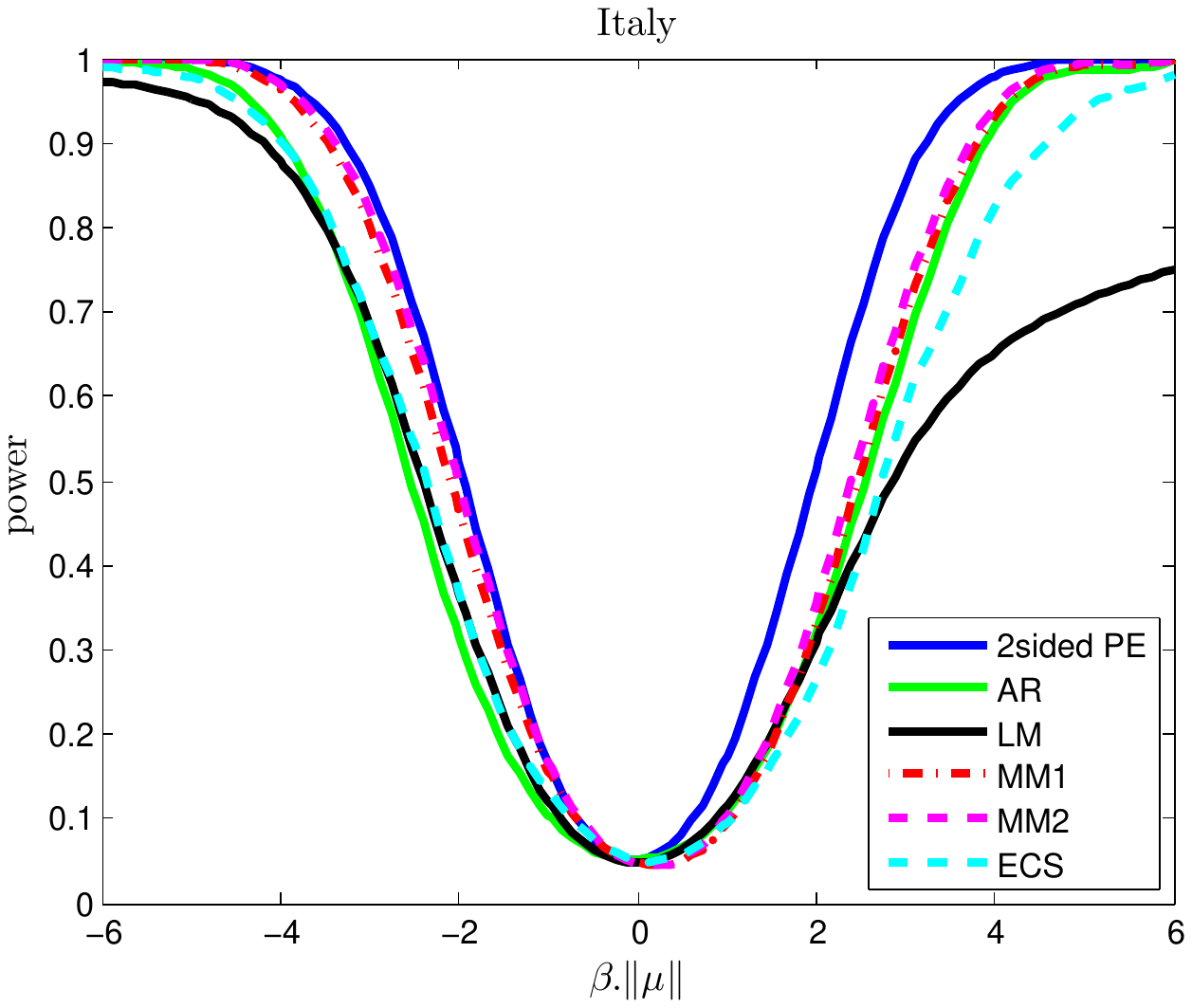}
\end{subfigure}%
\begin{subfigure}[b]{.5\textwidth}
  \centering
  \includegraphics[trim = 45mm 80mm 45mm 80mm, clip, width=5.5cm]{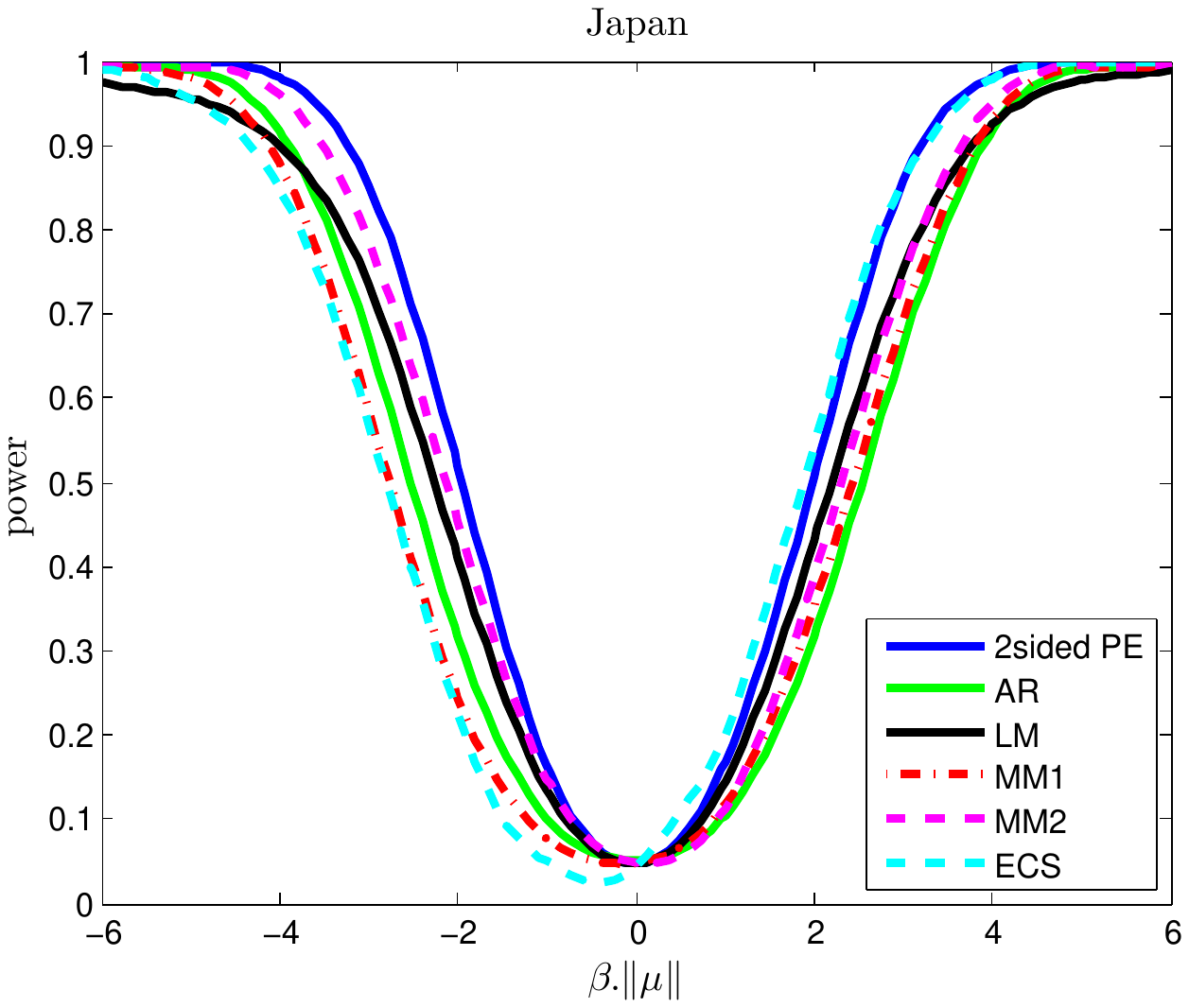}
\end{subfigure}
\end{figure}

\clearpage

\begin{figure}[tbh]
\ContinuedFloat  
\begin{subfigure}[b]{.5\textwidth}
  \centering
  \includegraphics[trim = 45mm 80mm 45mm 80mm, clip, width=5.5cm]{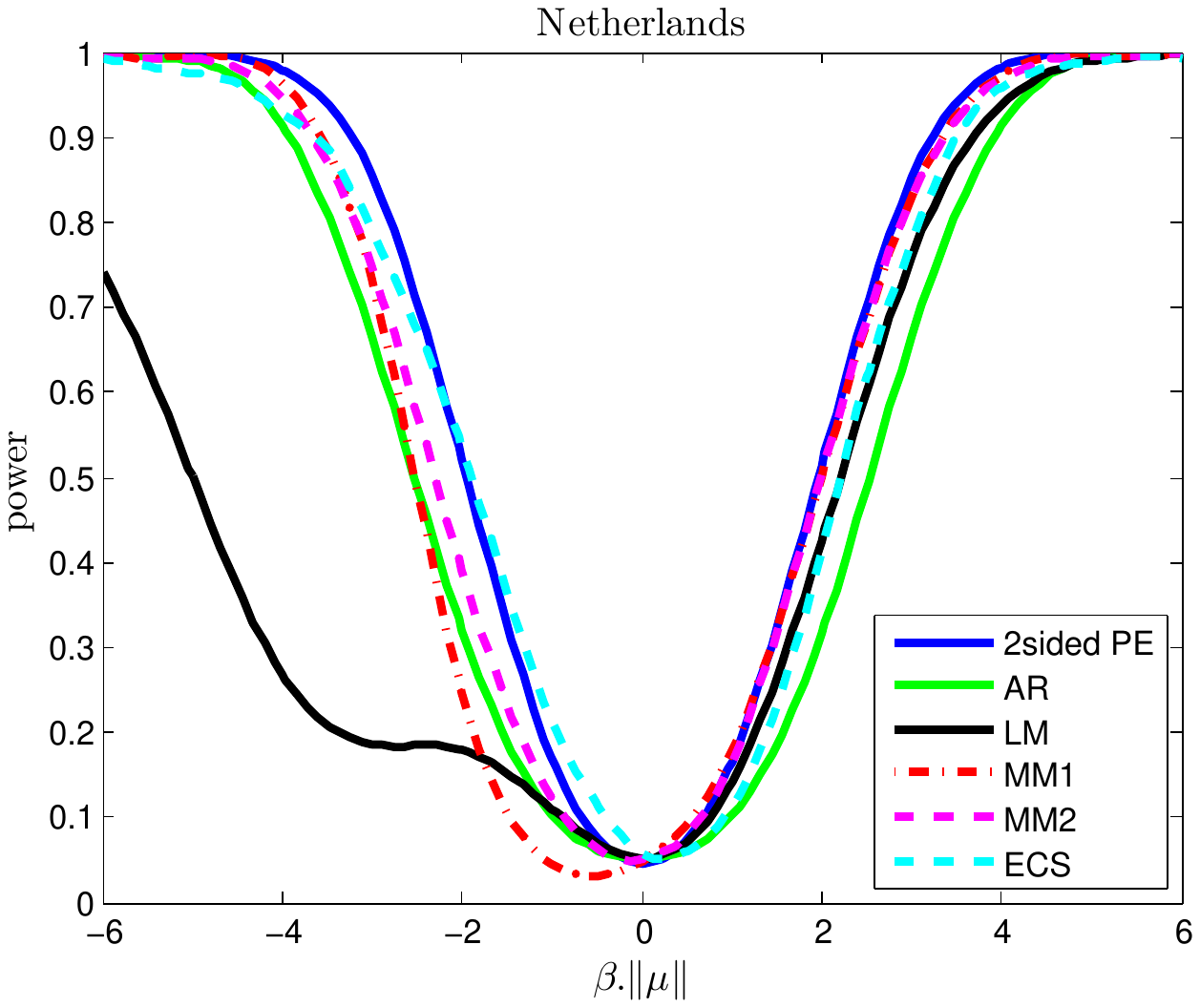}
\end{subfigure}%
\begin{subfigure}[b]{.5\textwidth}
  \centering
  \includegraphics[trim = 45mm 80mm 45mm 80mm, clip, width=5.5cm]{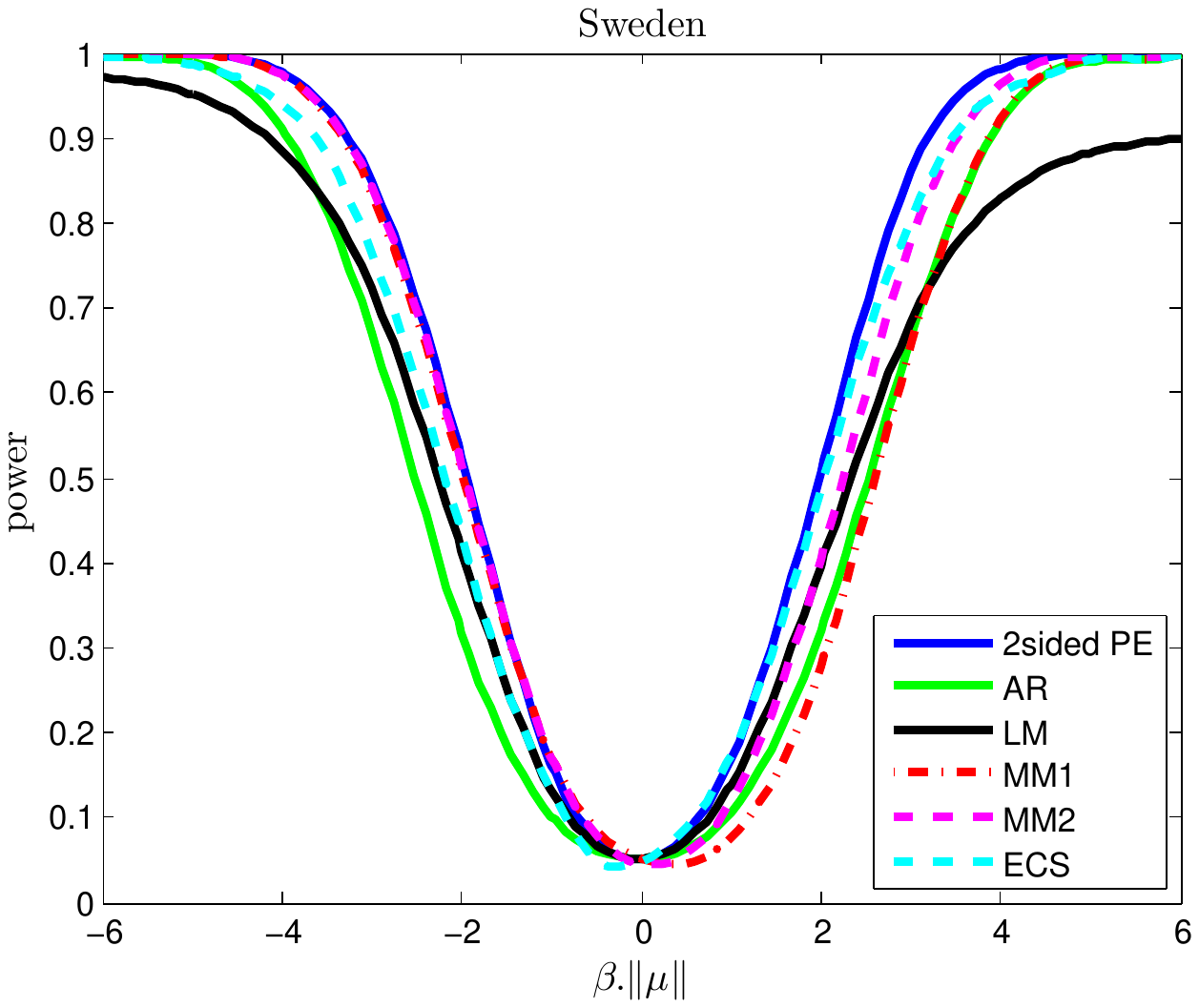}
\end{subfigure}
\begin{subfigure}[b]{.5\textwidth}
  \centering
  \includegraphics[trim = 45mm 80mm 45mm 80mm, clip, width=5.5cm]{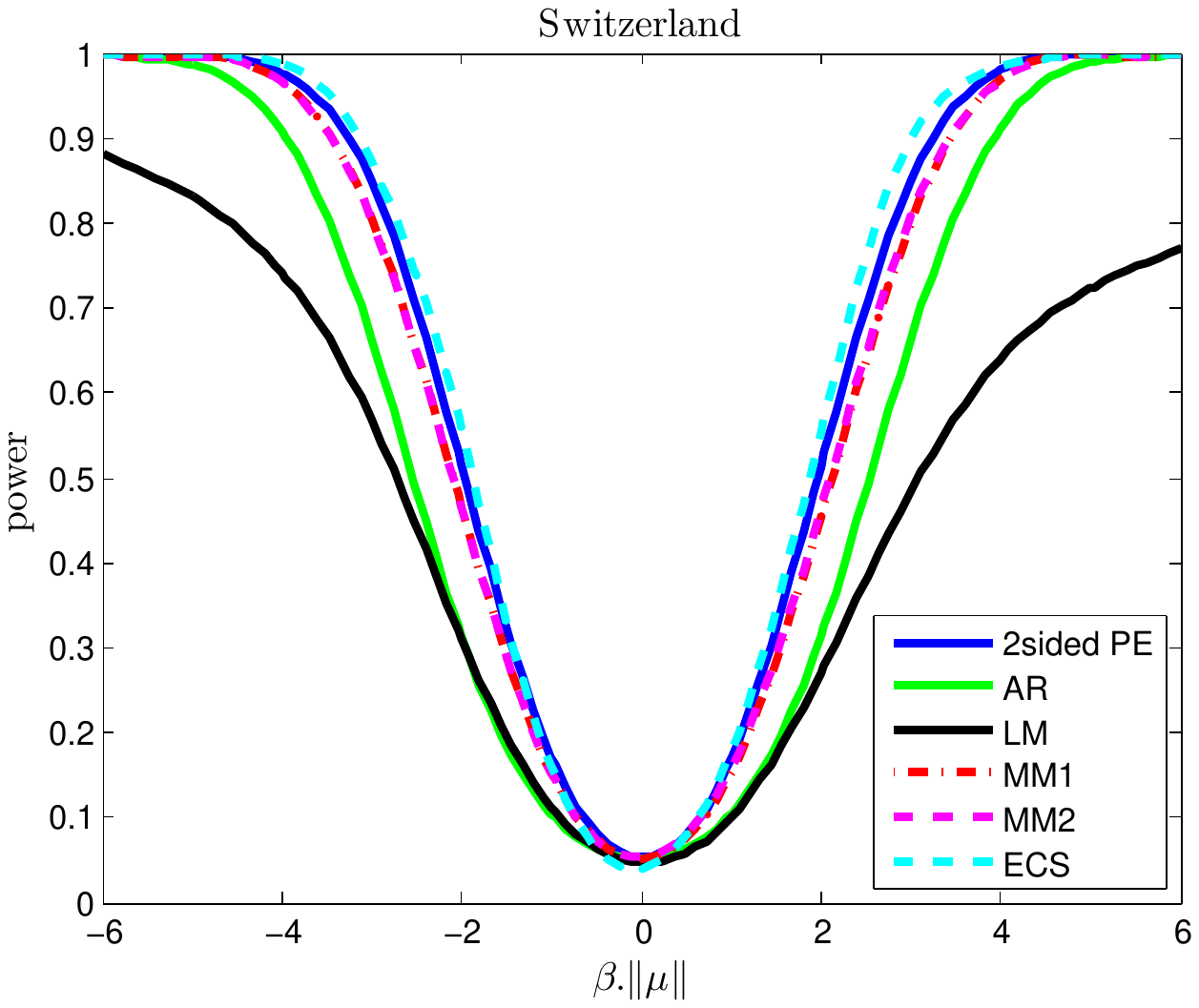}
\end{subfigure}%
\begin{subfigure}[b]{.5\textwidth}
  \centering
  \includegraphics[trim = 45mm 80mm 45mm 80mm, clip, width=5.5cm]{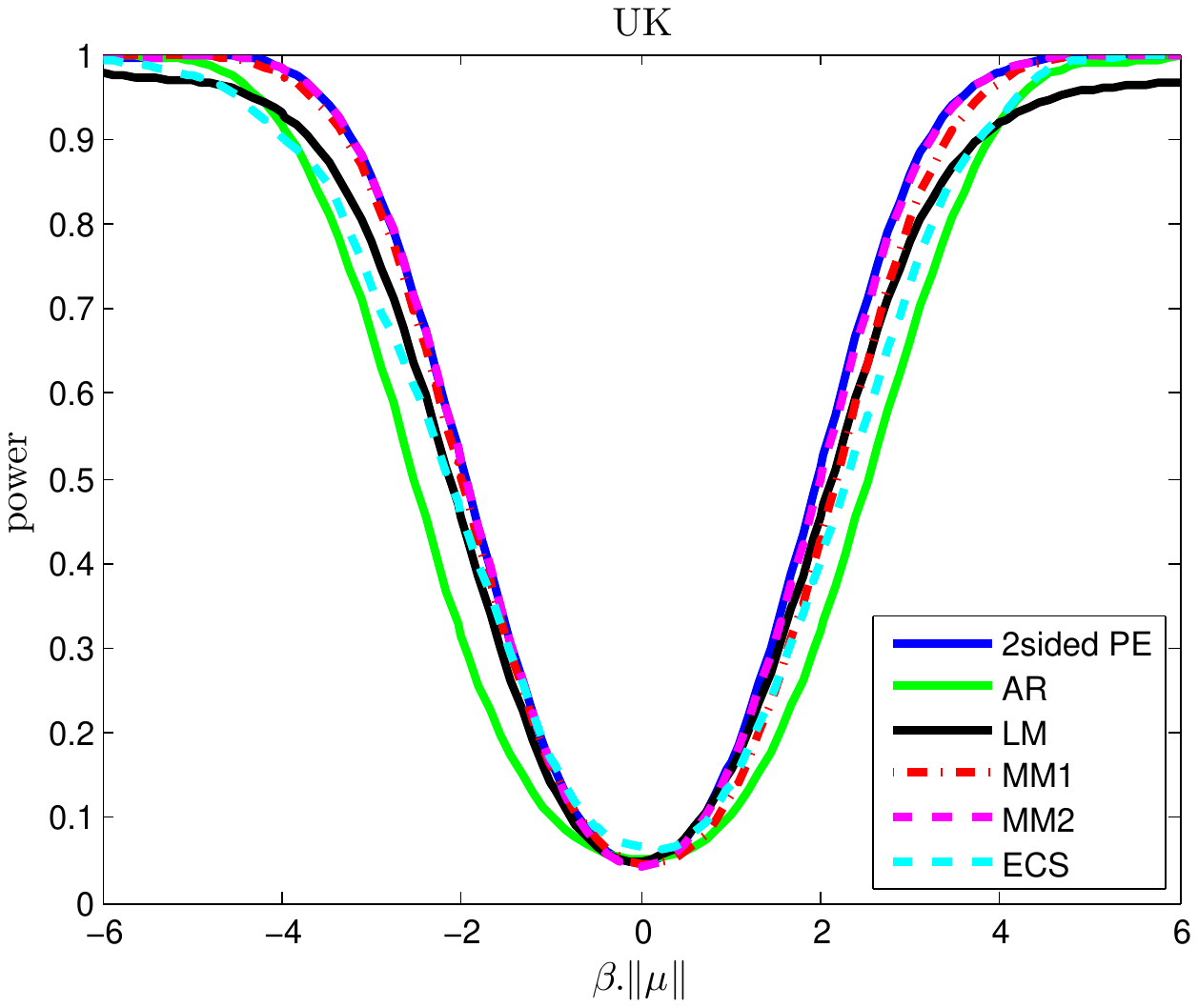}
\end{subfigure}
\begin{subfigure}[b]{.5\textwidth}
  \centering
  \includegraphics[trim = 45mm 80mm 45mm 80mm, clip, width=5.5cm]{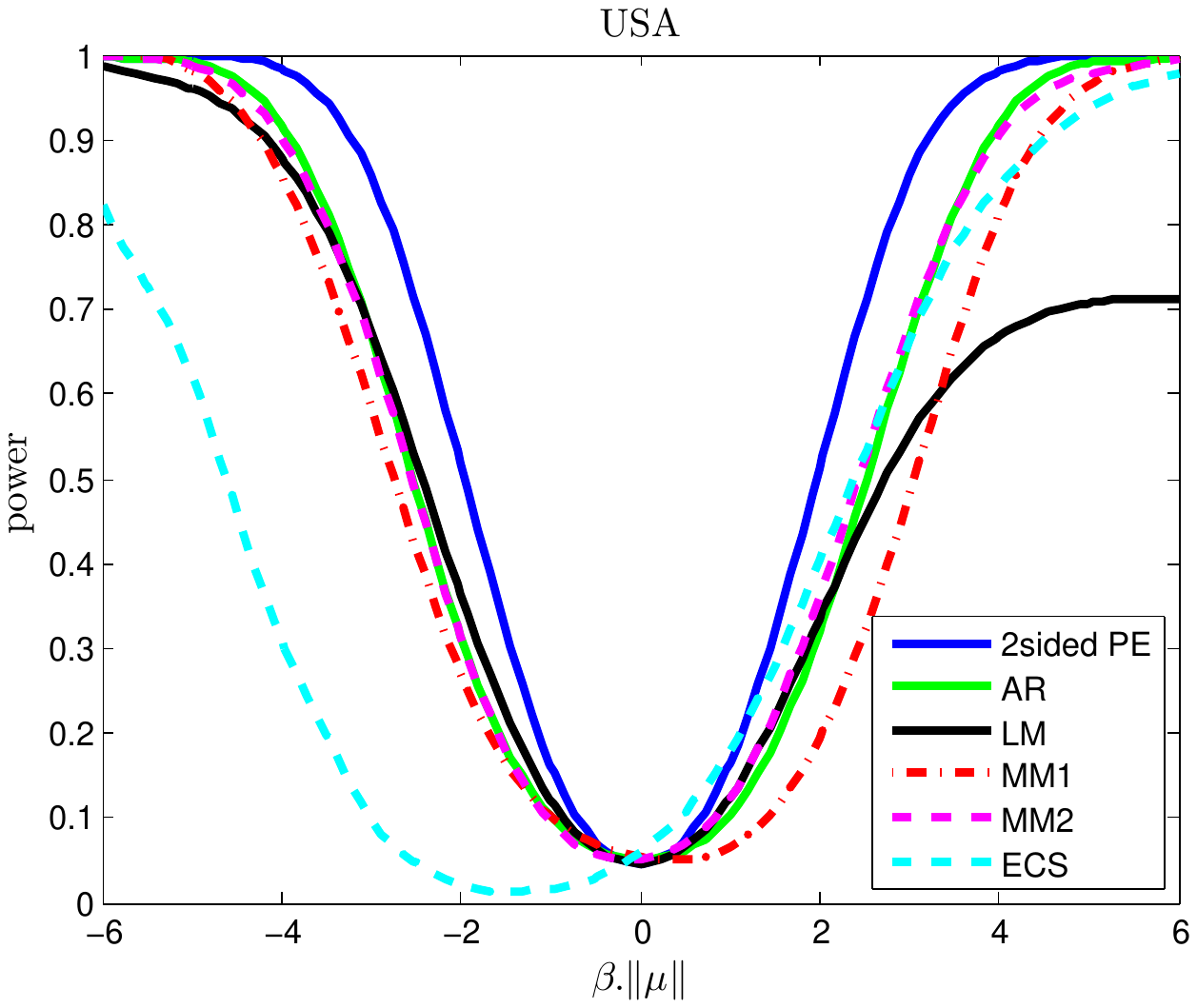}
\end{subfigure}
\end{figure}

We then compare power among two-sided tests which have arguably better
performance. Figure 4 plots power curves for the two-sided power envelope,
MM1-SU, MM2-SU, CQLR, CQLR-kron, and PI-CLC tests. All tests are adequate
for two-sided hypothesis testing. The PI-CLC and CQLR-kron test show some
improvements over the CQLR test for some, but not all, countries. The MM1-SU
test behaves near the MM2-SU test for several countries, but it has
considerably lower power for Japan and the United States\footnote{%
Conceivably, this power loss can be due to numerical integration over the
whole real line. Power may be improved by transforming the parameter $\beta $
to the quantity $\theta =\tan ^{-1}\left( d_{\beta }/c_{\beta }\right) $.
This improvement is left for future work.}. The MM2-SU test outperforms
these tests and when it occasionally has less power, the power loss is
small. This application based on real data supports our theoretical
contribution and the use of the MM2-SU test in practice.

\begin{figure}[tbh]
\caption{Power Comparison (two-sided tests)}%
\begin{subfigure}[b]{.5\textwidth}
  \centering
  \includegraphics[trim = 45mm 80mm 45mm 80mm, clip, width=5.5cm]{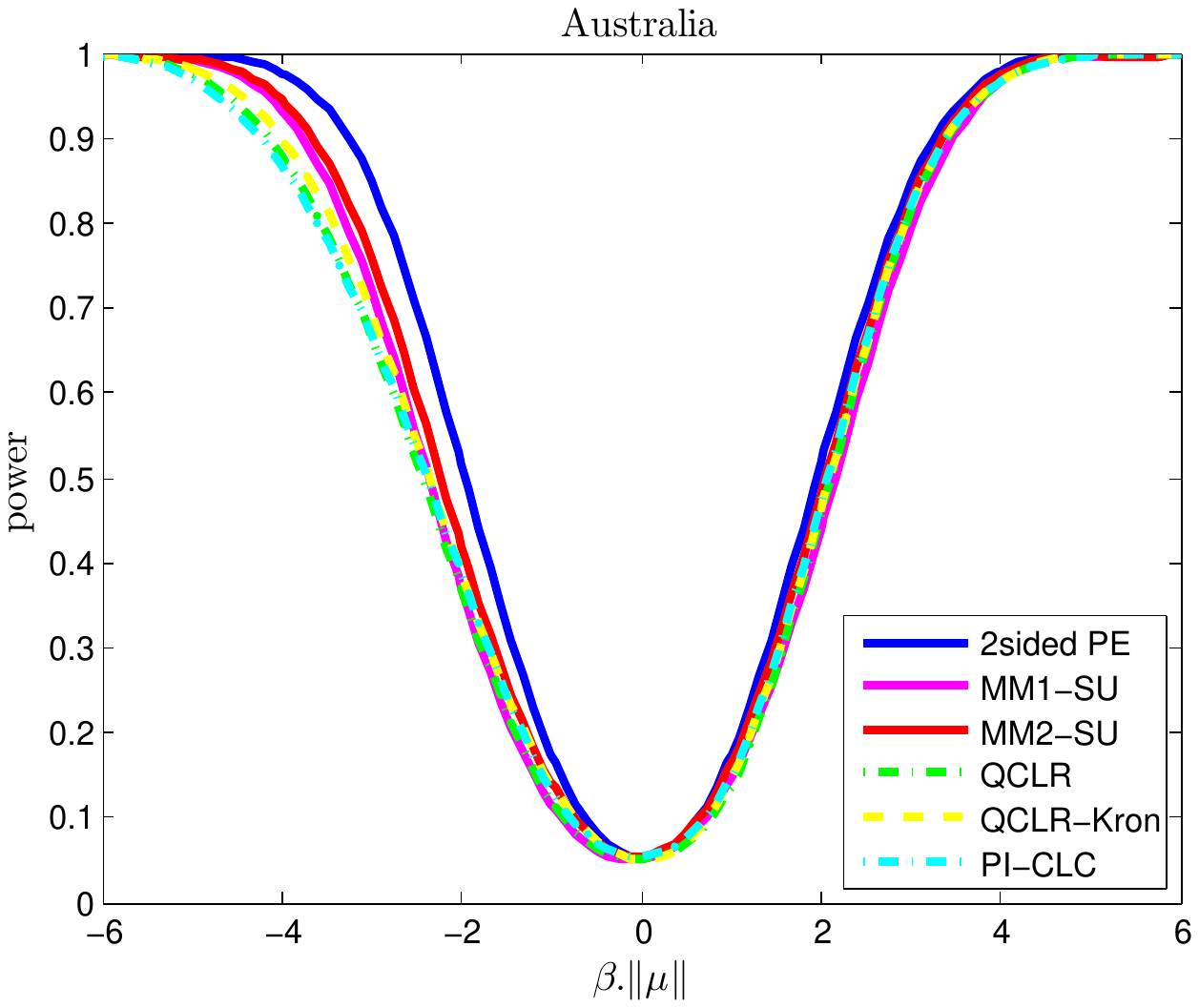}
\end{subfigure}%
\begin{subfigure}[b]{.5\textwidth}
  \centering
  \includegraphics[trim = 45mm 80mm 45mm 80mm, clip, width=5.5cm]{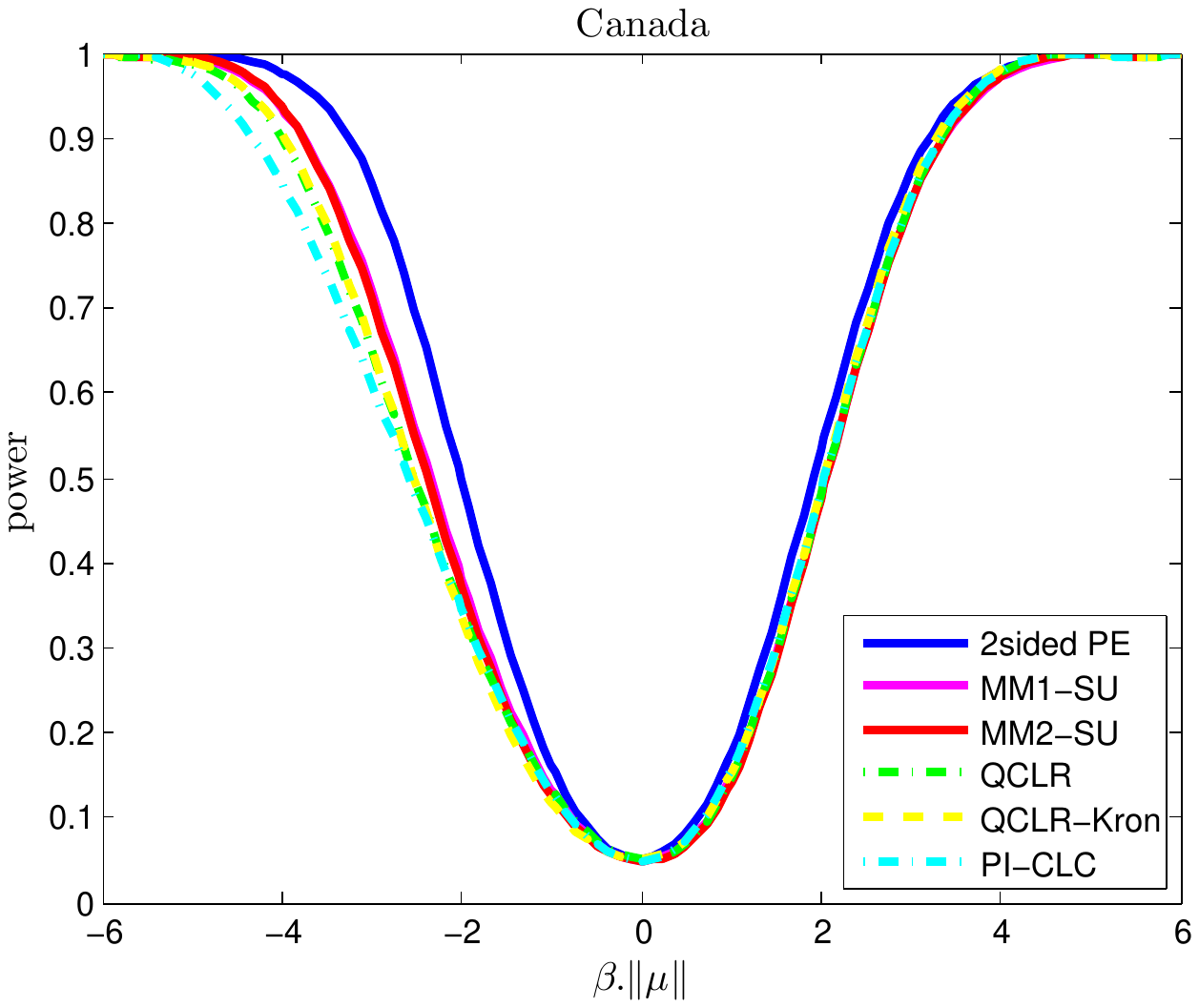}
\end{subfigure}
\begin{subfigure}[b]{.5\textwidth}
  \centering
  \includegraphics[trim = 45mm 80mm 45mm 80mm, clip, width=5.5cm]{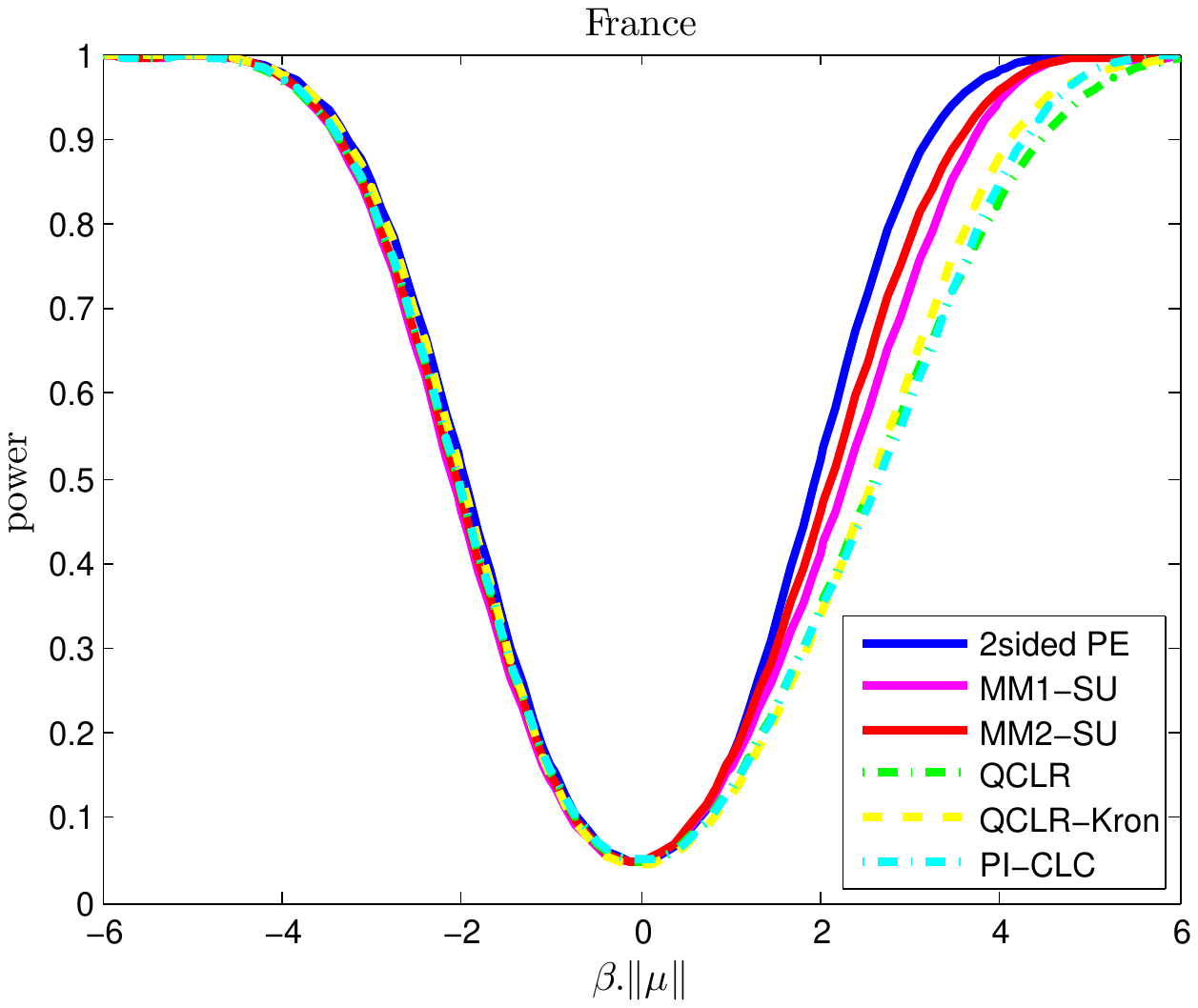}
\end{subfigure}%
\begin{subfigure}[b]{.5\textwidth}
  \centering
  \includegraphics[trim = 45mm 80mm 45mm 80mm, clip, width=5.5cm]{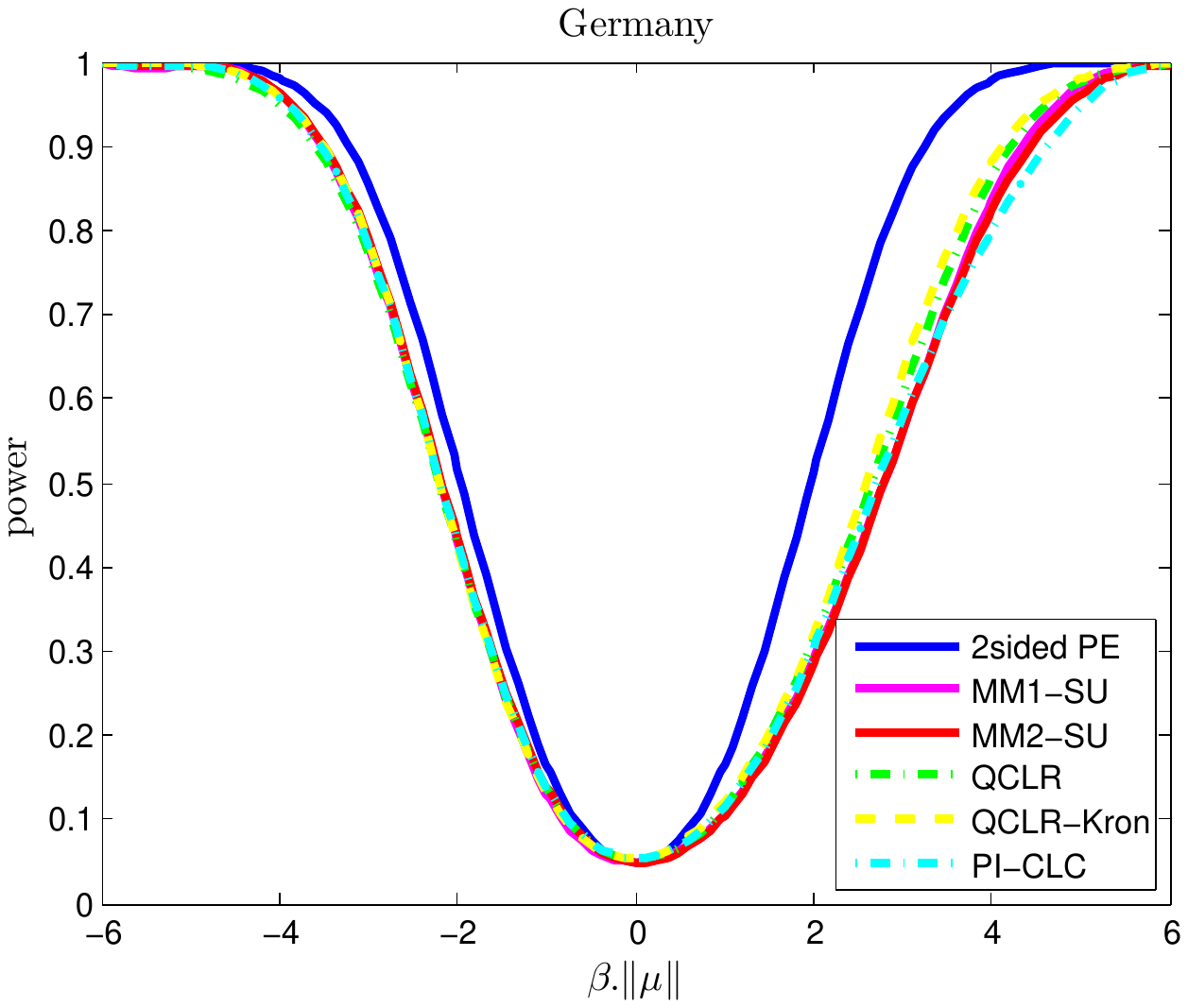}
\end{subfigure}
\begin{subfigure}[b]{.5\textwidth}
  \centering
  \includegraphics[trim = 45mm 80mm 45mm 80mm, clip, width=5.5cm]{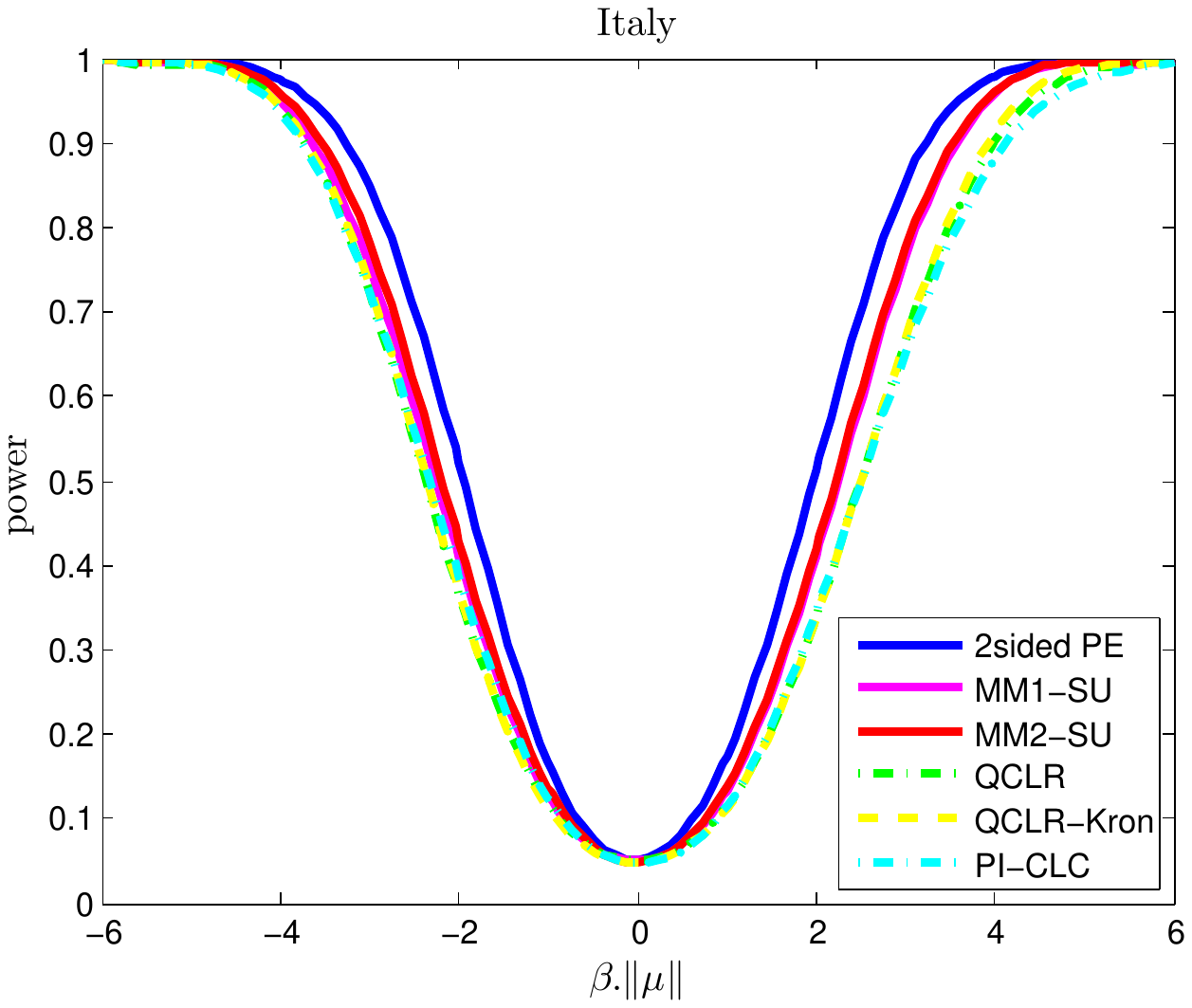}
\end{subfigure}%
\begin{subfigure}[b]{.5\textwidth}
  \centering
  \includegraphics[trim = 45mm 80mm 45mm 80mm, clip, width=5.5cm]{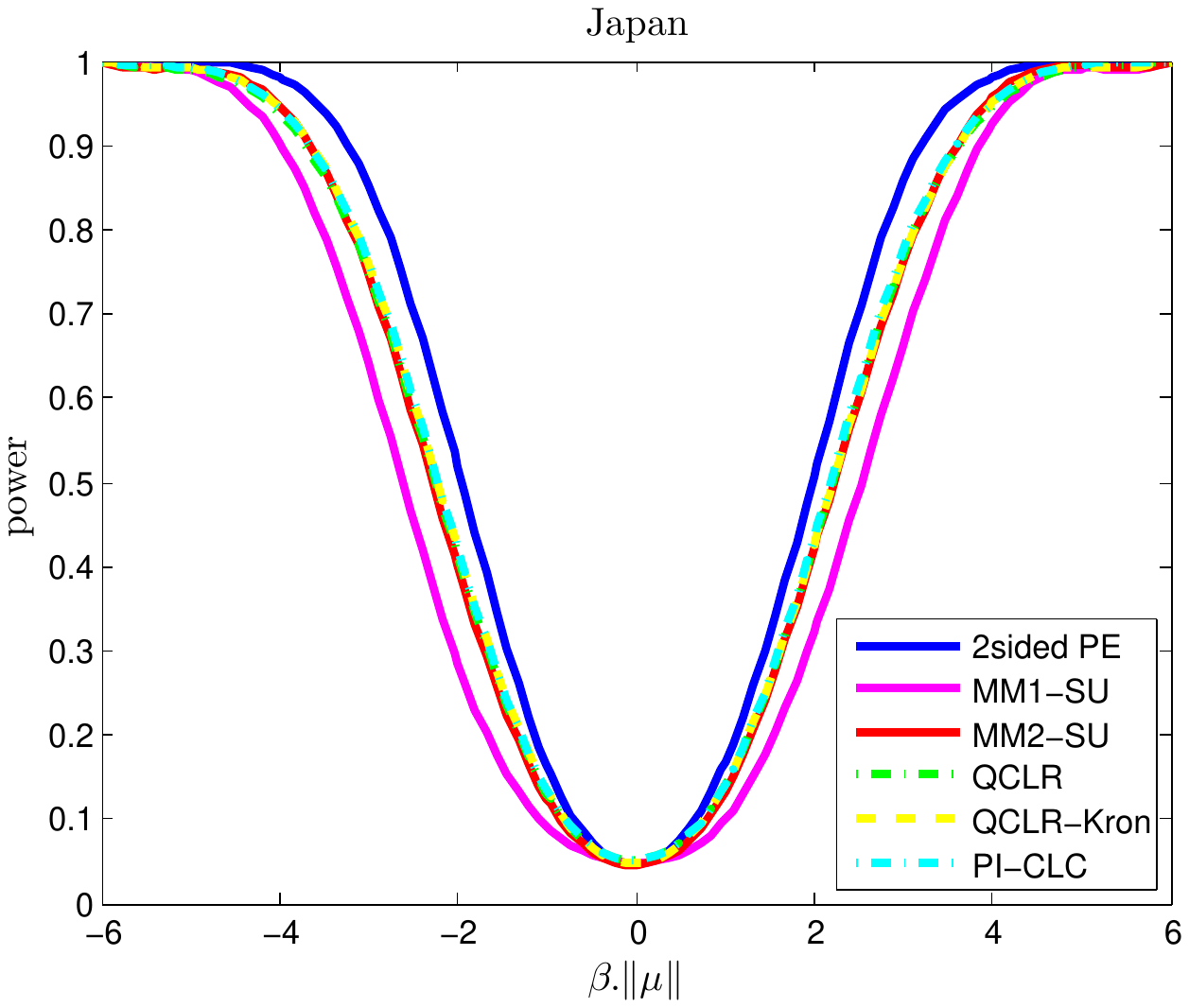}
\end{subfigure}
\end{figure}

\clearpage

\begin{figure}[tbh]
\ContinuedFloat  
\begin{subfigure}[b]{.5\textwidth}
  \centering
  \includegraphics[trim = 45mm 80mm 45mm 80mm, clip, width=5.5cm]{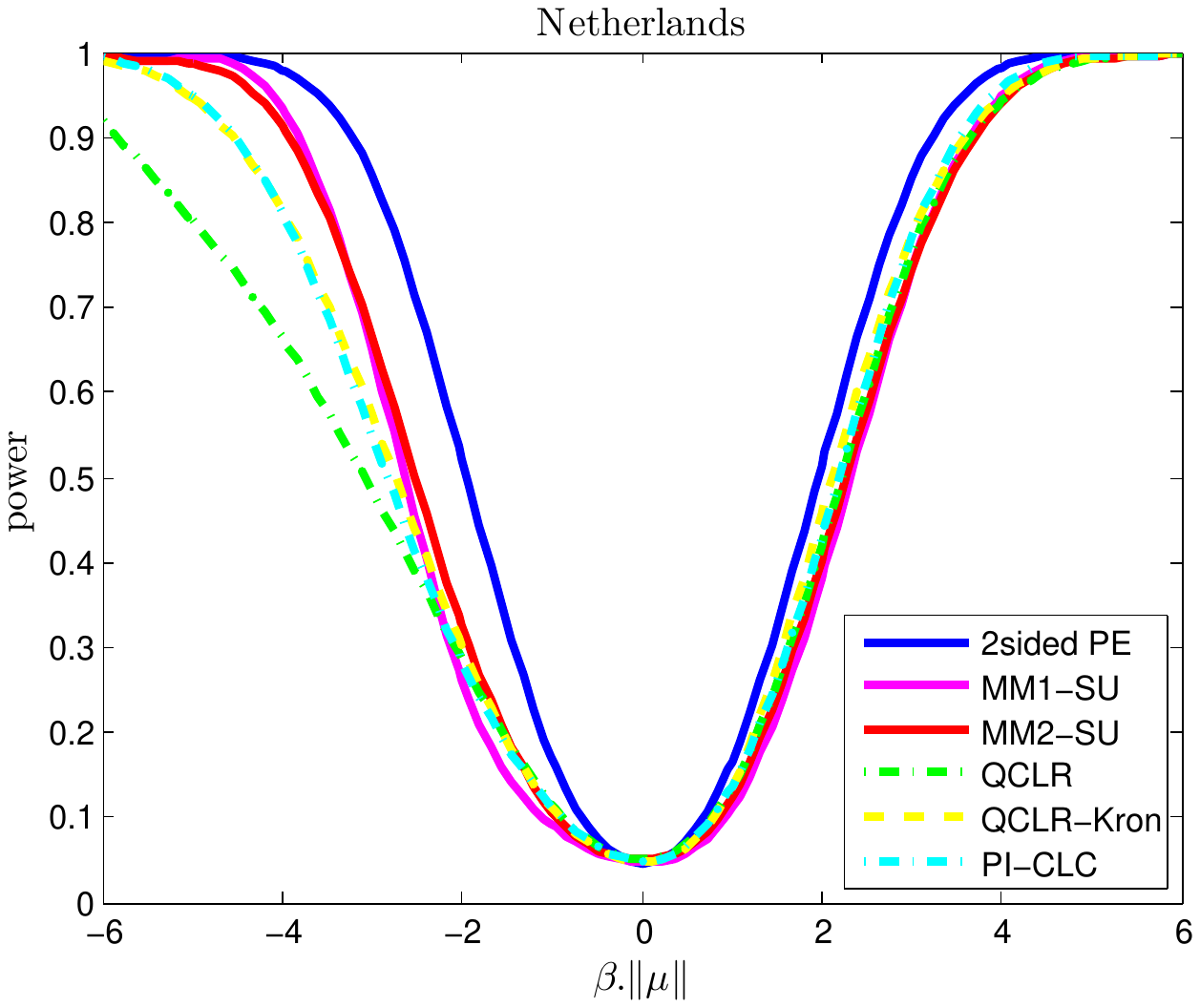}
\end{subfigure}%
\begin{subfigure}[b]{.5\textwidth}
  \centering
  \includegraphics[trim = 45mm 80mm 45mm 80mm, clip, width=5.5cm]{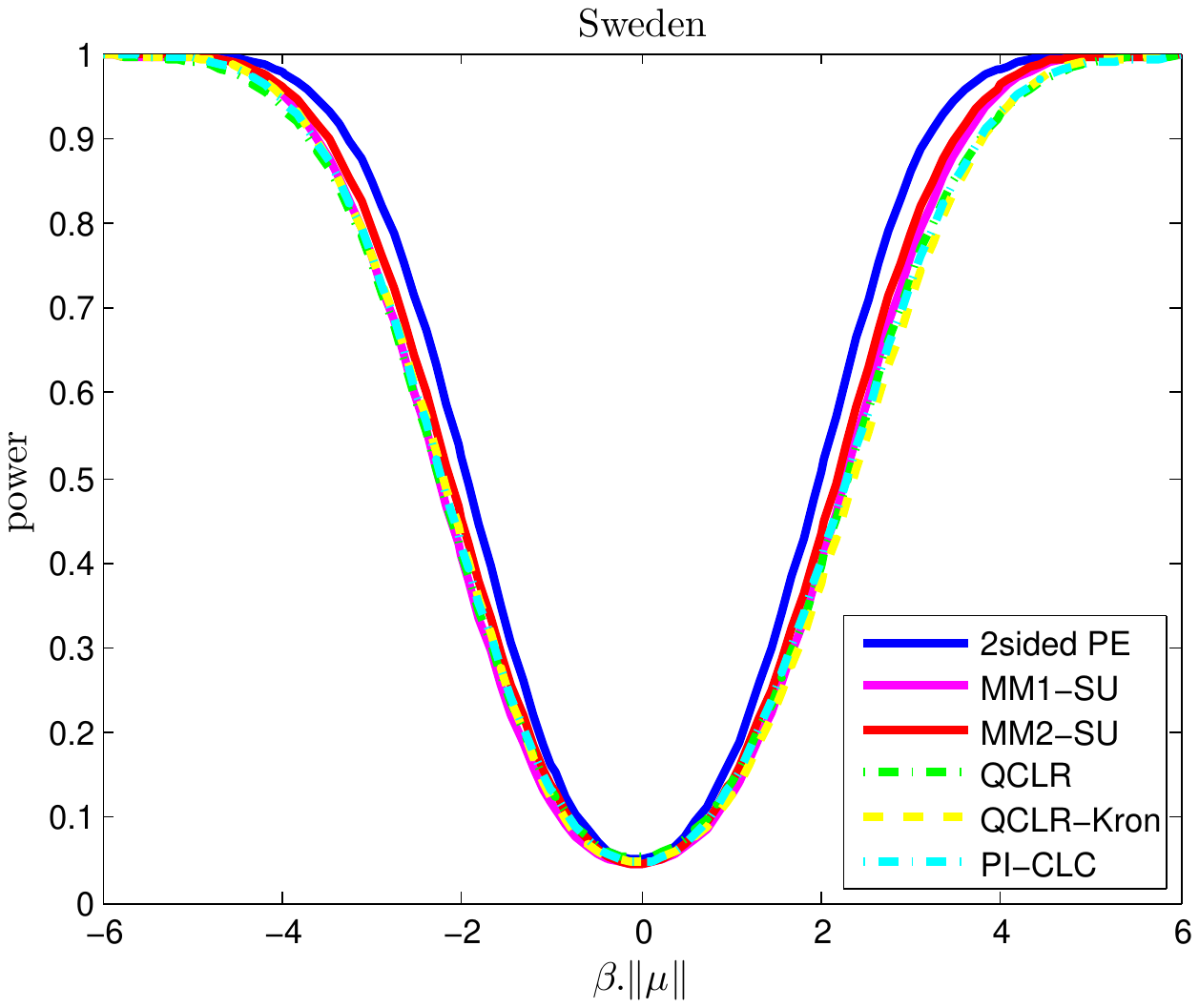}
\end{subfigure}
\begin{subfigure}[b]{.5\textwidth}
  \centering
  \includegraphics[trim = 45mm 80mm 45mm 80mm, clip, width=5.5cm]{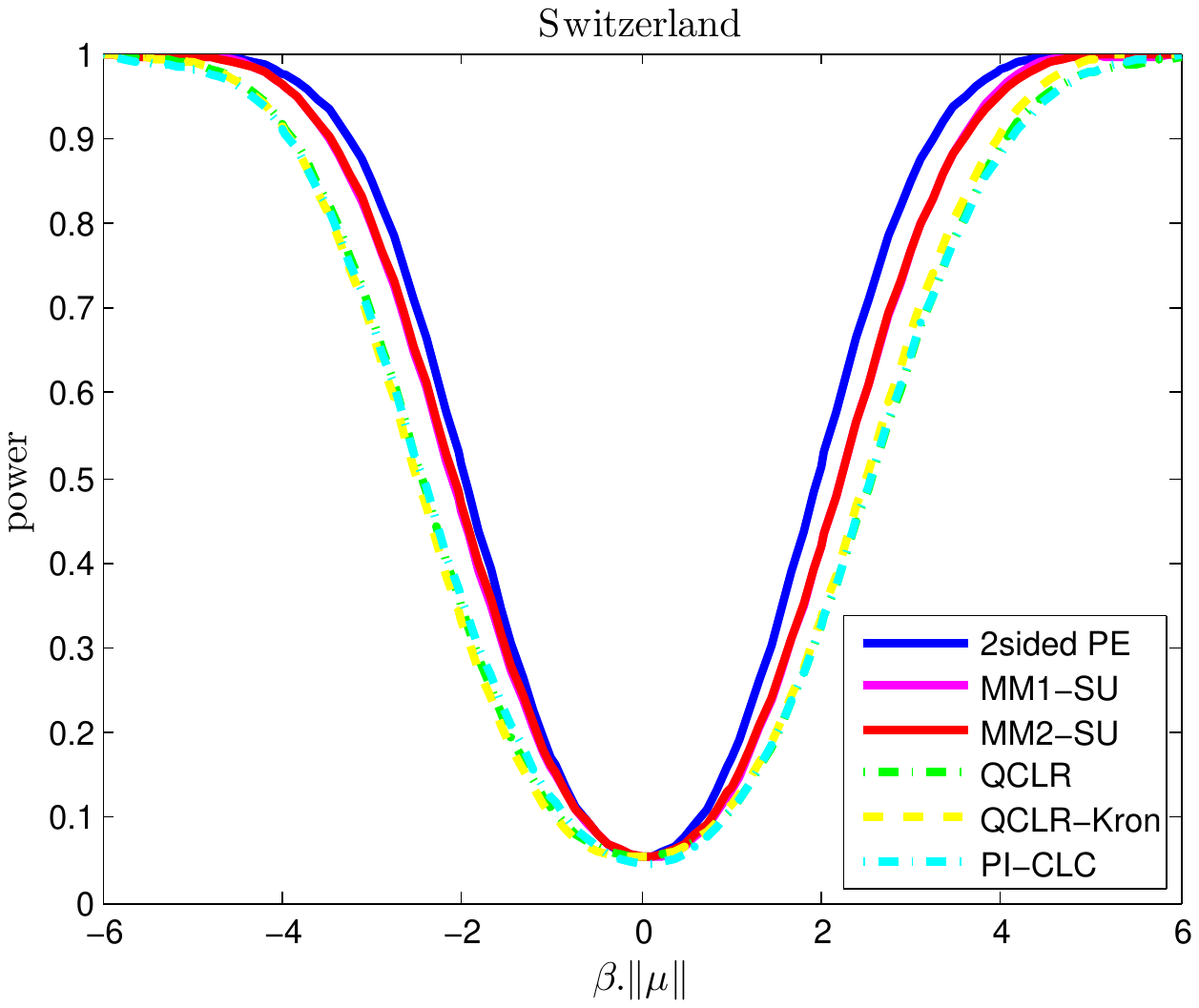}
\end{subfigure}%
\begin{subfigure}[b]{.5\textwidth}
  \centering
  \includegraphics[trim = 45mm 80mm 45mm 80mm, clip, width=5.5cm]{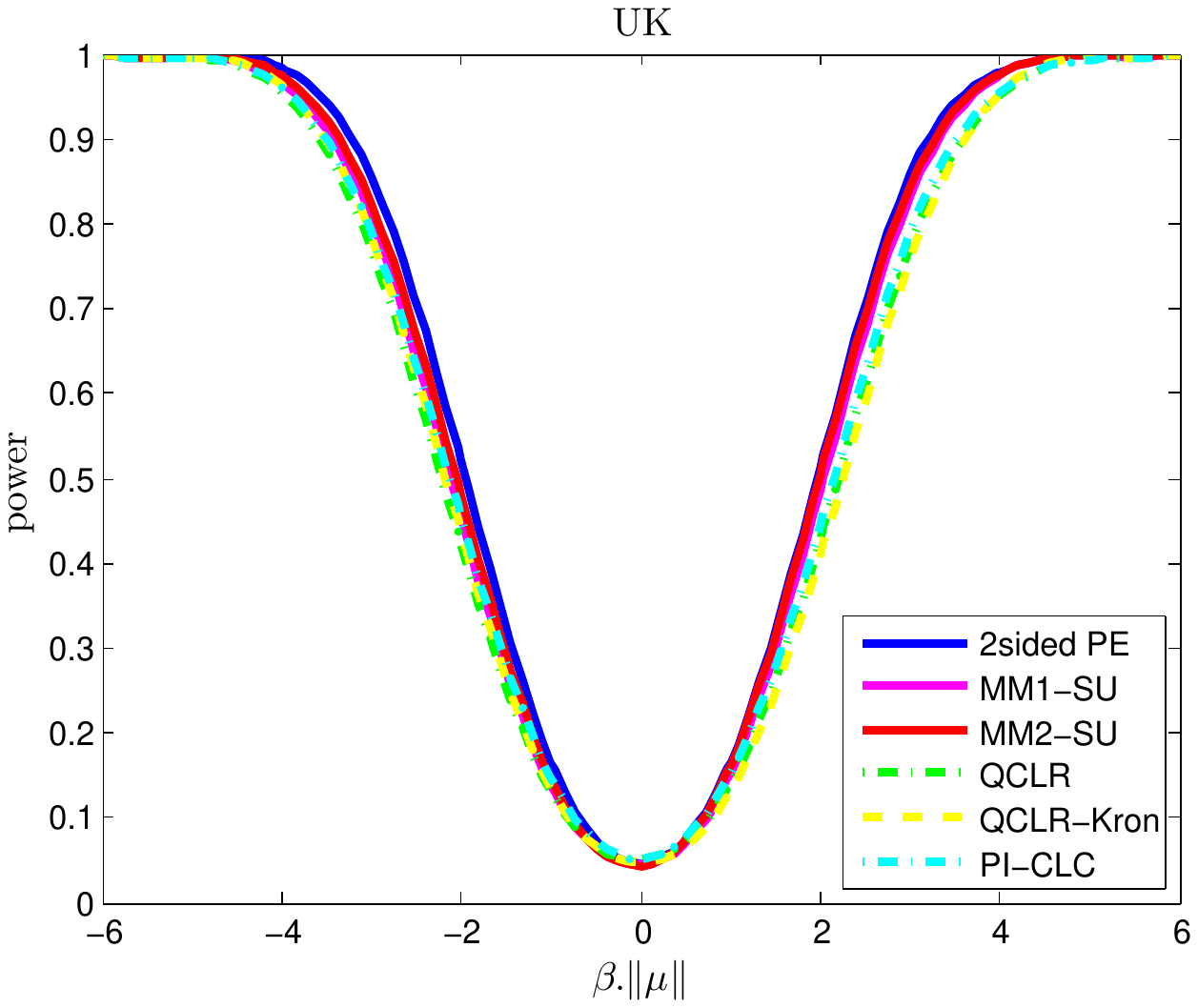}
\end{subfigure}
\begin{subfigure}[b]{.5\textwidth}
  \centering
  \includegraphics[trim = 45mm 80mm 45mm 80mm, clip, width=5.5cm]{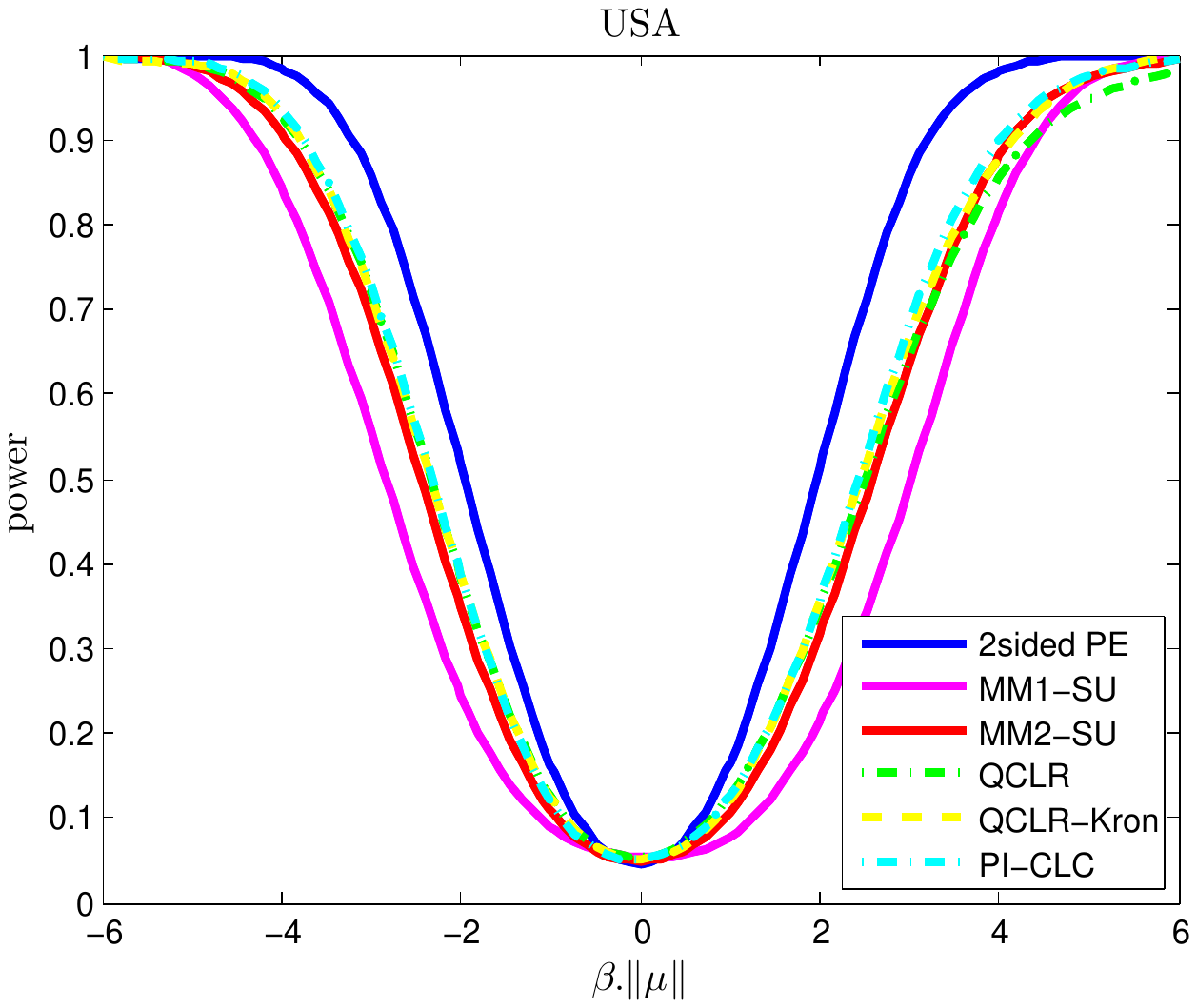}
\end{subfigure}
\end{figure}

\section{Concluding Remarks \label{Conclusion Sec}}

In this paper, we study the instrumental variable (IV) model with one
endogenous regressor and heteroskedastic and autocorrelated (HAC)\ errors.
The HAC-IV model with a known variance matrix is simpler than the model with
an unknown but consistently estimable long-run variance. However, inference
in both models is approximately the same whether or not the instruments are
weakly correlated with the endogenous variable. This simplification allows
us to develop a theory of optimal two-sided tests when the error stochastic
process is of unknown form.

We find that a test that has correct size and is optimal under standard
asymptotics may still have unacceptably low power in finite samples. This
issue appears in several econometric models. For the HAC-IV model, we solve
this problem by finding weighted-average power tests satisfying additional
two-sided conditions. In this paper, we consider two possibilities: the
locally unbiased (LU)\ and strongly unbiased (SU) conditions. While the
local condition yields admissible tests, the stronger condition is easier to
implement. Better yet, the MM1-SU and MM2-SU tests have power numerically
very close to their LU\ versions. Numerical simulations also show that the
MM2-SU test outperforms other tests proposed for the HAC-IV model.

The only other paper that satisfactorily addresses optimality of two-sided
tests in the HAC-IV model is that of I. \citet{Andrews15}. He explores
linear combinations of the Anderson-Rubin and score statistics, with weights
dependent on the conditioning statistic $T$. A class of these conditional
linear combination (CLC) tests is unbiased and admissible in the conditional
problem. By proposing a minimax regret criterion, he delivers a test which
plugs in a nuisance-parameter estimator. There is some power gained by
broadening the focus beyond those three statistics. On the other hand, we
impose $k$ additional constraints which are related to the SU condition. It
would be interesting to reduce the required computational time while
maintaining the power gains of the MM2-SU test by reducing the number of
boundary conditions when finding a WAP test.

Finally, the asymptotic theory based on Laplace approximations, developed in
this paper, is easily adaptable to other econometric models. For the HAC-IV
model, it relies on priors for the parameters $\beta $ and $\pi $ being
insensitive to the sample size. For the MM1 and MM2 weights, this implies
that the tuning parameters $\sigma ^{2}$ and $\zeta $ (used in the prior for 
$\mu =\left( Z^{\prime }Z\right) ^{1/2}\pi $) eventually grow at the sample
size $n$. Some power gains with weak instruments may be possible when the
tuning parameters are held constant. Another alternative is to find an
automatic rate for $\sigma ^{2}$ and $\zeta $ using a plug-in method. For
example, we could let these parameters be proportional to either $\left\Vert
T\right\Vert ^{2}$ or $n\cdot \left\Vert \pi \left( \beta _{0}\right)
\right\Vert ^{2}$. These quantities are stochastically bounded under weak
instruments and grow at the rate $n$ under strong instruments (which assures
asymptotic optimality). Since the constrained MLE $\pi \left( \beta
_{0}\right) $ is a one-to-one transformation of $T$, these modifications of
WAP-SU tests are still similar and uncorrelated with the pivotal statistic $%
S $ (hence, satisfy the SU\ Condition)\footnote{%
See \citeauthor{Moreira01} (\citeyear{Moreira01}, \citeyear{Moreira09a}) for
selecting among similar tests without creating size distortions; the
argument uses completeness of $T$ and is applicable to the SU condition as
well.}. We will consider this possibility in future work.

\bibliographystyle{econometrica}
\bibliography{References}

\end{document}